\renewcommand{\i}{\imath}
\newcommand{\et}{\hspace{-0.08in}{\bf .}\hspace{0.1in}}
\newcommand{\BOX}{\hbox {$\sqcap$ \kern -1em $\sqcup$}}
\newcommand{\qed}{\hskip 2em \hbox{\BOX} \vskip 2ex}
\renewcommand{\hom}{{\rm hom}}

\newcommand{\To}{\Rightarrow}
\newcommand{\Set}{{\rm Set}}
\newcommand{\Cat}{{\rm Cat}}
\newcommand{\Top}{{\rm Top}}
\newcommand{\Grp}{{\rm Grp}}
\newcommand{\CommRing}{{\rm CommRing}}
\newcommand{\Vect}{{\rm Vect}}
\newcommand{\Diff}{{\rm Diff}}
\renewcommand{\to}{\rightarrow}
\newcommand{\tensor}{\otimes}
\newcommand{\maps}{\colon}
\newcommand{\inv}{{\rm inv}}
\newcommand{\op}{{\rm op}}
\newcommand{\iso}{\cong}
\newcommand{\id}{{\rm id}}
\newcommand{\tr}{{\rm tr}}

\newcommand{\U}{{\rm U}}
\newcommand{\SU}{{\rm SU}}
\newcommand{\SO}{{\rm SO}}
\newcommand{\ISO}{{\rm ISO}}
\newcommand{\g}{{\mathfrak g}}
\newcommand{\h}{{\mathfrak h}}
\renewcommand{\u}{{\mathfrak u}}

\newcommand{\T}{{\cal T}}
\newcommand{\Mon}{{\rm MonCat}}
\newcommand{\Bitors}{{\rm Bitors}}

\newcommand{\Aut}{{\rm Aut}}

\newcommand{\R}{{\mathbb R}}

\newcommand{\Z}{{\mathbb Z}}

\newcommand{\cgl}{{\sf{c2g}}}
\newcommand{\wg}{{\rm W2G }}
\newcommand{\cg}{{\rm C2G }}
\newcommand{\sg}{{\rm S2G }}
\newcommand{\xb}{\bar{x}}
\newcommand{\yb}{\bar{y}}
\newcommand{\scs}{\scriptstyle}

\newcommand{\ten}{\otimes }
\newcommand{\imp}{{\rm Imp}}
\newcommand{\fex}{\widetilde{F(e_{x})}}
\newcommand{\fix}{\widetilde{F(i_{x})}}
\newcommand{\feix}{\widetilde{F(e_x)}^{-1}}

\newtheorem{thm}{Theorem}    
\newtheorem{cor}[thm]{Corollary}
\newtheorem{lem}[thm]{Lemma}

\newtheorem{prop}[thm]{Proposition}
\newtheorem{defn}[thm]{Definition}
\newtheorem{example}[thm]{Example}

\documentclass{article}
\usepackage{amsfonts,amssymb}
\usepackage[all]{xy}

\SelectTips{cm}{}

\title{Higher-Dimensional Algebra V: 2-Groups}
   \author{John C.\ Baez\\
   Department of Mathematics,  University of California\\
   Riverside, California 92521 \\
   USA \\
   \medskip
   email: baez@math.ucr.edu \\
   \medskip\\
   Aaron D.\ Lauda \\
   Department of Pure Mathematics and Mathematical Statistics \\
   University of Cambridge\\
   Cambridge CB3 0WB\\
   UK \\
   \medskip
   email: a.lauda@dpmms.cam.ac.uk \\}
    \date{October 1, 2004}
\begin{document}
\bibliographystyle{plain}
\maketitle

\begin{abstract}

\noindent A 2-group is a `categorified' version of a group, in
which the underlying set $G$ has been replaced by a category and
the multiplication map $m \maps G \times G \to G$ has been
replaced by a functor.  Various versions of this notion have
already been explored; our goal here is to provide a detailed
introduction to two, which we call `weak' and `coherent' 2-groups.
A weak 2-group is a weak monoidal category in which every morphism
has an inverse and every object $x$ has a `weak inverse': an
object $y$ such that $x \tensor y \iso 1 \iso y \tensor x$. A
coherent 2-group is a weak 2-group in which every object $x$ is
equipped with a {\it specified} weak inverse $\xb$ and
isomorphisms $i_x \maps 1 \to x \tensor \xb$, $e_x \maps \xb
\tensor x \to 1$ forming an adjunction.  We describe 2-categories of
weak and coherent 2-groups and an `improvement' 2-functor
that turns weak 2-groups into coherent ones, and prove that this
2-functor is a 2-equivalence of 2-categories.  We internalize the
concept of coherent 2-group, which gives a quick way to define Lie
2-groups.  We give a tour of examples, including the `fundamental
2-group' of a space and various Lie 2-groups.  We also explain how
coherent 2-groups can be classified in terms of 3rd cohomology
classes in group cohomology.  Finally, using this classification,
we construct for any connected and simply-connected compact simple Lie
group $G$ a family of 2-groups $G_\hbar$ ($\hbar \in \Z$) having $G$ as
its group of objects and $\U(1)$ as the group of automorphisms of
its identity object.  These 2-groups are built using Chern--Simons
theory, and are closely related to the Lie 2-algebras $\g_\hbar$
($\hbar \in \R$) described in a companion paper.
\end{abstract}

\section{Introduction} \label{introductionsection}

Group theory is a powerful tool in all branches of science where
symmetry plays a role.  However, thanks in large part to the
vision and persistence of Ronald Brown \cite{Brown2}, it has
become clear that group theory is just the tip of a larger subject
that deserves to be called `higher-dimensional group theory'. For
example, in many contexts where we are tempted to use groups, it
is actually more natural to use a richer sort of structure, where
in addition to group elements describing symmetries, we also have
isomorphisms between these, describing {\it symmetries between
symmetries}.  One might call this structure a `categorified'
group, since the underlying set $G$ of a traditional group is
replaced by a category, and the multiplication function $m \maps G
\times G \to G$ is replaced by a functor.   However, to hint at a
sequence of further generalizations where we use $n$-categories
and $n$-functors, we prefer the term `2-group'.

There are many different ways to make the notion of a 2-group
precise, so the history of this idea is complex, and we can
only briefly sketch it here.  A crucial first step
was J.\ H.\ C.\ Whitehead's \cite{Whitehead} concept of
`crossed module', formulated around 1946 without the aid of category theory.
In 1950, Mac Lane and Whitehead \cite{MW} proved that a crossed module
captures in algebraic form all the homotopy-invariant information about
what is now called a `connected pointed homotopy 2-type' --- roughly
speaking, a nice connected space equipped with a basepoint and
having homotopy groups that vanish above $\pi_2$.  By the 1960s it was
clear to Verdier and others that crossed modules are essentially the same
as `categorical groups'.   In the present paper we call these `strict
2-groups', since they are categorified groups in which the group laws hold
strictly, as {\it equations}.

Brown and Spencer \cite{BS} published a proof that crossed modules are
equivalent to categorical groups in 1976.  However, Grothendieck was
already familiar with these ideas, and in 1975 his student Hoang Xuan Sinh
wrote her thesis \cite{Sinh} on a more general concept, namely
`gr-categories', in which the group laws hold only {\it up to isomorphism}.
In the present paper we call these `weak' or `coherent' 2-groups, depending
on the precise formulation.

While influential, Sinh's thesis was never published, and is now quite
hard to find.  Also, while the precise relation between 2-groups,
crossed modules and group cohomology was greatly clarified in the 1986
draft of Joyal and Street's paper on braided tensor categories \cite{JS},
this section was omitted from the final published version.
So, while the basic facts about 2-groups are familiar to most experts in
category theory, it is difficult for beginners to find an introduction
to this material.  This is becoming a real nuisance as 2-groups
find their way into ever more branches of mathematics, and lately even
physics.  The first aim of the present paper is to fill this gap.

So, let us begin at the beginning.
Whenever one categorifies a mathematical concept, there are some
choices involved.  For example, one might define a 2-group simply
to be a category equipped with functors describing
multiplication, inverses and the identity, satisfying the usual
group axioms `on the nose' --- that is, as equations between
functors.  We call this a `strict' 2-group.  Part of the charm
of strict 2-groups is that they
can be defined in a large number of equivalent ways, including:
\begin{itemize}
 \item a strict monoidal category in which all objects
  and morphisms are invertible,
 \item a strict 2-category with one object in which all 1-morphisms and
 2-morphisms are invertible,
 \item a group object in $\Cat$ (also called a `categorical group'),
 \item a category object in $\Grp$,
 \item a crossed module.
\end{itemize}
There is an excellent review article by Forrester-Barker that
explains most of these notions and why they are equivalent
\cite{FB}.

Strict 2-groups have been applied in a variety of contexts, from
homotopy theory \cite{Brown,BS} and topological quantum field theory
\cite{Yetter} to nonabelian cohomology \cite{Breen,Breen2,Giraud}, the
theory of nonabelian gerbes \cite{Breen2,BM}, categorified gauge field
theory \cite{Attal,Baez,GP,Pfeiffer}, and even quantum gravity
\cite{CS,CY}.  However, the strict version of the 2-group concept is
not the best for all applications.  Rather than imposing the group axioms
as equational laws, it is sometimes better to `weaken' them: in other
words, to require only that they hold up to specified isomorphisms
satisfying certain laws of their own.  This leads to the concept of a
`coherent 2-group'.

For example, given objects $x,y,z$ in a strict 2-group we have
\[ (x \tensor y) \tensor z = x \tensor (y \tensor z)  \]
where we write multiplication as $\tensor$.   In a coherent
2-group, we instead specify an isomorphism called the
`associator':
\[
a_{x,y,z} \maps \xymatrix@1{(x \ten y) \ten z \;
\ar[r]^<<<<<{\sim} & \; x \ten (y \ten z)} .
\]
Similarly, we replace the left and right unit laws
\[        1 \tensor x = x, \qquad x \tensor 1 = x \]
by isomorphisms
\[
\ell_x \maps \xymatrix@1@=16pt{1 \tensor x \; \ar[r]^>>>>{\sim} &
\; x} , \qquad
  r_x \maps \xymatrix@1@=16pt{x \tensor 1\; \ar[r]^>>>>{\sim} &\; x}
\]
and replace the equations
\[ x \tensor x^{-1} = 1, \qquad x^{-1} \tensor x = 1 \]
by isomorphisms called the `unit' and `counit'.  Thus, instead of
an inverse in the strict sense, the object $x$ only has a
specified `weak inverse'.  To emphasize this fact, we denote this
weak inverse by $\xb$.

Next, to manipulate all these isomorphisms with some of the same
facility as equations, we require that they satisfy conditions
known as `coherence laws'.  The coherence laws for the associator
and the left and right unit laws were developed by Mac Lane
\cite{MacLane2} in his groundbreaking work on monoidal categories,
while those for the unit and counit are familiar from the
definition of an adjunction in a monoidal category \cite{JS}.
Putting these ideas together, one obtains Ulbrich and Laplaza's
definition of a `category with group structure'
\cite{Laplaza,Ulbrich}.  Finally, a `coherent 2-group' is a
category $G$ with group structure in which all morphisms are
invertible.  This last condition ensures that there is a covariant
functor
\[      \inv \maps G \to G   \]
sending each object $x \in G$ to its weak inverse $\xb$; otherwise
there will only be a contravariant functor of this sort.

In this paper we compare this sort of 2-group to a simpler sort,
which we call a `weak 2-group'.  This is a weak monoidal category
in which every morphism has an inverse and every object $x$ has a
weak inverse: an object $y$ such that $y \tensor x \iso 1$ and $x
\tensor y \iso 1$. Note that in this definition, we do not {\it
specify} the weak inverse $y$ or the isomorphisms from $y \tensor
x$ and $x \tensor y$ to $1$, nor do we impose any coherence laws
upon them.  Instead, we merely demand that they {\it exist}.
Nonetheless, it turns out that any weak 2-group can be improved to
become a coherent one!  While this follows from a theorem of
Laplaza \cite{Laplaza}, it seems worthwhile to give an expository
account here, and to formalize this process as a 2-functor
\[ \imp \maps \wg \to \cg \]
where $\wg$ and $\cg$ are suitable strict 2-categories of weak and
coherent 2-groups, respectively.

On the other hand, there is also a forgetful 2-functor
\[  {\rm F} \maps \cg \to \wg . \]
One of the goals of this paper is to show that $\imp$ and ${\rm
F}$ fit together to define a 2-equivalence of strict 2-categories.
This means that the 2-category of weak 2-groups and the 2-category
of coherent 2-groups are `the same' in a suitably weakened sense.
Thus there is ultimately not much difference between weak and
coherent 2-groups.

To show this, we start in Section \ref{weaksection} by defining
weak 2-groups and the 2-category $\wg$.  In Section
\ref{coherentsection} we define coherent 2-groups and the
2-category $\cg$.  To do calculations in 2-groups, it turns out
that certain 2-dimensional pictures called `string diagrams' can
be helpful, so we explain these in Section \ref{stringsection}.
In Section \ref{improvementsection} we use string diagrams to
define the `improvement' 2-functor $\imp \maps \wg \to \cg$ and
prove that it extends to a 2-equivalence of strict 2-categories.
This result relies crucially on the fact that morphisms in $\cg$
are just weak monoidal functors, with no requirement that they
preserve weak inverses.   In Section \ref{preservationsection} we
justify this choice, which may at first seem questionable, by
showing that weak monoidal functors automatically preserve the
specified weak inverses, up to a well-behaved isomorphism.

In applications of 2-groups to geometry and physics, we expect the
concept of {\it Lie 2-group} to be particularly important. This is
essentially just a 2-group where the set of objects and the set of
morphisms are manifolds, and all relevant maps are smooth.  Until
now, only strict Lie 2-groups have been defined \cite{Baez}.  In
section \ref{internalizationsection} we show that the concept of
`coherent 2-group' can be defined in any 2-category with finite
products.  This allows us to efficiently define coherent Lie
2-groups, topological 2-groups and the like.

In Section \ref{examplesection} we discuss examples of 2-groups.
These include various sorts of `automorphism 2-group' for an object
in a 2-category, the `fundamental 2-group' of a topological space, and
a variety of strict Lie 2-groups.  We also describe a way to classify
2-groups using group cohomology.  As we explain, coherent 2-groups ---
and thus also weak 2-groups --- can be classified up to equivalence in
terms of a group $G$, an action $\alpha$ of $G$ on an abelian group $H$,
and an element $[a]$ of the 3rd cohomology group of $G$ with coefficients
in $H$.  Here $G$ is the group of objects in a `skeletal' version of
the 2-group in question: that is, an equivalent 2-group containing
just one representative from each isomorphism class of objects.
$H$ is the group of automorphisms of the identity object, the action
$\alpha$ is defined using conjugation, and the 3-cocycle $a$ comes
from the associator in the skeletal version.  Thus, $[a]$ can be thought
of as the obstruction to making the 2-group simultaneously both skeletal
and strict.

In a companion to this paper, called HDA6 \cite{HDA6} for short,
Baez and Crans prove a Lie algebra analogue of this result: a
classification of `semistrict Lie 2-algebras'.  These are categorified
Lie algebras in which the antisymmetry of the Lie bracket holds on the
nose, but the Jacobi identity holds only up to a natural
isomorphism called the `Jacobiator'.  It turns out that semistrict Lie
2-algebras are classified up to equivalence
by a Lie algebra $\g$, a representation $\rho$ of $\g$ on an abelian Lie
algebra $\h$, and an element $[j]$ of the 3rd Lie algebra cohomology group
of $\g$ with coefficients in $\h$.   Here the cohomology class $[j]$
comes from the Jacobiator in a skeletal version of the Lie 2-algebra in
question.  A semistrict Lie 2-algebra in which the Jacobiator is the
identity is called `strict'.  Thus, the class $[j]$ is the obstruction
to making a Lie 2-algebra simultaneously skeletal and strict.

Interesting examples of Lie 2-algebras that cannot be made both skeletal
and strict arise when $\g$ is a finite-dimensional simple Lie algebra
over the real numbers.  In this
case we may assume without essential loss of generality that $\rho$ is
irreducible, since any representation is a direct sum of irreducibles.
When $\rho$ is irreducible, it turns out that $H^3(\g,\rho) = \{0\}$
unless $\rho$ is the trivial representation on the 1-dimensional abelian
Lie algebra $\u(1)$, in which case we have
\[         H^3(\g,\u(1)) \iso \R . \]
This implies that for any value of $\hbar \in \R$ we obtain a skeletal
Lie 2-algebra $\g_\hbar$ with $\g$ as its Lie algebra of objects,
$\u(1)$ as the endomorphisms of its zero object, and $[j]$ proportional
to $\hbar \in \R$.   When $\hbar = 0$, this Lie 2-algebra is just
$\g$ with identity morphisms adjoined to make it into a strict
Lie 2-algebra.  But when $\hbar \ne 0$, this Lie 2-algebra is
not equivalent to a skeletal strict one.

In short, the Lie algebra $\g$ sits inside a one-parameter family of
skeletal Lie 2-algebras $\g_\hbar$, which are strict only for $\hbar = 0$.
This is strongly reminiscent of some other well-known deformation
phenomena arising from the third cohomology of a simple Lie algebra.
For example, the universal enveloping algebra of $\g$ gives a
one-parameter family of quasitriangular Hopf algebras $U_\hbar \g$,
called `quantum groups'.  These Hopf algebras are cocommutative
only for $\hbar = 0$.  The theory of `affine Lie algebras' is based on
a closely related phenomenon: the Lie algebra of smooth functions
$C^\infty(S^1,\g)$ has a one-parameter family of central extensions,
which only split for $\hbar = 0$.  There is also a group version of
this phenomenon, which involves an integrality condition: the loop group
$C^\infty(S^1,G)$ has a one-parameter family of central extensions,
one for each $\hbar \in \Z$.  Again, these central extensions split
only for $\hbar = 0$.

All these other phenomena are closely connected to Chern--Simons
theory, a topological quantum field theory whose action is the secondary
characteristic class associated to an element of $H^4(BG,\Z) \iso \Z$.
The relation to Lie algebra cohomology comes from the existence
of an inclusion $H^4(BG,\Z) \hookrightarrow H^3(\g,\u(1)) \iso \R$.

Given all this, it is tempting to seek a 2-group analogue of the
Lie 2-algebras $\g_\hbar$.  Indeed, such an analogue exists!  Suppose
that $G$ is a connected and simply-connected compact simple Lie group.
In Section \ref{chernsimonssection} we construct a family of skeletal
2-groups $G_\hbar$, one for each $\hbar \in \Z$, each having $G$ as its
group of objects and $\U(1)$ as the group of automorphisms of its identity
object.  The associator in these 2-groups depends on $\hbar$, and
they are strict only for $\hbar = 0$.

Unfortunately, for reasons we shall explain, these 2-groups are not Lie
2-groups except for the trivial case $\hbar = 0$.   However, the
construction of these 2-groups uses Chern--Simons theory in an essential
way, so we feel confident that they are related to all the other
deformation phenomena listed above.  Since the rest of these phenomena
are important in mathematical physics, we hope these 2-groups $G_\hbar$
will be relevant as well.  A full understanding of them may require a
generalization of the concept of Lie 2-group presented in this paper.

{\bf Note}: in all that follows, we write the composite of
morphisms $f \maps x \to y$ and $g \maps y \to z$ as $fg \maps x
\to z$. We use the term `weak 2-category' to refer to a
`bicategory' in B\'enabou's sense \cite{Benabou}, and the term
`strict 2-category' to refer to what is often called simply a
`2-category' \cite{Street2}.

\section{Weak 2-groups} \label{weaksection}

Before we define a weak 2-group, recall that a {\bf weak monoidal
category} consists of:
\begin{description}
 \item[(i)] a category $M$,
 \item[(ii)] a functor $m \maps M \times M \to M$, where we
 write $m(x,y)=x \ten y$ and $m(f,g)=f \ten g$ for objects $x, y,
 \in M$ and morphisms $f, g$ in $M$,
 \item[(iii)] an `identity object' $1 \in M$,
 \item[(iv)] natural isomorphisms
  \[ a_{x,y,z} \maps (x \ten y) \ten z \to x \ten (y \ten z), \]
  \[ \ell_x \maps 1 \ten x \to x , \]
  \[ r_x \maps x \ten 1 \to x,  \]
\end{description}
such that the following diagrams commute for all objects $w,x,y,z
\in M$:
\[
\xy
 (0,20)*+{(w \ten x) \ten (y \ten z)}="1";
 (40,0)*+{w \ten (x \ten (y \ten z))}="2";
 (25,-20)*+{ \quad w \ten ((x \ten y) \ten z)}="3";
 (-25,-20)*+{(w \ten (x \ten y)) \ten z}="4";
 (-40,0)*+{((w \ten x) \ten y) \ten z}="5";
     {\ar^{a_{w,x,y \ten z}}     "1";"2"}
     {\ar_{1_w \ten a _{x,y,z}}  "3";"2"}
     {\ar^{a _{w,x \ten y,z}}    "4";"3"}
     {\ar_{a _{w,x,y} \ten 1_z}  "5";"4"}
     {\ar^{a _{w \ten x,y,z}}    "5";"1"}
\endxy
\\
\]
\vskip 1em
\[
\xymatrix{
 (x \ten 1) \ten y
     \ar[rr]^{a_{x,1,y}}
     \ar[dr]_{r_x \ten 1_y}
  &&  x \ten (1 \ten y)
     \ar[dl]^{1_x \ten \ell_y } \\
 & x \ten y   }  \\
\]
A {\bf strict monoidal category} is the special case where
$a_{x,y,z},\ell_x,r_x$ are all identity morphisms.  In this case
we have
\[ (x \tensor y) \tensor z = x \tensor (y \tensor z) , \]
\[        1 \tensor x = x, \qquad x \tensor 1 = x .\]
As mentioned in the Introduction, a {\bf strict 2-group} is a
strict monoidal category where every morphism is invertible and
every object $x$ has an inverse $x^{-1}$, meaning that
\[      x \tensor x^{-1} = 1, \qquad x^{-1} \tensor x = 1  .\]

Following the principle that it is wrong to impose equations
between objects in a category, we can instead start with a weak
monoidal category and require that every object has a `weak'
inverse.  With these changes we obtain the definition of `weak
2-group':

\begin{defn} \label{weakinv} \et
If $x$ is an object in a weak monoidal category, a {\bf weak
inverse} for $x$ is an object $y$ such that $x \tensor y \iso 1$
and $y \tensor x \iso 1$.  If $x$ has a weak inverse, we call it
{\bf weakly invertible}.
\end{defn}

\begin{defn} \label{weak} \et
A {\bf weak 2-group} is a weak monoidal category where all objects
are weakly invertible and all morphisms are invertible.
\end{defn}

\noindent  In fact, Joyal and Street \cite{JS} point out that
when every object in a weak monoidal category has a `one-sided'
weak inverse, every object is weakly invertible in the above sense.
Suppose for example that every object $x$ has an object $y$
with $y \tensor x \iso 1$.  Then $y$ has an object $z$ with
$z \tensor y \iso 1$, and
\[        z \iso z \tensor 1 \iso z \tensor (y \tensor x)
          \iso (z \tensor y) \tensor x \iso 1 \tensor x \iso x ,\]
so we also have $x \tensor y \iso 1$.

Weak 2-groups are the objects of a strict 2-category $\wg$; now
let us describe the morphisms and 2-morphisms in this 2-category.
Notice that the only {\it structure} in a weak 2-group is that of
its underlying weak monoidal category; the invertibility
conditions on objects and morphisms are only {\it properties}.
With this in mind, it is natural to define a morphism between weak
2-groups to be a weak monoidal functor.  Recall that a {\bf weak
monoidal functor} $F \maps C \to C'$ between monoidal categories
$C$ and $C'$ consists of:
\begin{description}
\item[(i)] a functor $F \maps C \to C'$, \item[(ii)] a natural
isomorphism $F_{2} \maps F(x) \otimes F(y) \to F(x \otimes y)$,
where for brevity we suppress the subscripts indicating the
dependence of this isomorphism on $x$ and $y$, \item[(iii)] an
isomorphism $F_{0} \maps 1' \to F(1)$, where $1$ is the identity
object of $C$ and $1'$ is the identity object of $C'$,
\end{description}
such that the following diagrams commute for all objects $x,y,z
\in C$:
\[
\xymatrix@!C{
 (F(x) \ten F(y)) \ten F(z)
   \ar[r]^>>>>>>>{F_{2} \ten 1}
   \ar[d]_{a_{F(x), F(y), F(z)}}
& F(x \ten y) \ten F(z)
   \ar[r]^{F_{2}}
& F((x \ten y) \ten z)
   \ar[d]^{F(a_{x,y,z})}   \\
 F(x) \ten (F(y) \ten F(z))
   \ar[r]^>>>>>>>{1 \ten F_{2}}
& F(x) \ten F(y \ten z)
   \ar[r]^{F_{2}}
& F(x \ten (y \ten z))
 }
\]
\vskip 1em
\[
\xymatrix{
 1' \ten F(x)
   \ar[r]^{\ell'_{F(x)}}
   \ar[d]_{F_{0} \ten 1}
&  F(x)  \\
 F(1) \ten F(x)
   \ar[r]^{F_{2}}
&  F(1 \ten x)
   \ar[u]_{F(\ell_{x})}
}
\]
\vskip 1em
\[
\xymatrix{
 F(x) \ten 1'
   \ar[r]^{r'_{F(x)}}
   \ar[d]_{1 \ten F_{0}}
&  F(x)    \\
 F(x) \ten F(1)
   \ar[r]^{F_{2}}
&  F(x \ten 1)
   \ar[u]_{F(r_x)}
}
\]

A weak monoidal functor preserves tensor products and the identity
object up to specified isomorphism.   As a consequence, it also
preserves weak inverses:

\begin{prop} \et
If $F \maps C \to C'$ is a weak monoidal functor and $y \in C$ is
a weak inverse of $x \in C$, then $F(y)$ is a weak inverse of
$F(x)$ in $C'$.
\end{prop}

{\bf Proof. } Since $y$ is a weak inverse of $x$, there exist
isomorphisms $\gamma \maps x \ten y \to 1$ and $\xi \maps y \ten x
\to 1$. The proposition is then established by composing the
following isomorphisms:
\[
\xymatrix{
 F(y) \ten F(x)
   \ar[r]^<<<<<{\sim}
   \ar[d]_{F_{2}}
&  1'       \\
 F(y \ten x)
   \ar[r]_{F(\xi)}
&  F(1)
   \ar[u]_{F_{0}^{-1}}
} \qquad \qquad \xymatrix{
 F(x) \ten F(y)
   \ar[r]^<<<<<{\sim}
   \ar[d]_{F_{2}}
&  1'       \\
 F(x \ten y)
   \ar[r]_{F(\gamma)}
&  F(1)
   \ar[u]_{F_{0}^{-1}}
}
\]
\hskip 30em \qed

We thus make the following definition:

\begin{defn} \et
A {\bf homomorphism} $F \maps C \to C'$ between weak 2-groups is a
weak monoidal functor.
\end{defn}
\noindent The composite of weak monoidal functors is again a weak
monoidal functor \cite{moncat}, and composition satisfies
associativity and the unit laws.  Thus, 2-groups and the
homomorphisms between them form a category.

Although they are not familiar from traditional group theory, it
is natural in this categorified context to also consider
`2-homomorphisms' between homomorphisms.  Since a homomorphism
between weak 2-groups is just a weak monoidal functor, it makes
sense to define 2-homomorphisms to be monoidal natural
transformations.  Recall that if $F,G \maps C \to C'$ are weak
monoidal functors, then a {\bf monoidal natural transformation}
$\theta \maps F \To G$ is a natural transformation
such that the following diagrams commute for all $x,y \in C$.
\[
\xymatrix{
 F(x) \ten F(y)
  \ar[rr]^{\theta_{x} \ten \theta_{y}}
  \ar[d]_{F_{2}}
&&  G(x) \ten G(y)
   \ar[d]^{G_{2}}     \\
 F(x \ten y)
   \ar[rr]^{\theta_{x \ten y}}
&&  G(x \ten y) }
\]

\[
\\
\\
\xymatrix{
 1'
   \ar[d]_{F_{0}}
   \ar[dr]^{G_{0}}  \\
 F(1)
   \ar[r]^{\theta_{1}}
&  G(1) }
\]
Thus we make the following definitions:

\begin{defn} \et
A {\bf 2-homomorphism} $\theta \maps F \To G$ between
homomorphisms \hfill \break $F,G \maps C \to C'$ of weak 2-groups
is a monoidal natural transformation.
\end{defn}

\begin{defn} \et Let {\bf W2G} be the strict 2-category consisting of
weak 2-groups, \break homomorphisms between these, and 2-homomorphisms
between those.  \end{defn}

\noindent There is a strict 2-category $\Mon$ with weak monoidal
categories as objects, weak monoidal functors as 1-morphisms, and
monoidal natural transformations as 2-morphisms \cite{moncat}.
$\wg$ is a strict 2-category because it is a sub-2-category of $\Mon$.

\section{Coherent 2-groups} \label{coherentsection}

In this section we explore another notion of 2-group. Rather than
requiring that objects be weakly invertible, we will require that
every object be equipped with a specified adjunction.  Recall that
an {\bf adjunction} is a quadruple $(x, \xb, i_x, e_x)$ where $i_x
\maps 1 \to x \ten \xb$ (called the {\bf unit}) and $e_x \maps \xb
\ten x \to 1$ (called the {\bf counit}) are morphisms such that
the following diagrams commute:
\[
\xymatrix@C=2.3pc{
  1 \ten x \ar[r]^<<<<<<{i_x \ten 1} \ar[d]_{\ell_x}
  & (x \ten \xb) \ten x \ar[r]^<<<<<{a_{x, \bar{x},x}}
  & x\ten( \bar{x}\ten x) \ar[d]^{1 \ten e_x} \\
  x \ar[rr]_{r^{-1}_x}
   && x \ten 1  }
\]

\[
\xymatrix@C2.3pc{
 \xb \ten 1
     \ar[r]^<<<<<<{1 \ten i_x}
     \ar[d]_{r_{\xb}}
   & \xb \ten (x \ten \xb)
     \ar[r]^{a^{-1}_{\xb, x, \xb}}
   & (\xb \ten x)\ten \xb
     \ar[d]^{e_x \ten 1} \\
 \xb
     \ar[rr]_{\ell^{-1}_{\xb}}
   && 1 \ten \xb   }
\]
When we express these laws using string diagrams in Section
\ref{stringsection}, we shall see that they give ways to
`straighten a zig-zag' in a piece of string.  Thus, we refer to
them as the first and second {\bf zig-zag identities},
respectively.

An adjunction $(x,\xb,i_x,e_x)$ for which the unit and counit are
invertible is called an {\bf adjoint equivalence}.  In this case
$x$ and $\xb$ are weak inverses.  Thus, specifying an adjoint
equivalence for $x$ ensures that $\xb$ is weakly invertible ---
but it does so by providing $x$ with extra {\it structure}, rather
than merely asserting a {\it property} of $x$.  We now make the
following definition:

\begin{defn}  \label{coherent} \et
A { \bf coherent 2-group } is a weak monoidal category $C$ in
which every morphism is invertible and every object $x \in C$ is
equipped with an adjoint equivalence $(x,\bar{x},i_{x},e_{x})$.
\end{defn}

\noindent Coherent 2-groups have been studied under many names.
Sinh \cite{Sinh} called them `gr-categories' when she initiated
work on them in 1975, and this name is also used by Saavedra
Rivano \cite{Rivano} and Breen \cite{Breen2}. As noted in the
Introduction, a coherent 2-group is the same as one of Ulbrich
and Laplaza's `categories with group structure'
\cite{Laplaza,Ulbrich} in which all morphisms are invertible.  It
is also the same as an `autonomous monoidal category' \cite{JS}
with all morphisms invertible, or a `bigroupoid' \cite{HKK} with
one object.

As we did with weak 2-groups, we can define a homomorphism between
coherent 2-groups.  As in the weak 2-group case we can begin by
taking it to be a weak monoidal functor, but now we must consider
what additional structure this must have to preserve each adjoint
equivalence $(x,\xb,i_x,e_x)$, at least up to a specified
isomorphism. At first it may seem that an additional structural
map is required. That is, given a weak monoidal functor $F$
between 2-groups, it may seem that we must include a natural
isomorphism
\[ F_{-1} \maps \overline{F(x)} \to F(\xb) \]
relating the weak inverse of the image of $x$ to the image of the
weak inverse $\xb$.  In Section \ref{preservationsection} we shall
show this is not the case: $F_{-1}$ can be constructed from the
data already present!  Moreover, it automatically satisfies the
appropriate coherence laws.  Thus we make the following
definitions:

\begin{defn} \label{coherenthomo} \et
A {\bf homomorphism} $F \maps C \to C'$ between coherent 2-groups
is a weak monoidal functor.
\end{defn}

\begin{defn} \label{coherent2homo} \et
A {\bf 2-homomorphism} $\theta \maps F \To G$ between
homomorphisms \break $F,G \maps C \to C'$ of coherent 2-groups is
a monoidal natural transformation.
\end{defn}

\begin{defn} \et
Let {\bf C2G} be the strict 2-category consisting of coherent
2-groups, homomorphisms between these, and 2-homomorphisms between
those.  \end{defn}

\noindent It is clear that $\cg$ forms a strict 2-category since
it is a sub-2-category of $\Mon$.

We conclude this section by stating the theorem that justifies the
term `coherent 2-group'.  This result is analogous to Mac Lane's
coherence theorem for monoidal categories.  A version of this
result was proved by Ulbrich \cite{Ulbrich} and Laplaza
\cite{Laplaza} for a structure called a {\bf category with group
structure}: a weak monoidal category equipped with an adjoint
equivalence for every object.  Through a series of lemmas, Laplaza
establishes that there can be at most one morphism between any two
objects in the free category with group structure on a set of
objects.  Here we translate this into the language of 2-groups and
explain the significance of this result.

Let $\cgl$ be the category of coherent 2-groups where the
morphisms are the functors that \textit{strictly} preserve the
monoidal structure and specified adjoint equivalences for each
object. Clearly there exists a forgetful functor $U \maps \cgl \to
\Set$ sending any coherent 2-group to its underlying set.  The
interesting part is:

\begin{prop} \et
The functor $U \maps \cgl \to \Set$ has a left adjoint $F \maps
\Set \to \cgl$.
\end{prop}

Since $a, \ell, r, i$ and $e$ are all isomorphism, the free
category with group structure on a set $S$ is the same as the free
coherent 2-group on $S$, so Laplaza's construction of $F(S)$
provides most of what we need for the proof of this theorem. In
Laplaza's words, the construction of $F(S)$ for a set $S$ is
``long, straightforward, and rather deceptive", because it hides
the essential simplicity of the ideas involved. For this reason,
we omit the proof of this theorem and refer the interested reader
to Laplaza's paper.

It follows that for any coherent 2-group $G$ there exists a
homomorphism of 2-groups $e_G \maps F(U(G)) \to G$ that
\textit{strictly} preserves the monoidal structure and chosen
adjoint equivalences.  This map allows us to interpret formal
expressions in the free coherent 2-group $F(U(G))$ as actual
objects and morphisms in $G$.  We now state the coherence theorem:

\begin{thm} \et
There exists at most one morphism between any pair of objects in
$F(U(G))$.
\end{thm}

This theorem, together with the homomorphism $e_G$, makes precise
the rough idea that there is at most one way to build an
isomorphism between two tensor products of objects and their weak
inverses in $G$ using $a, \ell, r, i$, and $e$.

\section{String diagrams} \label{stringsection}

Just as calculations in group theory are often done using
1-dimensional symbolic expressions such as
\[ x(yz)x^{-1} = (xyx^{-1})(xzx^{-1})  , \]
calculations in 2-groups are often done using 2-dimensional
pictures called string diagrams.  This is one of the reasons for
the term `higher-dimensional algebra'.  String diagrams for
2-categories \cite{Street} are Poincar\'e dual to the more
traditional globular diagrams in which objects are represented as
dots, 1-morphisms as arrows and 2-morphisms as 2-dimensional
globes.  In other words, in a string diagram one draws objects in
a 2-category as 2-dimensional regions in the plane, 1-morphisms as
1-dimensional `strings' separating regions, and 2-morphisms as
0-dimensional points (or small discs, if we wish to label them).

To apply these diagrams to 2-groups, first let us assume our
2-group is a strict monoidal category, which we may think of as a
strict 2-category with a single object, say $\bullet$.  A morphism
$f \maps x \to y $ in the monoidal category corresponds to a
2-morphism in the 2-category, and we convert the globular picture
of this into a string diagram as follows:
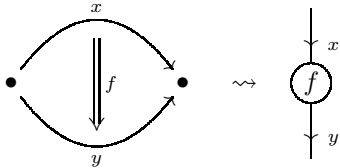
\[
\xymatrix{
 \bullet
   \ar@/^2pc/[rr]_{\quad}^{x}="1"
   \ar@/_2pc/[rr]_{y}="2"
&& \bullet
 \ar@{}"1";"2"|(.2){\,}="7"
  \ar@{}"1";"2"|(.8){\,}="8"
 \ar@{=>}"7" ;"8"^{f}
} \quad \rightsquigarrow \quad \xy (0,0)*{f};
(0,0)*\xycircle(2.65,2.65){-}="f"; (0,10)**\dir{-}
?(.5)*\dir{<}+(3,0)*{\scriptstyle x}; "f";(0,-10)**\dir{-}
?(.75)*\dir{>}+(3,0)*{\scriptstyle y};
\endxy
\]

We can use this idea to draw the composite or tensor product of
morphisms.   Composition of morphisms $f \maps x \to y$ and $g
\maps y \to z$ in the strict monoidal category corresponds to
vertical composition of 2-morphisms in the strict 2-category with
one object.   The globular picture of this is:
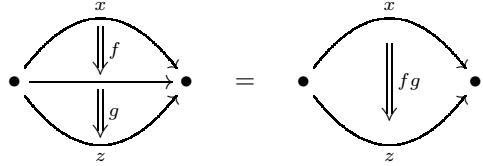
\[
\xymatrix{
 \bullet
   \ar@/^2pc/[rr]^{x}_{}="0"
   \ar[rr]^{}="1"
   \ar@/_2pc/[rr]_{z}_{}="2"
   \ar@{=>}"0";"1" ^{f}
   \ar@{=>}"1";"2" ^{g}
&&  \bullet } \quad = \quad \xymatrix{
 \bullet
   \ar@/^2pc/[rr]^{x}_{}="0"
   \ar@/_2pc/[rr]_{z}_{}="2"
   \ar@{}"0";"2"|(.2){\,}="7"
  \ar@{}"0";"2"|(.8){\,}="8"
   \ar@{=>}"7";"8" ^{fg}
&&  \bullet }
\]
and the Poincar\'e dual string diagram is:
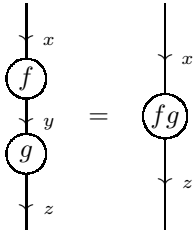
\[
\xy (0,5)*{f}; (0,5)*\xycircle(2.65,2.65){-}="f"; (0,15)**\dir{-}
?(.5)*\dir{<}+(3,0)*{\scriptstyle x};
"f";(0,-5)*\xycircle(2.65,2.65){-}="g"; **\dir{-}
?(.4)*\dir{<}+(3,0)*{\scriptstyle y}; "f";"g";(0,-15)**\dir{-}
?(.75)*\dir{>}+(3,0)*{\scriptstyle z}; (0,-5)*{g}
\endxy
\quad = \quad \xy  (0,0)*\xycircle(2.95,2.95){-}="fg"; (0,15)*{}
**\dir{-} ?(.5)*\dir{<}+(3,0)*{\scriptstyle x}; "fg";(0,-15)*{};
**\dir{-} ?(.4)*\dir{<}+(3,0)*{\scriptstyle z}; (0,0)*{fg};
\endxy
\]
Similarly, the tensor product of morphisms $f \maps x \to y$ and
$g \maps x' \to y'$ corresponds to horizontal composition of
2-morphisms in the 2-category. The globular picture is:
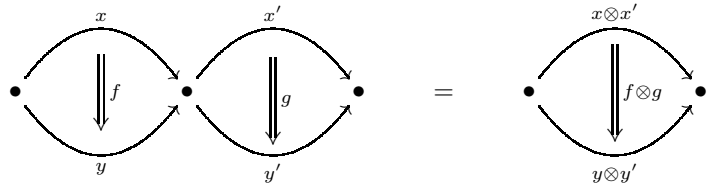
\[
\xymatrix{
 \bullet
   \ar@/^2pc/[rr]^{x}_{}="0"
   \ar@/_2pc/[rr]_{y}_{}="2"
   \ar@{}"0";"2"|(.2){\,}="7"
  \ar@{}"0";"2"|(.8){\,}="8"
   \ar@{=>}"7";"8" ^{f}          
&&  \bullet
   \ar@/^2pc/[rr]^{x'}_{}="0"
   \ar@/_2pc/[rr]_{y'}="2"
   \ar@{}"0";"2"|(.2){\,}="7"
  \ar@{}"0";"2"|(.8){\,}="8"
   \ar@{=>}"7";"8" ^{g}
&&  \bullet } \qquad = \qquad \xymatrix{
 \bullet
   \ar@/^2pc/[rr]_{\quad}^{x \ten x'}="1"
   \ar@/_2pc/[rr]_{y \ten y'}="2"
&& \bullet
 \ar@{}"1";"2"|(.2){\,}="7"
  \ar@{}"1";"2"|(.8){\,}="8"
 \ar@{=>}"7" ;"8"^{f \ten g}
}
\]
and the Poincar\'e dual string diagram is:
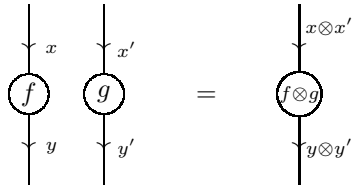
\[
\xy
 (-5,0)*{f};
(-5,12)*{}; (-5,0)*\xycircle(2.65,2.65){-}="1_x"; **\dir{-}
?(.5)*\dir{<}+(3,0)*{\scriptstyle x}; "1_x";(-5,-12)*{}; **\dir{-}
?(.4)*\dir{<}+(3,0)*{\scriptstyle y}; (5,12)*{};
(5,0)*\xycircle(2.65,2.65){-}="1_x"; **\dir{-}
?(.5)*\dir{<}+(3,0)*{\scriptstyle x'}; "1_x";(5,-12)*{}; **\dir{-}
?(.4)*\dir{<}+(3,0)*{\scriptstyle y'}; (5,0)*{g};
\endxy
\qquad = \qquad \xy (0,12)*{};
(0,0)*\xycircle(2.95,2.95){-}="1_x"; **\dir{-}
?(.5)*\dir{<}+(4,3)*{\scriptstyle x \ten x'}; "1_x";(0,-12)*{};
**\dir{-} ?(.4)*\dir{<}+(4,0)*{\scriptstyle y \ten y'}; (0,0)*{
\scriptstyle f \ten g};
\endxy
\]

We also introduce abbreviations for identity morphisms and the
identity object.  We draw the identity morphism $1_x \maps x \to
x$ as a straight vertical line:
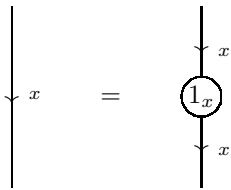
\[
\xy (0,12)*{};(0,-12)*{}; **\dir{-}
?(.47)*\dir{<}+(3,1)*{\scriptstyle x}
\endxy
\qquad = \qquad \xy (0,12)*{};
(0,0)*\xycircle(2.65,2.65){-}="1_x"; **\dir{-}
?(.5)*\dir{<}+(3,0)*{\scriptstyle x}; "1_x";(0,-12)*{}; **\dir{-}
?(.4)*\dir{<}+(3,0)*{\scriptstyle x}; (0,0)*{1_x};
\endxy
\]
The identity object will not be drawn in the diagrams, but merely
implied. As an example of this, consider how we obtain the string
diagram for $i_{x} \maps 1 \to x \ten \xb$:
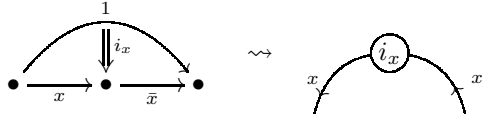
\[
\vcenter{ \xymatrix{
 \bullet \ar[r]_{x} \ar@/^2pc/[rr]^{1}_{}="1"
 & \bullet \ar[r]_{\xb} \ar@{=>}^<<<{i_{x}}"1";{}
 & \bullet
} } \quad \rightsquigarrow \quad \xy
 (0,0)*{i_{x}};
 (0,0)*\xycircle(2.5,2.5){-}="i";
 (-10,-8)*{};
 **\crv{(-8,0)} ?(.16)*\dir{<}+(-1,2)*{\scriptstyle x};
 (10,-8)*{};"i";
 **\crv{(8,0)} ?(.75)*\dir{<}+(3,1)*{\scriptstyle x};
\endxy
\]
Note that we omit the incoming string corresponding to the
identity object $1$.  Also, we indicate weak inverse objects with
arrows `going backwards in time', following this rule:
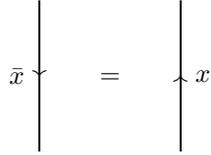
\[
 \xy
 (0,10)*{};
 (0,-10)*{}
**\dir{-}?(.5)*\dir{>}+(-3,0)*{\xb};
 \endxy
\qquad = \qquad
 \xy
 (0,10)*{};
 (0,-10)*{}
**\dir{-}?(.5)*\dir{<}+(3,0)*{x};
 \endxy
\]
In calculations, it is handy to draw the unit $i_x$ in an even
more abbreviated form:
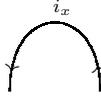
\[
\xy
 (6,0)*{};
 (-6,0)*{}
   **\crv{(6,12) & (-6,12)} ?(.1)*\dir{>} ?(.92)*\dir{>};
 (1,11)*{\scriptstyle i_{x}};
\endxy
\]
where we omit the disc surrounding the morphism label `$i_x$', and
it is understood that the downward pointing arrow corresponds to
$x$ and the upward pointing arrow to $\xb$.  Similarly, we draw
the morphism $e_x$ as
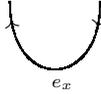
\[
\xy
 (6,0)*{};
 (-6,0)*{}
   **\crv{(6,-12) & (-6,-12)} ?(.1)*\dir{>} ?(.92)*\dir{>};
 (1,-11)*{\scriptstyle e_{x}};
\endxy
\]

In a strict monoidal category, where the associator and the left
and right unit laws are identity morphisms, one can interpret any
string diagram as a morphism in a unique way.  In fact, Joyal and
Street have proved some rigorous theorems to this effect
\cite{JS0}.  With the help of Mac Lane's coherence theorem
\cite{MacLane2} we can also do this in a weak monoidal category. To
do this, we interpret any string of objects and $1$'s as a tensor
product of objects where all parentheses start in front and all
1's are removed. Using the associator and left/right unit laws to
do any necessary reparenthesization and introduction or
elimination of 1's, any string diagram then describes a morphism
between tensor products of this sort. The fact that this morphism
is unambiguously defined follows from Mac Lane's coherence
theorem.

For a simple example of string diagram technology in action,
consider the zig-zag identities.  To begin with, these say that
the following diagrams commute:
\[
\xymatrix@C=2.4pc{
  1 \tensor x \ar[r]^<<<<<{i_x \tensor 1} \ar[d]_{\ell_{x}}
  & (x \tensor \bar{x}) \tensor x \ar[r]^{a_{x, \bar{x},x} }
  & x \tensor ( \bar{x} \tensor x) \ar[d]^{1 \tensor e_x} \; \;\\
  x \ar[rr]_{r^{-1}_x}
   && x \tensor 1  }
\]
\[
\xymatrix@C=2.4pc{
 \bar{x} \ten 1 \ar[r]^<<<<<<{1 \ten i_x} \ar[d]_{r_{\bar{x}}}
  & \bar{x} \ten (x \ten \bar{x}) \ar[r]^{a^{-1}_{\bar{x},x, \bar{x}}}
  & (\bar{x} \ten x) \ten \bar{x} \ar[d]^{e_x \ten 1} \\
 \bar{x} \ar[rr]_{\ell^{-1}_{\bar{x}}}
  && 1 \ten x   } \quad
\]
In globular notation these diagrams become:
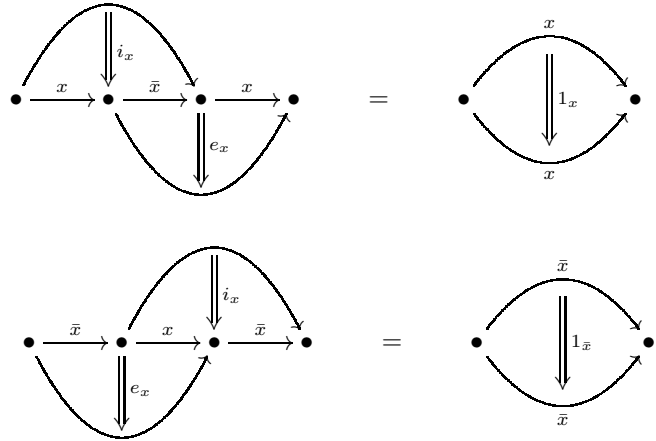
\[
\xymatrix{ \bullet
  \ar[r]^x
  \ar@/^3pc/[rr]^{\quad}_{}="0"
& \bullet
  \ar[r]^{\bar{x}}
  \ar@/_3pc/[rr]^{\quad}_{}="2"
  \ar@{=>}"0"; {} ^{i_x}
& \bullet
  \ar[r]^x
  \ar@{=>}{}; "2" ^{e_x}
& \bullet } \qquad = \qquad \xymatrix{
 \bullet
   \ar@/^2pc/[rr]_{\quad}^{x}="1"
   \ar@/_2pc/[rr]_{x}="2"
&& \bullet
 \ar@{}"1";"2"|(.2){\,}="7"
  \ar@{}"1";"2"|(.8){\,}="8"
 \ar@{=>}"7" ;"8"^{1_x}   }  \quad
\]
\[
\xymatrix{ \bullet
  \ar[r]^{\bar{x}}
  \ar@/_3pc/[rr]^{\quad}_{}="0"
& \bullet
  \ar[r]^x
  \ar@/^3pc/[rr]^{\quad}_{}="2"
  \ar@{=>};"0" ^{e_x}
& \bullet
  \ar[r]^{\bar{x}}
  \ar@{=>}"2";{} ^{i_x}
& \bullet } \qquad = \qquad \xymatrix{
 \bullet
   \ar@/^2pc/[rr]_{\quad}^{\xb}="1"
   \ar@/_2pc/[rr]_{\xb}="2"
&& \bullet
 \ar@{}"1";"2"|(.2){\,}="7"
  \ar@{}"1";"2"|(.8){\,}="8"
 \ar@{=>}"7" ;"8"^{1_{\xb}}   }
\]
Taking Poincar\'e duals, we obtain the zig-zag identities in
string diagram form:
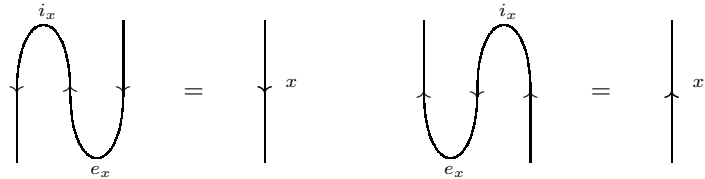
\[
\xy0;/r.28pc/: (-6,-8)*{}="1E"; (-6,0)*{}="1"; (0,0)*{}="2";
(6,0)*{}="3"; (6,8)*{}="3B";
 "2";"1" **\crv{(0,10)& (-6,10)}
     ?(.02)*\dir{>}  ?(1)*\dir{>};
 "3";"2" **\crv{(6,-10)& (0,-10)} ;
 "1";"1E" **\dir{-};
 "3B";"3" **\dir{-}?(1)*\dir{>} ;
(-2.5,9)*{\scriptstyle i_x}; (3.5,-9)*{\scriptstyle e_x};
\endxy
\qquad = \xy0;/r.28pc/: (-6,8)*{}; (0,8)*{}; (0,-8)*{}; **\dir{-}
?(.5)*\dir{<}; (6,-8)*{}?(.5)*\dir{<}+(3,1)*{\scriptstyle x};
\endxy
\qquad  \qquad \xy0;/r.28pc/: (-6,8)*{}="1E"; (-6,0)*{}="1";
(0,0)*{}="2"; (6,0)*{}="3"; (6,-8)*{}="3B";
 "2";"1" **\crv{(0,-10)& (-6,-10)}
     ?(.02)*\dir{>}  ?(1)*\dir{>};
 "3";"2" **\crv{(6,10)& (0,10)}
       ;
 "1";"1E" **\dir{-};
 "3B";"3" **\dir{-}?(1)*\dir{>};
(-2.5,-9)*{\scriptstyle e_x}; (3.5,9)*{\scriptstyle i_x};
\endxy
\qquad = \xy0;/r.28pc/: (-6,8)*{}; (0,-8)*{}; (0,8)*{}; **\dir{-}
?(.5)*\dir{<}; (6,-8)*{}?(.5)*\dir{<}+(3,1)*{\scriptstyle x};
\endxy
\]
This picture explains their name!  The zig-zag identities simply
allow us to straighten a piece of string.

In most of our calculations we only need string diagrams where all
strings are labelled by $x$ and $\xb$.  In this case we can omit
these labels and just use downwards or upwards arrows to
distinguish between $x$ and $\xb$.  We draw $i_x$ as
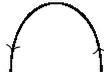
\[
\xy
 (6,0)*{};
 (-6,0)*{}
   **\crv{(6,12) & (-6,12)} ?(.1)*\dir{>} ?(.92)*\dir{>};
\endxy
\]
and draw $e_x$ as
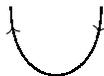
\[
\xy
 (6,0)*{};
 (-6,0)*{}
   **\crv{(6,-12) & (-6,-12)} ?(.1)*\dir{>} ?(.92)*\dir{>};
\endxy
\]
The zig-zag identities become just:
\[
\xy0;/r.28pc/: (-6,-8)*{}="1E"; (-6,0)*{}="1"; (0,0)*{}="2";
(6,0)*{}="3"; (6,8)*{}="3B";
 "2";"1" **\crv{(0,10)& (-6,10)}
     ?(.02)*\dir{>}  ?(1)*\dir{>};
 "3";"2" **\crv{(6,-10)& (0,-10)} ;
 "1";"1E" **\dir{-};
 "3B";"3" **\dir{-}?(1)*\dir{>} ;
(-2.5,9)*{}; (3.5,-9)*{};
\endxy
\qquad = \xy0;/r.28pc/: (-6,8)*{}; (0,8)*{}; (0,-8)*{}; **\dir{-}
?(.5)*\dir{<}; (6,-8)*{}?(.5)*\dir{<}+(3,1)*{};
\endxy
\qquad  \qquad \xy0;/r.28pc/: (-6,8)*{}="1E"; (-6,0)*{}="1";
(0,0)*{}="2"; (6,0)*{}="3"; (6,-8)*{}="3B";
 "2";"1" **\crv{(0,-10)& (-6,-10)}
     ?(.02)*\dir{>}  ?(1)*\dir{>};
 "3";"2" **\crv{(6,10)& (0,10)}
       ;
 "1";"1E" **\dir{-};
 "3B";"3" **\dir{-}?(1)*\dir{>};
(-2.5,-9)*{}; (3.5,9)*{};
\endxy
\qquad = \xy0;/r.28pc/: (-6,8)*{}; (0,-8)*{}; (0,8)*{}; **\dir{-}
?(.5)*\dir{<}; (6,-8)*{}?(.5)*\dir{<}+(3,1)*{};
\endxy
\]

We also obtain some rules for manipulating string diagrams just
from the fact that $i_x$ and $e_x$ have inverses. For these, we
draw $i_x^{-1}$ as
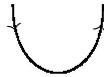
\[
\xy
 (-6,0)*{};
 (6,0)*{}
   **\crv{(-6,-12) & (6,-12)} ?(.1)*\dir{>} ?(.92)*\dir{>};
\endxy
\]
and $e_x^{-1}$ as
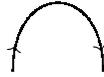
\[
\xy
 (-6,0)*{};
 (6,0)*{}
   **\crv{(-6,12) & (6,12)} ?(.1)*\dir{>} ?(.92)*\dir{>};
\endxy
\]
The equations $i_x i_x^{-1} = 1_1$ and $e_x^{-1} e_x = 1_1$ give
the rules
\[
\xy
 (5,0)*{}="1";
 (-5,0)*{}="2";
   "1";"2" **\crv{(5,12)&(-5,12)}
     ?(0)*\dir{>}  ?(1)*\dir{>} ;
   "2";"1" **\crv{ (-5,-12)& (5,-12)}
\endxy
\quad = \quad \qquad \qquad \qquad \qquad \qquad \xy
 (5,0)*{}="1";
 (-5,0)*{}="2";
   "2";"1" **\crv{(-5,12)&(5,12)}
     ?(0)*\dir{>}  ?(1)*\dir{>} ;
   "1";"2" **\crv{(5,-12) & (-5,-12)}
\endxy
\quad = \quad \qquad \qquad
\]
which mean that in a string diagram, a loop of either form may be
removed or inserted without changing the morphism described by the
diagram. Similarly, the equations $e_x e_x^{-1} = 1_{\bar{x}
\tensor x}$ and $i_x^{-1}i_x = 1_{x \tensor \bar{x}}$ give the
rules
\[
\xy
 (5,10)*{};
 (-5,10)*{}
   **\crv{(5,-1) & (-5,-1)} ?(.12)*\dir{>} ?(.92)*\dir{>};
 (-5,-10)*{};
 (5,-10)*{}
   **\crv{(-5,1) & (5,1)} ?(.12)*\dir{>} ?(.92)*\dir{>};
\endxy
\quad = \quad \xy (-4,-10)*{};(-4,10)**\dir{-} ?(.52)*\dir{>};
(4,10)*{};(4,-10)**\dir{-} ?(.5)*\dir{>};
\endxy
\quad \qquad \qquad \qquad \quad \xy
 (-5,10)*{};
 (5,10)*{}
   **\crv{(-5,-1) & (5,-1)} ?(.12)*\dir{>} ?(.92)*\dir{>};
 (5,-10)*{};
 (-5,-10)*{}
   **\crv{(5,1) & (-5,1)} ?(.12)*\dir{>} ?(.92)*\dir{>};
\endxy
\quad = \quad \xy (-4,10)*{};(-4,-10)**\dir{-} ?(.5)*\dir{>};
(4,-10)*{};(4,10)**\dir{-} ?(.52)*\dir{>};
\endxy
\]
Again, these rules mean that in a string diagram we can modify any
portion as above without changing the morphism in question.

By taking the inverse of both sides in the zig-zag identities, we
obtain extra zig-zag identities involving $i_x^{-1}$ and
$e_x^{-1}$:
\[
\xy0;/r.28pc/: (6,8)*{}="1E"; (6,0)*{}="1"; (0,0)*{}="2";
(-6,0)*{}="3"; (-6,-8)*{}="3B";
 "2";"1" **\crv{(0,-10)& (6,-10)}
     ?(.02)*\dir{>}  ?(1)*\dir{>};
 "3";"2" **\crv{(-6,10)& (0,10)};
 "1";"1E" **\dir{-};
 "3B";"3" **\dir{-}?(1)*\dir{>};
\endxy
\qquad = \xy0;/r.28pc/: (6,8)*{}; (0,-8)*{}; (0,8)*{}; **\dir{-}
?(.5)*\dir{<}; (-6,-8)*{}?(.5)*\dir{<}+(3,1)*{};
\endxy
\qquad  \qquad \xy0;/r.28pc/: (6,-8)*{}="1E"; (6,0)*{}="1";
(0,0)*{}="2"; (-6,0)*{}="3"; (-6,8)*{}="3B";
 "2";"1" **\crv{(0,10)& (6,10)}
     ?(.02)*\dir{>}  ?(1)*\dir{>};
 "3";"2" **\crv{(-6,-10)& (0,-10)} ;
 "1";"1E" **\dir{-};
 "3B";"3" **\dir{-}?(1)*\dir{>} ;
\endxy
\qquad = \xy0;/r.28pc/: (6,8)*{}; (0,8)*{}; (0,-8)*{}; **\dir{-}
?(.5)*\dir{<}; (-6,-8)*{}?(.5)*\dir{<}+(3,1)*{};
\endxy
\]
Conceptually, this means that whenever $(x,\xb,i_x,e_x)$ is an
adjoint equivalence, so is $(\xb,x,e_x^{-1},i_x^{-1})$.

In the calculations to come we shall also use another rule, the
`horizontal slide':
\[
\vcenter{\xy (5,6)*{};
 (5,0)*{};
 (-5,0)*{}
   **\crv{(5,-11) & (-5,-11)} ?(.1)*\dir{>} ?(.92)*\dir{>};
 (1,-10)*{\scriptstyle e_{x}};
\endxy} \quad
\vcenter{\xy
 (-5,0)*{};
 (5,0)*{}
   **\crv{(-5,11) & (5,11)} ?(.1)*\dir{>} ?(.92)*\dir{>};
 (1,10)*{\scriptstyle e_{y}^{-1}};
\endxy}
\qquad = \qquad \xy
 (5,13)*{};
 (-5,13)*{}
   **\crv{(5,2) & (-5,2)} ?(.12)*\dir{>} ?(.92)*\dir{>};
 (-5,-13)*{};
 (5,-13)*{}
   **\crv{(-5,-2) & (5,-2)} ?(.12)*\dir{>} ?(.92)*\dir{>};
(0,3)*{\scriptstyle e_x};
     (1,-3)*{\scriptstyle e_y^{-1}};
\endxy
\]
This follows from general results on the isotopy-invariance of the
morphisms described by string diagrams \cite{JS}, but it also
follows directly from the interchange law relating vertical and
horizontal composition in a 2-category:
\[
\xy (-23,0)*{}; (23,0)*{}; (-22,0)*+{\bullet}="L";
(0,0)*+{\bullet}="M"; (22,0)*+{\bullet}="R";
(-11,10)*+{\bullet}="LT"; (11,-10)*+{\bullet}="RB";
 (-11,7)="L1";
 (-11,-7)="L2";
 (11,7)="R1";
 (11,-7)="R2";
 "L";"M" **\crv{(-17,-12) & (-5,-12)};
             ?(.95)*\dir{>};
 "M";"R" **\crv{ (5,12)& (17,12)};
             ?(.95)*\dir{>};
 "L";"LT" **\crv{(-20,7.5)};
             ?(.92)*\dir{>};
 "LT";"M" **\crv{(-2,7.5)};
             ?(.92)*\dir{>};
 "M";"RB" **\crv{(2,-7.5)};
             ?(.92)*\dir{>};
 "RB";"R" **\crv{(20,-7.5)};
             ?(.92)*\dir{>};
 {\ar@{=>}^{e_x} "L1";"L2"};
 {\ar@{=>}^{e_y^{-1}} "R1";"R2"};
     (-18,9)*{\scriptstyle \xb};
     (-4,9)*{\scriptstyle x};
     (4,-9)*{\scriptstyle \bar{y}};
     (18,-9)*{\scriptstyle y};
     (-11,-11)*{\scriptstyle 1};
     (11,11)*{\scriptstyle 1};
\endxy
\qquad = \qquad \xy (-23,0)*{}; (23,0)*{}; (-22,0)*+{\bullet}="L";
(0,0)*+{\bullet}="M"; (22,0)*+{\bullet}="R";
(-11,10)*+{\bullet}="LT"; (11,-10)*+{\bullet}="RB";
 (-11,8)="L1";
 (-11,2)="L2";
 (11,-2)="R1";
 (11,-8)="R2";
 (-11,-8)="L4";
 (-11,-2)="L3";
 (11,2)="R4";
 (11,7.5)="R3";
 "L";"M" **\crv{(-17,-12) & (-5,-12)};
             ?(.95)*\dir{>};
 "M";"R" **\crv{ (5,12)& (17,12)};
             ?(.95)*\dir{>};
 "L";"LT" **\crv{(-20,7.5)};
             ?(.92)*\dir{>};
 "LT";"M" **\crv{(-2,7.5)};
             ?(.92)*\dir{>};
 "M";"RB" **\crv{(2,-7.5)};
             ?(.92)*\dir{>};
 "RB";"R" **\crv{(20,-7.5)};
             ?(.92)*\dir{>};
 {\ar "L";"M"};
 {\ar "M";"R"};
 {\ar@{=>}^{e_x} "L1";"L2"};
 {\ar@{=>}^{e_y^{-1}} "R1";"R2"};
 {\ar@{=>}^{1_1} "L3";"L4"};
 {\ar@{=>}^{1_1} "R3";"R4"};
     (-18,9)*{\scriptstyle \xb};
     (-4,9)*{\scriptstyle x};
     (4,-9)*{\scriptstyle \bar{y}};
     (18,-9)*{\scriptstyle y};
     (-11,-11)*{\scriptstyle 1};
     (11,11)*{\scriptstyle 1};
\endxy
\]
\[
\xy (-23,0)*{}; (23,0)*{};
\endxy
\qquad = \qquad \xy (-23,0)*{}; (23,0)*{}; (-22,0)*+{\bullet}="L";
(0,0)*+{\bullet}="M"; (-11,10)*+{\bullet}="LT";
(-11,-10)*+{\bullet}="RB";
 (-11,8)="L1";
 (-11,2)="L2";
 (-11,-8)="L4";
 (-11,-2)="L3";
 "L";"LT" **\crv{(-20,7.5)};
             ?(.92)*\dir{>};
 "L";"RB" **\crv{(-20,-7.5)};
             ?(.92)*\dir{>};
 "LT";"M" **\crv{(-2,7.5)};
             ?(.92)*\dir{>};
 "RB";"M" **\crv{(-2,-7.5)};
             ?(.92)*\dir{>};
 {\ar "L";"M"};
 {\ar@{=>}^{e_x} "L1";"L2"};
 {\ar@{=>}^{e_y^{-1}} "L3";"L4"};
     (-18,9)*{\scriptstyle \xb};
     (-4,9)*{\scriptstyle x};
     (-18,-9)*{\scriptstyle \bar{y}};
     (-4,-9)*{\scriptstyle y};
\endxy
\]

\noindent We will also be using other slightly different versions
of the horizontal slide, which can be proved the same way.

As an illustration of how these rules are used, we give a string
diagram proof of a result due to Saavedra Rivano \cite{Rivano},
which allows a certain simplification in the definition of
`coherent 2-group':

\begin{prop} \et
Let $C$ be a weak monoidal category, and let $x,\xb \in C$ be
objects equipped with isomorphisms $i_x \maps 1 \to x \ten \xb$
and $e_x \maps \xb \ten x \to 1$.  If the quadruple
$(x,\xb,i_x,e_x)$ satisfies either one of the zig-zag identities,
it automatically satisfies the other as well.
\end{prop}

\textbf{Proof. } Suppose the first zig-zag identity holds:
\[
\xy (-6,-8)*{}="1E"; (-6,0)*{}="1"; (0,0)*{}="2"; (6,0)*{}="3";
(6,8)*{}="3B";
 "2";"1" **\crv{(0,10)& (-6,10)}
     ?(.0)*\dir{>}  ?(1)*\dir{>};
 "3";"2" **\crv{(6,-10)& (0,-10)};
 "1";"1E" **\dir{-};
 "3B";"3" **\dir{-}?(1)*\dir{>} ;
\endxy
\qquad = \xy (-6,8)*{}; (0,8)*{}; (0,-8)*{}; **\dir{-}
?(.47)*\dir{<}; (6,-8)*{};
\endxy
\]

\noindent Then the second zig-zag identity may be shown as
follows:
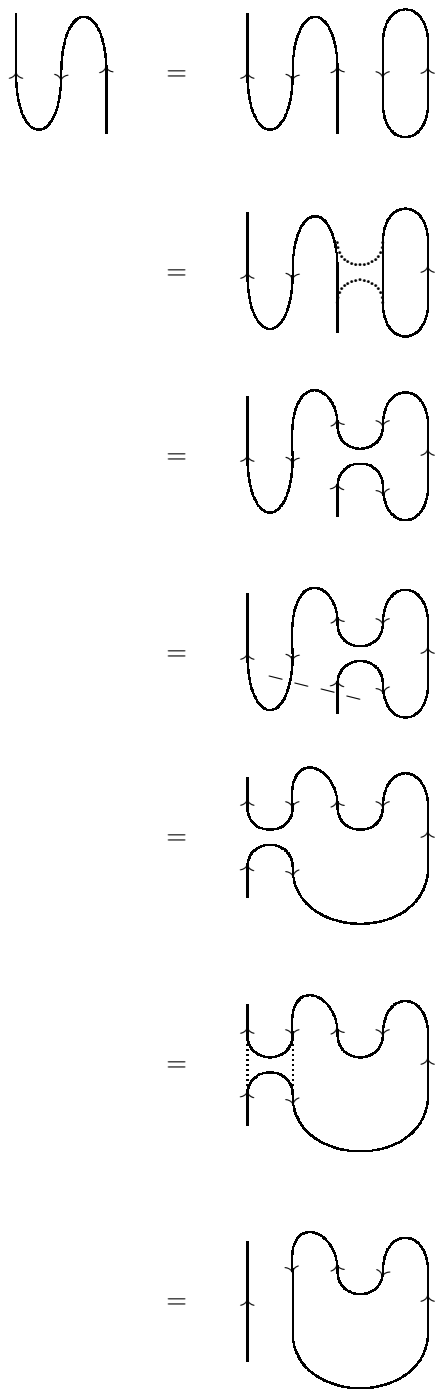
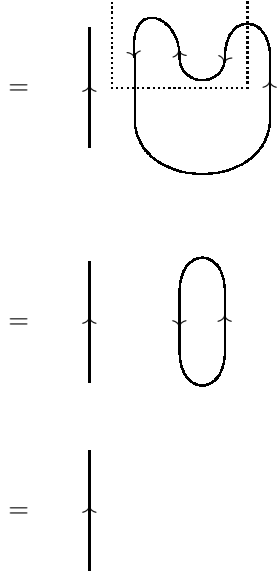
\[
\xy (-6,8)*{}="1E"; (-6,0)*{}="1"; (0,0)*{}="2"; (6,0)*{}="3";
(6,-8)*{}="3B";
 "2";"1" **\crv{(0,-10)& (-6,-10)}
     ?(.03)*\dir{>}  ?(1)*\dir{>};
 "3";"2" **\crv{(6,10)& (0,10)}
     ?(.03)*\dir{>}  ;
 "1";"1E" **\dir{-};
 "3B";"3" **\dir{-};
\endxy
\qquad = \qquad \xy (-6,8)*{}="1E"; (-6,0)*{}="1"; (0,0)*{}="2";
(6,0)*{}="3"; (6,-8)*{}="3B";
 "2";"1" **\crv{(0,-10)& (-6,-10)}
     ?(.03)*\dir{>}  ?(1)*\dir{>};
 "3";"2" **\crv{(6,10)& (0,10)}
     ?(.03)*\dir{>}  ;
 "1";"1E" **\dir{-};
 "3B";"3" **\dir{-};
   (12,4)*{}="1";
   (12,-4)*{}="A";
   (18,-4)*{}="B";
   (18,4)*{}="B'";
     "B'";"1" **\crv{(18,10)& (12,10)};
     "A";"B" **\crv{(12,-10)&(18,-10)};
     "B";"B'" **\dir{-} ?(.60)*\dir{>};
     "1";"A" **\dir{-};  ?(.57)*\dir{>};
\endxy
\]

\[
\xy (-6,0)*{}="1E"; (6,0)*{}="3B";
\endxy
\qquad = \qquad \xy (-6,8)*{}="1E"; (-6,0)*{}="1"; (0,0)*{}="2";
(6,0)*{}="3"; (6,-8)*{}="3B";
 "2";"1" **\crv{(0,-10)& (-6,-10)}
     ?(.03)*\dir{>}  ?(1)*\dir{>};
 "3";"2" **\crv{(6,10)& (0,10)};
 "1";"1E" **\dir{-};
 "3B";"3" **\dir{-};
   (12,4)*{}="1";
   (12,-4)*{}="A";
   (18,-4)*{}="B";
   (18,4)*{}="B'";
     "B'";"1" **\crv{(18,10)& (12,10)};
     "A";"B" **\crv{(12,-10)&(18,-10)};
     "B";"B'" **\dir{-} ?(.60)*\dir{>};
     "1";"A" **\dir{-};
  (6,4)*{}="D1";
  (12,4)*{}="D2";
  (6,-4)*{}="F1";
  (12,-4)*{}="F2";
   "D1";"D2" **\crv{~*=<2pt>{.}(6,0)&(12,0)};
   "F1";"F2" **\crv{~*=<2pt>{.}(6,0)&(12,0)};
\endxy
\]
\[
\xy (-6,0)*{}="1E"; (6,0)*{}="3B";
\endxy
\qquad = \qquad \xy (-6,8)*{}="1E"; (-6,0)*{}="1"; (0,0)*{}="2";
(6,4)*{}="3"; (6,-8)*{}="3B";
 "2";"1" **\crv{(0,-10)& (-6,-10)}
     ?(.03)*\dir{>}  ?(1)*\dir{>};
 "3";"2" **\crv{(6,10)& (-1,12)} ?(.05)*\dir{>};
 "1";"1E" **\dir{-};
   (12,4)*{}="1";
   (12,-4)*{}="A";
   (18,-4)*{}="B";
   (18,4)*{}="B'";
     "B'";"1" **\crv{(18,10)& (12,10)};
     "A";"B" **\crv{(12,-10)&(18,-10)} ?(.07)*\dir{>};
     "B";"B'" **\dir{-} ?(.60)*\dir{>};
  (6,4)*{}="D1";
  (12,4)*{}="D2";
  (6,-4)*{}="F1";
  (12,-4)*{}="F2";
   "D1";"D2" **\crv{(6,0)&(12,0)} ?(1)*\dir{<};
   "F1";"F2" **\crv{(6,0)&(12,0)}
     ?(.03)*\dir{>};
       "3B";"F1" **\dir{-};
\endxy
\]

\[
\xy (-6,0)*{}="1E"; (6,0)*{}="3B";
\endxy
\qquad = \qquad \xy (-6,8)*{}="1E"; (-6,0)*{}="1"; (0,0)*{}="2";
(6,4)*{}="3"; (6,-8)*{}="3B";
 "2";"1" **\crv{(0,-10)& (-6,-10)}
     ?(.03)*\dir{>}  ?(1)*\dir{>};
 "3";"2" **\crv{(6,10)& (-1,12)} ?(.05)*\dir{>};
 "1";"1E" **\dir{-};
   (12,4)*{}="1";
   (12,-4)*{}="A";
   (18,-4)*{}="B";
   (18,4)*{}="B'";
     "B'";"1" **\crv{(18,10)& (12,10)};
     "A";"B" **\crv{(12,-10)&(18,-10)} ?(.07)*\dir{>};
     "B";"B'" **\dir{-} ?(.60)*\dir{>};
  (6,4)*{}="D1";
  (12,4)*{}="D2";
  (6,-4)*{}="F1";
  (12,-4)*{}="F2";
   "D1";"D2" **\crv{(6,0)&(12,0)} ?(1)*\dir{<};
   "F1";"F2" **\crv{(6,0)&(12,0)}
     ?(.03)*\dir{>};
       "3B";"F1" **\dir{-};
 (-3,-3)*{}="X";
 (9,-6)*{}="XX";
 "X";"XX" **\dir{--};
\endxy
\]
\[
\xy (-6,0)*{}="1E"; (6,0)*{}="3B";
\endxy
\qquad = \qquad \xy (-6,8)*{}="1E"; (-6,4)*{}="1"; (0,4)*{}="2";
(6,4)*{}="3"; (-6,-8)*{}="3B";
 "2";"1" **\crv{(0,0)& (-6,0)};
     ?(.0)*\dir{>};
 "3";"2" **\crv{(6,10)& (-1,12)} ?(.05)*\dir{>};
 "1";"1E" **\dir{-}; ?(.25)*\dir{>};
   (12,4)*{}="1";
   (0,-4)*{}="A";
   (18,-4)*{}="B";
   (18,4)*{}="B'";
     "B'";"1" **\crv{(18,10)& (12,10)};
     "A";"B" **\crv{(0,-14)&(18,-14)} ?(.04)*\dir{>};
     "B";"B'" **\dir{-} ?(.60)*\dir{>};
  (6,4)*{}="D1";
  (12,4)*{}="D2";
  (-6,-4)*{}="F1";
  (0,-4)*{}="F2";
   "D1";"D2" **\crv{(6,0)&(12,0)} ?(1)*\dir{<};
   "F1";"F2" **\crv{(-6,0)&(0,0)}
     ?(.03)*\dir{>};
       "3B";"F1" **\dir{-};
\endxy
\]

\[
\xy (-6,0)*{}="1E"; (6,0)*{}="3B";
\endxy
\qquad = \qquad \xy (-6,8)*{}="1E"; (-6,4)*{}="1"; (0,4)*{}="2";
(6,4)*{}="3"; (-6,-8)*{}="3B";
 "2";"1" **\crv{(0,0)& (-6,0)}
     ?(.0)*\dir{>};
 "3";"2" **\crv{(6,10)& (-1,12)} ?(.05)*\dir{>};
 "1";"1E" **\dir{-}; ?(.25)*\dir{>};
   (12,4)*{}="1A";
   (0,-4)*{}="A";
   (18,-4)*{}="B";
   (18,4)*{}="B'";
     "B'";"1A" **\crv{(18,10)& (12,10)};
     "A";"B" **\crv{(0,-14)&(18,-14)} ?(.04)*\dir{>};
     "B";"B'" **\dir{-} ?(.60)*\dir{>};
  (6,4)*{}="D1";
  (12,4)*{}="D2";
  (-6,-4)*{}="F1";
  (0,-4)*{}="F2";
   "D1";"D2" **\crv{(6,0)&(12,0)} ?(1)*\dir{<};
   "F1";"F2" **\crv{(-6,0)&(0,0)}
     ?(.03)*\dir{>};
       "3B";"F1" **\dir{-};
   "1";"F1" **\dir{.};
   "2";"F2" **\dir{.};
\endxy
\]

\[
\xy (-6,0)*{}="1E"; (6,0)*{}="3B";
\endxy
\qquad = \qquad \xy (-6,8)*{}="1E"; (-6,4)*{}="1"; (0,4)*{}="2";
(6,4)*{}="3"; (-6,-8)*{}="3B";
 "3";"2" **\crv{(6,10)& (-1,12)}
   ?(.05)*\dir{>} ?(1)*\dir{>};
 "1";"1E" **\dir{-};
   (12,4)*{}="1A";
   (0,-4)*{}="A";
   (18,-4)*{}="B";
   (18,4)*{}="B'";
     "B'";"1A" **\crv{(18,10)& (12,10)};
     "A";"B" **\crv{(0,-14)&(18,-14)};
     "B";"B'" **\dir{-} ?(.60)*\dir{>};
  (6,4)*{}="D1";
  (12,4)*{}="D2";
  (-6,-4)*{}="F1";
  (0,-4)*{}="F2";
   "D1";"D2" **\crv{(6,0)&(12,0)} ?(.95)*\dir{<};
       "3B";"F1" **\dir{-};
   "1";"F1" **\dir{-} ?(.47)*\dir{<};
   "2";"F2" **\dir{-};
\endxy
\]

\[
\xy (-6,0)*{}="1E"; (6,0)*{}="3B";
\endxy
\qquad = \qquad \xy (-6,8)*{}="1E"; (-6,4)*{}="1"; (0,4)*{}="2";
(6,4)*{}="3"; (-6,-8)*{}="3B";
 "3";"2" **\crv{(6,10)& (-1,12)}
   ?(.05)*\dir{>} ?(1)*\dir{>};
 "1";"1E" **\dir{-};
   (12,4)*{}="1A";
   (0,-4)*{}="A";
   (18,-4)*{}="B";
   (18,4)*{}="B'";
     "B'";"1A" **\crv{(18,10)& (12,10)};
     "A";"B" **\crv{(0,-14)&(18,-14)};
     "B";"B'" **\dir{-} ?(.60)*\dir{>};
  (6,4)*{}="D1";
  (12,4)*{}="D2";
  (-6,-4)*{}="F1";
  (0,-4)*{}="F2";
   "D1";"D2" **\crv{(6,0)&(12,0)} ?(.95)*\dir{<};
       "3B";"F1" **\dir{-};
   "1";"F1" **\dir{-} ?(.47)*\dir{<};
   "2";"F2" **\dir{-};
         (-3,12)*{}="X1";
         (-3,0)*{}="X4";
         (15,12)*{}="X2";
         (15,0)*{}="X3";
             "X1";"X2" **\dir{.};
             "X2";"X3" **\dir{.};
             "X3";"X4" **\dir{.};
             "X4";"X1" **\dir{.};
\endxy
\]

\[
\xy (-6,0)*{}="1E"; (6,0)*{}="3B";
\endxy
\qquad = \qquad \xy (-6,8)*{}="1E"; (-6,-8)*{}="1"; (18,0)*{};
 "1";"1E" **\dir{-} ?(.53)*\dir{>};
   (6,4)*{}="1";
   (6,-4)*{}="A";
   (12,-4)*{}="B";
   (12,4)*{}="B'";
     "B'";"1" **\crv{(12,10)& (6,10)};
     "A";"B" **\crv{(6,-10)&(12,-10)};
     "B";"B'" **\dir{-} ?(.60)*\dir{>};
     "1";"A" **\dir{-};  ?(.57)*\dir{>};
\endxy
\]

\[
\xy (-6,0)*{}="1E"; (6,0)*{}="3B";
\endxy
\qquad = \qquad \xy (-6,8)*{}="1E"; (-6,-8)*{}="1"; (18,0)*{};
 "1";"1E" **\dir{-} ?(.53)*\dir{>};
\endxy
\]
In this calculation, we indicate an application of the `horizontal
slide' rule by a dashed line.  Dotted curves or lines indicate
applications of the rule $e_x e_x^{-1} = 1_{\xb \tensor x}$.  A
box indicates an application of the first zig-zag identity. The
converse can be proven similarly. \qed

\section{Improvement} \label{improvementsection}

We now use string diagrams to show that any weak 2-group can be
improved to a coherent one.  There are shorter proofs, but none
quite so pretty, at least in a purely visual sense.  Given
a weak 2-group $C$ and any object $x \in C$, we can choose a weak
inverse $\xb$ for $x$ together with isomorphisms $i_x \maps 1
\rightarrow x \ten \xb$, $e_x \maps \xb \ten x \rightarrow 1$.
From this data we shall construct an adjoint equivalence
$(x,\xb,i'_x,e'_x)$.   By doing this for every object of $C$,
we make $C$ into a coherent 2-group.

\begin{thm} \label{improve} \et
Any weak 2-group $C$ can be given the structure of a coherent
2-group $\imp(C)$ by equipping each object with an adjoint
equivalence.
\end{thm}

\textbf{Proof. } First, for each object $x$ we choose a weak
inverse $\xb$ and isomorphisms $i_x \maps 1 \to x \ten \xb$, $e_x
\maps \xb \ten x \to 1$. From this data we construct an adjoint
equivalence $(x,\xb,i'_x,e'_x)$.  To do this, we set $e'_x = e_x$
and define $i'_x$ as the following composite morphism:

\[
\def\objectstyle{\scriptstyle}
\def\labelstyle{\scriptstyle}
\xy (-60,0)*+{1}="1"; (-52,0)*+{x\bar{x}}="2";
(-40,0)*+{x(1\bar{x})}="3"; (-25,0)*+{x((\bar{x}x)\bar{x})}="4";
(-5,0)*+{ x(\bar{x}(x\bar{x}))}="5";
(15,0)*+{(x\bar{x})(x\bar{x})}="6"; (35,0)*+{1(x\bar{x})}="7";
(50,0)*+{(1x)\bar{x}}="8"; (60,0)*+{x\bar{x} . }="9";
 {\ar^{i_x} "1";"2"};
 {\ar^{x\ell^{-1}_{\bar{x}}} "2";"3"};
 {\ar^{x  e^{-1}_{x} \bar{x}} "3";"4"};
 {\ar^{x a_{\bar{x},x,\bar{x}}} "4";"5"};
 {\ar^{a^{-1}_{x,\bar{x},x\bar{x}}} "5";"6"};
 {\ar^{i^{-1}_x  (x\bar{x})} "6";"7"};
 {\ar^{a^{-1}_{1,x,\bar{x}}} "7";"8"};
 {\ar^{\ell_x \bar{x}} "8";"9"};
\endxy
\]
where we omit tensor product symbols for brevity.

The above rather cryptic formula for $i'_x$ becomes much clearer
if we use pictures.  If we think of a weak 2-group as a one-object
2-category and write this formula in globular notation it becomes:
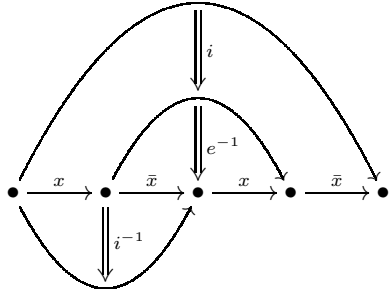
\[
\xymatrix{
  \bullet 
\ar[r]^x
 \ar@/^6pc/[rrrr]^{\quad}_{}="0"
 \ar@/_3pc/[rr]^{\quad}_{}="3"
&  \bullet 
\ar[r]^{\bar{x}}
   \ar@/^3pc/[rr]^{\quad}_{}="1"
        \ar@{=>} {} ;"3"  ^{i^{-1}}
&  \bullet 
\ar[r]^x
        \ar@{=>}"1"; {} ^{e^{-1}}
&  \bullet 
\ar[r]^{\bar{x}}
        \ar@{=>}"0"; "1" ^{i}
&  \bullet 
 }
\]
where we have suppressed associators and the left unit law for
clarity.  If we write it as a string diagram it looks even
simpler:
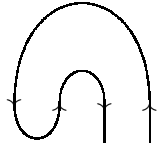
\[
\xy
 (-6,-0)*{}="1";
 (0,0)*{}="2";
 (6,0)*{}="3";
 (6,-5)*{}="3'";
 (12,0)*{}="4";
 (12,-5)*{}="4'";
   "4";"1" **\crv{(12,18)& (-6,18)};
   "1";"2" **\crv{(-6,-6)&(0,-6)}?(0)*\dir{>};
   "2";"3" **\crv{(0,6)&(6,6)};
   ?(.0)*\dir{>} ;
   "3";"3'" **\dir{-}
     ?(.1)*\dir{>};
   "4'";"4" **\dir{-}?(1)*\dir{>}
\endxy
\]
At this point one may wonder why we did not choose some other
isomorphism going from the identity to $x \ten \bar{x}$.  For
instance:
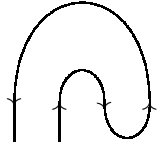
\[
\xy
 (6,-0)*{}="1";
 (0,0)*{}="2";
 (-6,0)*{}="3";
 (-6,-5)*{}="3'";
 (-12,0)*{}="4";
 (-12,-5)*{}="4'";
   "4";"1" **\crv{(-12,18)& (6,18)};
   "1";"2" **\crv{(6,-6)&(0,-6)}?(0)*\dir{<} ?(.97)*\dir{<} ;
   "2";"3" **\crv{(0,6)&(-6,6)};
   "3";"3'" **\dir{-}
     ?(.0)*\dir{<};
   "4'";"4" **\dir{-}?(1)*\dir{<}
\endxy
\]
is another morphism with the desired properties.  In fact, these
two morphisms are equal, as the following lemma shows.

In the calculations that follow, we denote an application of the
`horizontal slide' rule by a dashed line connecting the
appropriate zig and zag.  Dotted curves connecting two parallel
strings will indicate an application of the rules $e_x e_x^{-1} =
1_{\xb \tensor x}$ or $i_x^{-1}i_x = 1_{x \tensor \xb}$.
Furthermore, the rules $i_xi_x^{-1} = 1_1$ and $e_x^{-1}e_x = 1_1$
allow us to remove a closed loop any time one appears.

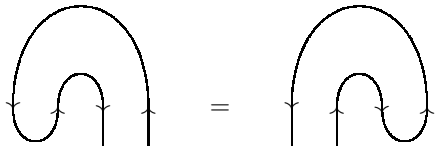
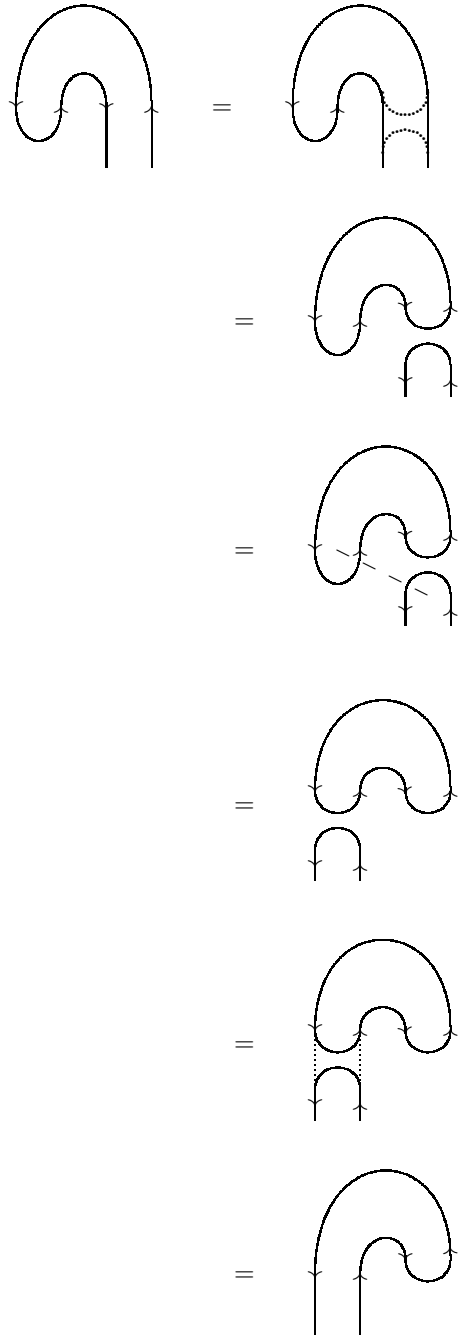
\begin{lem} \label{imove} \et
\[
\xy
 (-6,-0)*{}="1";
 (0,0)*{}="2";
 (6,0)*{}="3";
 (6,-5)*{}="3'";
 (12,0)*{}="4";
 (12,-5)*{}="4'";
   "4";"1" **\crv{(12,18)& (-6,18)};
   "1";"2" **\crv{(-6,-6)&(0,-6)}?(0)*\dir{>};
   "2";"3" **\crv{(0,6)&(6,6)};
   ?(.0)*\dir{>} ;
   "3";"3'" **\dir{-}
     ?(.1)*\dir{>};
   "4'";"4" **\dir{-}?(1)*\dir{>};
\endxy
   \qquad = \qquad
\xy
 (6,-0)*{}="1";
 (0,0)*{}="2";
 (-6,0)*{}="3";
 (-6,-5)*{}="3'";
 (-12,0)*{}="4";
 (-12,-5)*{}="4'";
   "4";"1" **\crv{(-12,18)& (6,18)};
   "1";"2" **\crv{(6,-6)&(0,-6)}?(0)*\dir{<} ?(.97)*\dir{<} ;
   "2";"3" **\crv{(0,6)&(-6,6)};
   "3";"3'" **\dir{-}
     ?(.0)*\dir{<};
   "4'";"4" **\dir{-}?(1)*\dir{<}
\endxy
\]
\end{lem}

\textbf{Proof. }
\[
\xy
 (-6,-0)*{}="1";
 (0,0)*{}="2";
 (6,0)*{}="3";
 (12,0)*{}="4";
 (12,-8)*{}="4'";
 (6,-8)*{}="0";
   "4";"1" **\crv{(12,18)& (-6,18)};
   "1";"2" **\crv{(-6,-6)&(0,-6)}?(0)*\dir{>};
   "2";"3" **\crv{(0,6)&(6,6)};
   ?(.0)*\dir{>} ;
   "3";"0" **\dir{-}
     ?(.03)*\dir{>};
   "4'";"4" **\dir{-}?(1)*\dir{>};
\endxy
\qquad = \qquad \xy
 (-6,-0)*{}="1";
 (0,0)*{}="2";
 (6,0)*{}="3";
 (12,0)*{}="4";
 (12,-8)*{}="4'";
 (6,-8)*{}="0";
   "4";"1" **\crv{(12,18)& (-6,18)};
   "1";"2" **\crv{(-6,-6) & (0,-6)};
     ?(.0)*\dir{>};
   "2";"3" **\crv{(0,6) & (6,6)};
     ?(.0)*\dir{>} ;
   "4'";"4" **\dir{-};
   "3";"0" **\dir{-};
 (6,2)*{}="A";
 (6,-6)*{}="A'";
 (12,2)="B";
 (12,-6)="B'";
 "A";"B" **\crv{~*=<2pt>{.}(6,-2)&(12,-2)};
   "A'";"B'" **\crv{~*=<2pt>{.}(6,-2)&(12,-2)};
\endxy
\]
\[
\xy
 (12,0)*{};
 (-12,0)*{};
\endxy
\qquad = \qquad \xy
 (-6,-0)*{}="1";
 (0,0)*{}="2";
 (6,0)*{}="3";
 (12,0)*{}="4";
 (12,-10)*{}="4'";
 (6,-10)*{}="0'";
 (6,2)*{}="A";
 (6,-6)*{}="A'";
 (12,2)="B";
 (12,-6)="B'";
   "B";"1" **\crv{(12,18)& (-6,18)}
     ?(.0)*\dir{>};
   "1";"2" **\crv{(-6,-6) & (0,-6)}
     ?(.0)*\dir{>};
   "2";"A" **\crv{(0,6) & (6,6)};
     ?(.0)*\dir{>} ;
   "4'";"B'" **\dir{-};
      ?(.5)*\dir{>};
   "A'";"0'" **\dir{-}; ?(.55)*\dir{>};
   "A";"B" **\crv{(6,-2)&(12,-2)}; ?(0.02)*\dir{>} ;
   "B'";"A'" **\crv{(12,-2)&(6,-2)};
\endxy
\]
\[
\xy
 (12,0)*{};
 (-12,0)*{};
\endxy
\qquad = \qquad \xy
 (-6,-0)*{}="1";
 (0,0)*{}="2";
 (6,0)*{}="3";
 (12,0)*{}="4";
 (12,-10)*{}="4'";
 (6,-10)*{}="0'";
 (6,2)*{}="A";
 (6,-6)*{}="A'";
 (12,2)="B";
 (12,-6)="B'";
   "B";"1" **\crv{(12,18)& (-6,18)}
     ?(.0)*\dir{>};
   "1";"2" **\crv{(-6,-6) & (0,-6)}
     ?(.0)*\dir{>};
   "2";"A" **\crv{(0,6) & (6,6)};
     ?(.0)*\dir{>} ;
   "4'";"B'" **\dir{-};
      ?(.5)*\dir{>};
   "A'";"0'" **\dir{-}; ?(.55)*\dir{>};
   "A";"B" **\crv{(6,-2)&(12,-2)}; ?(0.02)*\dir{>} ;
   "B'";"A'" **\crv{(12,-2)&(6,-2)};
        (-3,0)*{}="X";
 (9,-6)*{}="XX";
 "X";"XX" **\dir{--};
\endxy
\]

\[
\xy
 (12,0)*{};
 (-12,0)*{};
\endxy
\qquad = \qquad \xy
 (-6,2)*{}="1";
 (0,2)*{}="2";
 (6,0)*{}="3";
 (12,0)*{}="4";
 (0,-10)*{}="4'";
 (-6,-10)*{}="0'";
 (6,2)*{}="A";
 (-6,-6)*{}="A'";
 (12,2)="B";
 (0,-6)="B'";
   "B";"1" **\crv{(12,18)& (-6,18)}
     ?(.0)*\dir{>};
   "1";"2" **\crv{(-6,-2) & (0,-2)}
     ?(.0)*\dir{>};
   "2";"A" **\crv{(0,6) & (6,6)};
     ?(.0)*\dir{>} ;
   "4'";"B'" **\dir{-}
      ?(.53)*\dir{>};
   "A'";"0'" **\dir{-} ?(.5)*\dir{>};
   "A";"B" **\crv{(6,-2)&(12,-2)}
     ?(.03)*\dir{>} ;
   "B'";"A'" **\crv{(0,-2)&(-6,-2)};
\endxy
\]
\[
\xy
 (12,0)*{};
 (-12,0)*{};
\endxy
\qquad = \qquad \xy
 (-6,2)*{}="1";
 (0,2)*{}="2";
 (6,0)*{}="3";
 (12,0)*{}="4";
 (0,-10)*{}="4'";
 (-6,-10)*{}="0'";
 (6,2)*{}="A";
 (-6,-6)*{}="A'";
 (12,2)="B";
 (0,-6)="B'";
   "B";"1" **\crv{(12,18)& (-6,18)}
     ?(.0)*\dir{>};
   "1";"2" **\crv{(-6,-2) & (0,-2)}
     ?(.0)*\dir{>};
   "2";"A" **\crv{(0,6) & (6,6)};
     ?(.0)*\dir{>} ;
   "4'";"B'" **\dir{-}
      ?(.53)*\dir{>};
   "A'";"0'" **\dir{-} ?(.5)*\dir{>};
   "A";"B" **\crv{(6,-2)&(12,-2)}
     ?(.03)*\dir{>} ;
   "B'";"A'" **\crv{(0,-2)&(-6,-2)};
   "1";"A'" **\dir{.};
   "2";"B'" **\dir{.};
\endxy
\]
\[
\xy
 (-12,-0)*{};
 (12,0)*{};
\endxy
\qquad = \qquad \xy
 (-6,-0)*{}="1";
 (0,0)*{}="2";
 (0,-8)*{}="2'";
 (6,0)*{}="3";
 (-6,-8)*{}="0'";
 (6,2)*{}="A";
 (12,2)="B";
     "1";"0'" **\dir{-} ?(.05)*\dir{>};
     "B";"1" **\crv{(12,18)& (-6,18)};
       ?(.03)*\dir{>};
     "A";"B" **\crv{(6,-2)&(12,-2)};
       ?(.0)*\dir{>};
     "2";"A" **\crv{(0,6) & (6,6)};
       ?(.0)*\dir{>};
     "2'";"2" **\dir{-};
\endxy
\]
\hskip 30em \qed

Now let us show that $(x, \xb, i'_x, e'_x)$ satisfies the zig-zag
identities. Recall that these identities say that:
\[
\xy0;/r.28pc/: (-6,-8)*{}="1E"; (-6,0)*{}="1"; (0,0)*{}="2";
(6,0)*{}="3"; (6,8)*{}="3B";
 "2";"1" **\crv{(0,10)& (-6,10)}
     ?(.02)*\dir{>}  ?(1)*\dir{>};
 "3";"2" **\crv{(6,-10)& (0,-10)} ;
 "1";"1E" **\dir{-};
 "3B";"3" **\dir{-}?(1)*\dir{>} ;
(-2.5,9)*{\scriptstyle i'_x}; (3.5,-9)*{\scriptstyle e'_x};
\endxy
\qquad = \xy0;/r.28pc/: (-6,8)*{}; (0,8)*{}; (0,-8)*{}; **\dir{-}
?(.5)*\dir{<}; (6,-8)*{}?(.5)*\dir{<}+(3,1)*{};
\endxy
\]
and
\[
\xy0;/r.28pc/: (-6,8)*{}="1E"; (-6,0)*{}="1"; (0,0)*{}="2";
(6,0)*{}="3"; (6,-8)*{}="3B";
 "2";"1" **\crv{(0,-10)& (-6,-10)}
     ?(.02)*\dir{>}  ?(1)*\dir{>};
 "3";"2" **\crv{(6,10)& (0,10)}
       ;
 "1";"1E" **\dir{-};
 "3B";"3" **\dir{-}?(1)*\dir{>};
(-2.5,-9)*{\scriptstyle e'_x}; (3.5,9)*{\scriptstyle i'_x};
\endxy
\qquad = \xy0;/r.28pc/: (-6,8)*{}; (0,-8)*{}; (0,8)*{}; **\dir{-}
?(.5)*\dir{<}; (6,-8)*{}?(.5)*\dir{<}+(3,1)*{};
\endxy
\]
If we express $i'_x$ and $e'_x$ in terms of $i_x$ and $e_x$, these
equations become
\[
\vcenter{ \xy
 (-6,-0)*{}="1";
 (0,0)  *{}="2";
 (6,0)  *{}="3";
 (6,-5) *{}="3'";
 (12,0) *{}="4";
 (18,0) *={}="5";
 (18,12)*={}="5'";
   "5";"4" **\crv{(18,-6)& (12,-6)};
   "5'";"5" **\dir{-}            ?(1)*\dir{>};
   "4";"1" **\crv{(12,18)& (-6,18)};
                                ?(.0)*\dir{>}  ?(1)*\dir{>};
   "1";"2" **\crv{(-6,-6)&(0,-6)};
   "2";"3" **\crv{(0,6)&(6,6)};
                                ?(.0)*\dir{>} ;
   "3";"3'" **\dir{-}
                                 ?(0)*\dir{>};
   (16,-7)*{\scriptstyle e_{x}};        
   (4,7)*{\scriptstyle e_{x}^{-1}};     
   (4,16)*{\scriptstyle i_{x}};         
   (-2,-7)*{\scriptstyle i_{x}^{-1}};
\endxy}
\qquad = \qquad \xy (0,8)*{};(0,-10)**\dir{-} ?(.5)*\dir{>};
\endxy
\]
and
\[
\vcenter{ \xy
 (-6,-0)*{}="1";
 (0,0)*{}="2";
 (6,0)*{}="3";
 (12,0)*{}="4";
 (12,-12)*{}="4'";
 (-12,12)*{}="0'";
 (-12,0)*{}="0";
   "4";"1" **\crv{(12,18)& (-6,18)};
     ?(.0)*\dir{>}  ;
   "1";"2" **\crv{(-6,-6) & (0,-6)}
     ?(.0)*\dir{>};
   "2";"3" **\crv{(0,6) & (6,6)};
     ?(.0)*\dir{>} ;
   "4'";"4" **\dir{-};
   "3";"0" **\crv{(6,-18)&(-12,-18)}
     ?(.0)*\dir{>} ;
   "0";"0'" **\dir{-} ?(.05)*\dir{>};
         (3,16)*{\scriptstyle i_{x}};
         (-3,-16)*{\scriptstyle e_{x}};
         (-3,-7)*{\scriptstyle i_{x}^{-1}};
         (4,7)*{\scriptstyle e_{x}^{-1}};
\endxy}
\qquad = \qquad \xy (0,12)*{};(0,-12)**\dir{-} ?(.5)*\dir{<};
\endxy
\]
To verify these two equations we use string diagrams.  The first
equation can be shown as follows:
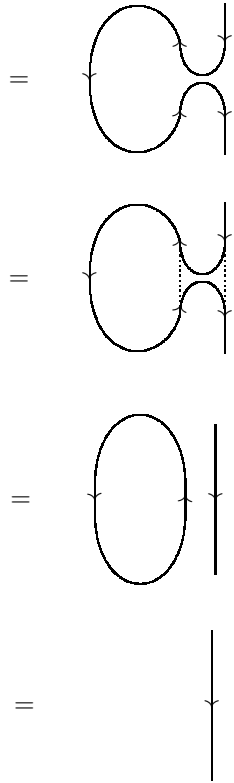
\[
\xy
 (-6,-0)*{}="1";
 (0,0)  *{}="2";
 (6,0)  *{}="3";
 (6,-5) *{}="3'";
 (12,0) *{}="4";
 (18,0) *={}="5";
 (18,12)*={}="5'";
   "5";"4" **\crv{(18,-6)& (12,-6)};
   "5'";"5" **\dir{-} ?(1)*\dir{>};
   "4";"1" **\crv{(12,18)& (-6,18)};
     ?(.0)*\dir{>}  ?(1)*\dir{>};
   "1";"2" **\crv{(-6,-6)&(0,-6)};
   "2";"3" **\crv{(0,6)&(6,6)};
     ?(.0)*\dir{>} ;
   "3";"3'" **\dir{-}
     ?(0)*\dir{>};
\endxy
\qquad = \qquad \xy
 (-6,-0)*{}="1";
 (0,0)*{}="2";
 (6,0)*{}="3";
 (6,-5)*{}="3'";
 (12,0)*{}="4";
 (18,0)*={}="5";
 (18,12)*={}="5'";
   "5";"4" **\crv{(18,-6)& (12,-6)};
   "5'";"5" **\dir{-} ?(1)*\dir{>};
   "4";"1" **\crv{(12,18)& (-6,18)}
     ?(.)*\dir{>}  ?(1)*\dir{>};
   "1";"2" **\crv{(-6,-6)&(0,-6)};
   "2";"3" **\crv{(0,6)&(6,6)};
     ?(.0)*\dir{>} ;
   "3";"3'" **\dir{-}
     ?(0)*\dir{>} ;
 (3,0)*{}="A";
 (15,0)*{}="B";
     "A";"B" **\dir{--};
\endxy
\]
\[
\xy
 (-6,-0)*{}="1";
 (0,0)  *{}="2";
 (6,0)  *{}="3";
 (6,-5) *{}="3'";
 (12,0) *{}="4";
 (18,0) *={}="5";
 (18,12)*={}="5'";
\endxy
\qquad = \qquad \xy
 (-6,-0)*{}="1";
 (6,-5)  *{}="2";
 (12,-5)  *{}="3";
 (12,-10) *{}="3'";
 (6,5) *{}="4";
 (12,5) *={}="5";
 (12,10)*={}="5'";
 (18,0) *={}="";
   "5";"4" **\crv{(12,-1)& (6,-1)}
                                ?(1)*\dir{>};
   "5'";"5" **\dir{-}           ?(1)*\dir{>};
   "4";"1" **\crv{(4,12)& (-6,12)};
                                ?(1)*\dir{>};
   "1";"2" **\crv{(-6,-12)&(4,-12)};
   "2";"3" **\crv{(6,1)&(12,1)};
                                ?(.04)*\dir{>} ;
   "3";"3'" **\dir{-}
                                ?(0)*\dir{>};
\endxy
\]
\[
\xy
 (-6,-0)*{}="";
 (18,0) *={}="";
\endxy
\qquad = \qquad \xy
 (-6,-0)*{}="1";
 (6,-5)  *{}="2";
 (12,-5)  *{}="3";
 (12,-10) *{}="3'";
 (6,5) *{}="4";
 (12,5) *={}="5";
 (12,10)*={}="5'";
 (18,0) *={}="";
   "5";"4" **\crv{(12,-1)& (6,-1)}
         ?(1)*\dir{>};
   "5'";"5" **\dir{-} ?(1)*\dir{>};
   "4";"1" **\crv{(4,12)& (-6,12)};
       ?(1)*\dir{>};
   "1";"2" **\crv{(-6,-12)&(4,-12)};
   "2";"3" **\crv{(6,1)&(12,1)};
     ?(.08)*\dir{>} ;
   "3";"3'" **\dir{-}
     ?(0)*\dir{>};
   "2";"4" **\dir{.};
   "3";"5" **\dir{.};
\endxy
\]
\[
\xy
 (-6,-0)*{}="";
 (18,0) *={}="";
\endxy
\qquad = \qquad \xy
 (6,0)*{}="1";
 (-6,0)*{}="2";
   "1";"2" **\crv{(7,15)&(-7,15)}
     ?(0.01)*\dir{>}  ?(1)*\dir{>} ;
   "2";"1" **\crv{ (-7,-15)& (7,-15)};
\endxy
\quad \xy (0,10)*{};(0,-10)**\dir{-} ?(.5)*\dir{>};
\endxy \qquad
\]
\[
\xy
 (-6,-0)*{}="1";
 (18,0) *={}="5";
\endxy
\qquad = \qquad \xy
 (6,0)*{}="";
 (-6,0)*{}="";
\endxy
\quad \xy (0,10)*{};(0,-10)**\dir{-} ?(.5)*\dir{>};
\endxy \qquad
\]
The second equation can be shown with the help of Lemma
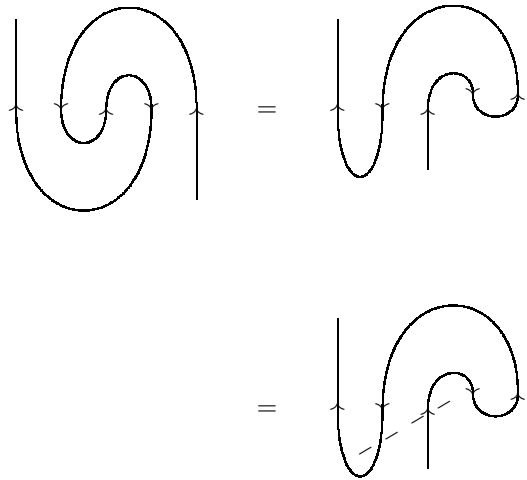
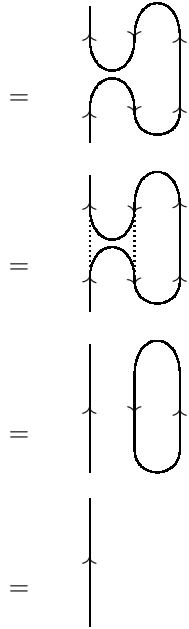
\ref{imove}:
\[
\xy
 (-6,-0)*{}="1";
 (0,0)*{}="2";
 (6,0)*{}="3";
 (12,0)*{}="4";
 (12,-12)*{}="4'";
 (-12,12)*{}="0'";
 (-12,0)*{}="0";
   "4";"1" **\crv{(12,18)& (-6,18)};
     ?(.0)*\dir{>}  ;
   "1";"2" **\crv{(-6,-6) & (0,-6)}
     ?(.0)*\dir{>};
   "2";"3" **\crv{(0,6) & (6,6)};
     ?(.0)*\dir{>} ;
   "4'";"4" **\dir{-};
   "3";"0" **\crv{(6,-18)&(-12,-18)}
     ?(.0)*\dir{>} ;
   "0";"0'" **\dir{-} ?(.05)*\dir{>};
\endxy
\qquad  = \qquad \xy
 (-6,-0)*{}="1";
 (0,0)*{}="2";
 (0,-8)*{}="2'";
 (6,0)*{}="3";
 (-12,12)*{}="0'";
 (-12,0)*{}="0";
 (6,2)*{}="A";
 (12,2)="B";
     "0";"0'" **\dir{-} ?(.05)*\dir{>};
     "1";"0" **\crv{(-6,-12) &(-12,-12)}
         ?(.0)*\dir{>};
     "B";"1" **\crv{(12,18)& (-6,18)};
       ?(.0)*\dir{>};
     "A";"B" **\crv{(6,-2)&(12,-2)};
       ?(.0)*\dir{>};
     "2";"A" **\crv{(0,6) & (6,6)};
       ?(.0)*\dir{>};
     "2'";"2" **\dir{-};
\endxy
\]
\[
\xy
 (-12,-0)*{};
 (12,0)*{};
\endxy
\qquad = \qquad \xy
 (-6,-0)*{}="1";
 (0,0)*{}="2";
 (0,-8)*{}="2'";
 (6,0)*{}="3";
 (-12,12)*{}="0'";
 (-12,0)*{}="0";
 (6,2)*{}="A";
 (12,2)="B";
     "0";"0'" **\dir{-} ?(.05)*\dir{>};
     "1";"0" **\crv{(-6,-12) &(-12,-12)}
         ?(.0)*\dir{>};
     "B";"1" **\crv{(12,18)& (-6,18)};
       ?(.0)*\dir{>};
     "A";"B" **\crv{(6,-2)&(12,-2)};
       ?(.0)*\dir{>};
     "2";"A" **\crv{(0,6) & (6,6)};
       ?(.0)*\dir{>};
     "2'";"2" **\dir{-};
 (-9,-6)*{}="X";
 (3,1)*{}="XX";
 "X";"XX" **\dir{--};
\endxy
\]
\[
\xy
 (-12,-0)*{};
 (12,0)*{};
\endxy
\qquad = \qquad \xy
 (-12,12)*{}="0'";
 (-12,8)*{}="0";
 (-6,8)*{}="1";
 (-12,-2)*{}="2";
 (-12,-6)*{}="2'";
 (6,-2)*{}="3";
 (-6,-2)*{}="A";
 (0,-2)*{}="B";
 (0,8)*{}="B'";
 (12,0)*{}="";
     "0";"0'" **\dir{-} ?(.05)*\dir{>};
     "1";"0" **\crv{(-6,2) &(-12,2)}
         ?(.0)*\dir{>};
     "B'";"1" **\crv{(0,14)& (-6,14)};
       ?(.0)*\dir{>};
     "A";"B" **\crv{(-6,-6)&(0,-6)};
       ?(.0)*\dir{>};
     "2";"A" **\crv{(-12,4) & (-6,4)};
       ?(.0)*\dir{>};
     "2'";"2" **\dir{-};
     "B";"B'" **\dir{-} ?(.06)*\dir{>};;
\endxy
\]
\[
\xy
 (-12,-0)*{};
 (12,0)*{};
\endxy
\qquad = \qquad \xy
 (-12,12)*{}="0'";
 (-12,8)*{}="0";
 (-6,8)*{}="1";
 (-12,-2)*{}="2";
 (-12,-6)*{}="2'";
 (6,-2)*{}="3";
 (-6,-2)*{}="A";
 (0,-2)*{}="B";
 (0,8)*{}="B'";
 (12,0)*{}="";
     "0";"0'" **\dir{-} ?(.05)*\dir{>};
     "1";"0" **\crv{(-6,2) &(-12,2)}
         ?(.05)*\dir{>};
     "B'";"1" **\crv{(0,14)& (-6,14)};
       ?(.03)*\dir{>};
     "A";"B" **\crv{(-6,-6)&(0,-6)};
       ?(.04)*\dir{>};
     "2";"A" **\crv{(-12,4) & (-6,4)};
       ?(.03)*\dir{>};
     "2'";"2" **\dir{-};
     "B";"B'" **\dir{-} ?(.06)*\dir{>};
     "0";"2" **\dir{.};
     "1";"A" **\dir{.};
\endxy
\]
\[
\xy
 (-12,-0)*{};
 (12,0)*{};
\endxy
\qquad = \qquad \xy
 (-12,12)*{}="0'";
 (-6,8)*{}="1";
 (-12,-5)*{}="2'";
 (6,-2)*{}="3";
 (-6,-2)*{}="A";
 (0,-2)*{}="B";
 (0,8)*{}="B'";
 (12,0)*{}="";
     "B'";"1" **\crv{(0,14)& (-6,14)};
     "A";"B" **\crv{(-6,-6)&(0,-6)};
     "B";"B'" **\dir{-} ?(.53)*\dir{>};
     "2'";"0'" **\dir{-} ?(.5)*\dir{>};
     "1";"A" **\dir{-};  ?(.5)*\dir{>};
\endxy
\]
\[
\xy
 (-12,-0)*{};
 (12,0)*{};
\endxy
\qquad = \qquad \xy
 (-12,12)*{}="0'";
 (-12,-5)*{}="2'";
 (12,0)*{}="";
     "2'";"0'" **\dir{-} ?(.53)*\dir{>};
\endxy
\]
\hskip 30em \qed

The `improvement' process of Theorem \ref{improve} can be made
into a 2-functor $\imp \maps \wg \to \cg$:

\begin{cor} \et
There exists a 2-functor $ \imp \maps \wg \to \cg$ which sends any
object $C \in \wg$ to $\imp(C) \in \cg$ and acts as the identity
on morphisms and 2-morphisms.
\end{cor}

\textbf{Proof. } This is a trivial consequence of Theorem
\ref{improve}.  Obviously all domains, codomains, identities and
composites are preserved, since the 1-morphisms and 2-morphisms
are unchanged as a result of Definitions \ref{coherenthomo} and
\ref{coherent2homo}. \qed

On the other hand, there is also a forgetful 2-functor ${\rm F}
\maps \cg \to \wg$, which forgets the extra structure on objects
and acts as the identity on morphisms and 2-morphisms.

\begin{thm} \et
The 2-functors $ \imp \maps \wg \to \cg$, ${\rm F} \maps \cg \to
\wg$ extend to define a 2-equivalence between the 2-categories
$\wg$ and $\cg$.
\end{thm}

\textbf{Proof. } The 2-functor $\imp$ equips each object of $\wg$
with the structure of a coherent 2-group, while $\rm F$ forgets
this extra structure.  Both act as the identity on morphisms and
2-morphisms.  As a consequence, $\imp$ followed by $\rm F$ acts as
the identity on $\wg$:
\[   \imp \circ {\rm F} = 1_{\wg} \]
(where we write the functors in order of application). To prove
the theorem, it therefore suffices to construct a natural
isomorphism
\[    e \maps {\rm F} \circ \imp \To 1_{\cg} .\]

To do this, note that applying $F$ and then $\imp$ to a coherent
2-group $C$ amounts to forgetting the choice of adjoint
equivalence for each object of $C$ and then making a new such
choice.  We obtain a new coherent 2-group $\imp({\rm F}(C))$, but
it has the same underlying weak monoidal category, so the identity
functor on $C$ defines a coherent 2-group isomorphism from
$\imp({\rm F}(C))$ to $C$.   We take this as $e_C \maps \imp({\rm
F}(C)) \to C$.

To see that this defines a natural isomorphism between 2-functors,
note that for every coherent 2-group homomorphism $f \maps C \to
C'$ we have a commutative square:
\[
\xymatrix{ \imp(F(C)) \ar[rr]^{\imp(F(f))} \ar[d]_{e_C}
&& \imp(F(C')) \ar[d]^{e_{C'}} \\
C \ar[rr]^{f} && C' }
\]
This commutes because $\imp(F(f)) = f$ as weak monoidal functors,
while $e_C$ and $e_{C'}$ are the identity as weak monoidal
functors. \qed

The significance of this theorem is that while we have been
carefully distinguishing between weak and coherent 2-groups, the
difference is really not so great.  Since the 2-category of weak
2-groups is 2-equivalent to the 2-category of coherent ones, one
can use whichever sort of 2-group happens to be more convenient at
the time, freely translating results back and forth as desired.
So, except when one is trying to be precise, one can relax and use
the term {\bf 2-group} for either sort.

Of course, we made heavy use of the axiom of choice in proving the
existence of the improvement 2-functor $\imp \maps \wg \to \cg$,
so constructivists will not consider weak and coherent 2-groups to
be equivalent.  Mathematicians of this ilk are urged to use
coherent 2-groups.  Indeed, even pro-choice mathematicians will
find it preferable to use coherent 2-groups when working in
contexts where the axiom of choice fails.  These are not at all
exotic.   For example, the theory of `Lie 2-groups' works well
with coherent 2-groups, but not very well with weak 2-groups, as
we shall see in Section \ref{internalizationsection}.

To conclude, let us summarize {\it why} weak and coherent 2-groups
are not really so different.  At first, the choice of a specified
adjoint equivalence for each object seems like a substantial extra
structure to put on a weak 2-group.  However, Theorem
\ref{improve} shows that we can always succeed in putting this
extra structure on any weak 2-group. Furthermore, while there are
many ways to equip a weak 2-group with this extra structure, there
is `essentially' just one way, since all coherent 2-groups with
the same underlying weak 2-group are isomorphic. It is thus an
example of what Kelly and Lack \cite{KL} call a `property-like
structure'.

Of course, the observant reader will note that this fact has
simply been built into our definitions!  The reason all coherent
2-groups with the same underlying weak 2-group are isomorphic is
that we have defined a homomorphism of coherent 2-groups to be a
weak monoidal functor, not requiring it to preserve the choice of
adjoint equivalence for each object.  This may seem like
`cheating', but in the next section we justify it by showing that
this choice is {\it automatically preserved up to coherent
isomorphism} by any weak monoidal functor.

\section{Preservation of weak inverses} \label{preservationsection}

Suppose that $F \maps C \to C'$ is a weak monoidal functor between
coherent 2-groups.  To show that $F$ automatically preserves the
specified weak inverses up to isomorphism, we now construct an
isomorphism
\[    (F_{-1})_x \maps \overline{F(x)} \to F(\xb)  \]
for each object $x \in C$.  This isomorphism is uniquely
determined if we require the following coherence laws:
\begin{description}
 \item[H1]
\[
\xy (-35,8)*+{F(x) \ten \overline{F(x)}}="TL"; (0,8)*+{F(x) \ten
F(\xb)}="TM"; (30,8)*+{F(x \ten \xb)}="TR"; (-35,-8)*+{1'}="BL";
(30,-8)*+{F(1)}="BR";
       {\ar^{1 \ten F_{-1}} "TL";"TM"};
       {\ar^{i_{F(x)}} "BL";"TL"};
       {\ar^{F_{2}} "TM";"TR"};
       {\ar_{F(i_x)} "BR";"TR"};
       {\ar_{F_{0}} "BL";"BR"};
\endxy
\]
 \item[H2]
\[
\xy (-35,8)*+{\overline{F(x)} \ten F(x)}="TL"; (0,8)*+{F(\xb) \ten
F(x)}="TM"; (30,8)*+{F(\xb \ten x)}="TR"; (-35,-8)*+{1'}="BL";
(30,-8)*+{F(1)}="BR";
       {\ar^{F_{-1} \ten 1} "TL";"TM"};
       {\ar_{e_{F(x)}} "TL";"BL"};
       {\ar^{F_{2}} "TM";"TR"};
       {\ar^{F(e_x)} "TR";"BR"};
       {\ar_{F_{0}} "BL";"BR"};
\endxy
\]
\end{description}
These say that $F_{-1}$ is compatible with units and counits. In
the above diagrams and in what follows, we suppress the subscript
on $F_{-1}$, just as we are already doing for $F_2$.

\begin{thm} \label{preservation} \et  Suppose that $F \maps C \to C'$ is a
homomorphism of coherent 2-groups.  Then for any object $x \in C$
there exists a unique isomorphism $F_{-1} \maps \overline{F(x)}
\to F(\xb)$ satisfying the coherence laws {\bf H1} and {\bf H2}.
\end{thm}

\textbf{Proof. } This follows from the general fact that pseudofunctors
between bicategories preserve adjunctions.  However, to illustrate
the use of string diagrams we prefer to simply take one of these laws, 
solve it for $F_{-1}$, and show that the result also satisfies the other 
law.  We start by writing the law \textbf{H1} in a more suggestive manner:
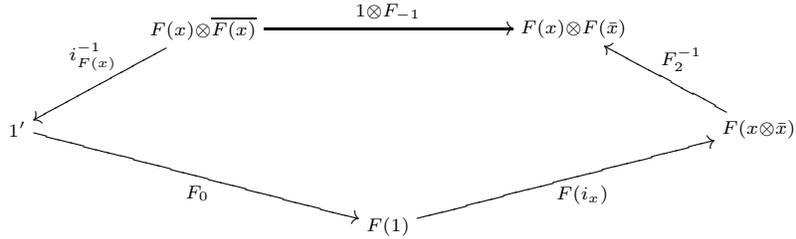
\[
\def\objectstyle{\scriptstyle}
\def\labelstyle{\scriptstyle}
\vcenter{ \xymatrix@!C{
 &  F(x) \ten \overline{F(x)}
   \ar[rr]^{1 \ten F_{-1}}
   \ar[dl]_{i_{F(x)}^{-1}}
 && F(x) \ten F(\xb) \\
 1' \ar[drr]_{F_{0}}
 &&&& F(x \ten \xb) \ar[ul]_{F_{2}^{-1}} \\
 && F(1) \ar[urr]_{F(i_{x})}
 }}
\]
If we assume this diagram commutes, it gives a formula for
\[   1 \tensor F_{-1}
  \maps \xymatrix@1{F(x) \ten \overline{F(x)}
         \ar[r]^<<<<<{\sim} & \;     F(x) \ten F(\xb) .  }
\]
Writing this formula in string notation, it becomes
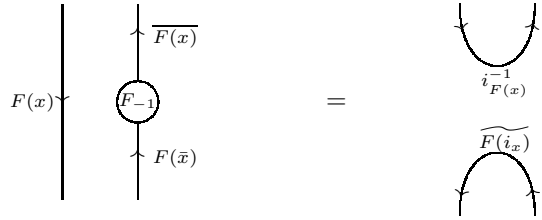
\[
\xy (-15,0)*{};(15,0)*{}; (-5,13)*{}; (-5,-13)*{};
**\dir{-}?(.5)*\dir{<}+(-4,0)*{\scriptstyle F(x)};
(5,0)*\xycircle(2.75,2.75){-}="f"; (5,13)**\dir{-}
?(.75)*\dir{>}+(5,-1)*{\scriptstyle \overline{F(x)}};
"f";(5,-13)**\dir{-} ?(.5)*\dir{<}+(5,-1)*{\scriptstyle F(\xb)};
(5,0)*{\scriptstyle F_{-1}};
\endxy
\qquad \qquad = \qquad \qquad \xy
 (-5,13)*{};
 (5,13)*{}
   **\crv{(-5,2) & (5,2)} ?(.12)*\dir{>} ?(.92)*\dir{>};
 (5,-15)*{};
 (-5,-15)*{}
   **\crv{(5,-4) & (-5,-4)} ?(.12)*\dir{>} ?(.92)*\dir{>};
(1,2.5)*{\scriptstyle i_{F(x)}^{-1}};
     (1,-5)*{\scriptstyle \fix};
\endxy
\]
where we set
\[   \fix = F_0 \circ F(i_x) \circ F_2^{-1} \maps 1' \to F(x) \ten F(\xb). \]
This equation can in turn be solved for $F_{-1}$, as follows:
\[
\xy (0,0)*{\scriptstyle F_{-1}}; (-10,0)*{};(10,0)*{};
(0,0)*\xycircle(2.75,2.75){-}="f"; (0,14)**\dir{-}
?(.75)*\dir{>}+(5,-1)*{\scriptstyle \overline{F(x)}};
"f";(0,-14)**\dir{-} ?(.5)*\dir{<}+(5,-1)*{\scriptstyle F(\xb)};
\endxy
\qquad = \qquad \xy (-15,0)*{};(15,0)*{}; (-8,10)*{}="TL";
(0,10)*{}="TR"; (-8,-10)*{}="BL"; (0,-10)*{}="BR";
     "TL";"BL" **\dir{-};
         ?(.5)*\dir{<};
     "TR";"BR" **\dir{-};
         ?(.5)*\dir{>};
     "TL";"TR" **\crv{(-8,18)& (0,18)};
     "BL";"BR" **\crv{(-8,-18)& (0,-18)};
(8,0)*{\scriptstyle F_{-1}}; (8,0)*\xycircle(2.75,2.75){-}="f";
(8,14)**\dir{-} ?(.75)*\dir{>}+(5,-1)*{\scriptstyle
\overline{F(x)}}; "f";(8,-14)**\dir{-}
?(.5)*\dir{<}+(5,-1)*{\scriptstyle F(\xb)};
 (3.5,5)*{\scriptstyle F(x)};
 (-12,5)*{\scriptstyle \overline{F(x)}};
\endxy
\]
\[
\xy (-10,0)*{};(10,0)*{};
\endxy
\qquad = \qquad \xy (-15,0)*{};(15,0)*{}; (-8,10)*{}="TL";
(0,10)*{}="TR"; (-8,-10)*{}="BL"; (0,-10)*{}="BR"; (8,10)*{}="RU";
(8,-10)*{}="RB";
     "TL";"BL" **\dir{-};
         ?(.5)*\dir{<};
     "TR";"BR" **\dir{-};
         ?(.5)*\dir{>};
     "TL";"TR" **\crv{(-8,18)& (0,18)};
     "BL";"BR" **\crv{(-8,-18)& (0,-18)};
     "RU";"TR" **\crv{~*=<2pt>{.}(8,4)& (0,4)};
     "RB";"BR" **\crv{~*=<2pt>{.}(8,-4)& (0,-4)};
(8,0)*{\scriptstyle F_{-1}}; (8,0)*\xycircle(2.75,2.75){-}="f";
(8,14)**\dir{-} ?(.75)*\dir{>}+(5,-1)*{\scriptstyle
\overline{F(x)}}; "f";(8,-14)**\dir{-}
?(.5)*\dir{<}+(5,-1)*{\scriptstyle F(\xb)};
 (-12,5)*{\scriptstyle \overline{F(x)}};
\endxy
\]
\[
\xy (-10,0)*{};(10,0)*{};
\endxy
\qquad = \qquad \xy (-15,0)*{};(15,0)*{}; (-8,10)*{}="TL";
(0,10)*{}="TR"; (-8,-10)*{}="BL"; (0,-10)*{}="BR"; (8,10)*{}="RU";
(8,16)*{}="RUEND"; (8,-16)*{}="RBEND"; (8,-10)*{}="RB";
     "TL";"BL" **\dir{-};
         ?(.5)*\dir{<};
     "RU";"RUEND" **\dir{-};
         ?(.15)*\dir{>}+(5,2)*{\scriptstyle \overline{F(x)}};
     "RB";"RBEND" **\dir{-};
         ?(0)*\dir{<}+(5,-2)*{\scriptstyle F(\xb)};
     "TL";"TR" **\crv{(-8,18)& (0,18)};
         ?(1)*\dir{>};
     "BL";"BR" **\crv{(-8,-18)& (0,-18)};
         ?(.95)*\dir{<}+(4,8)*{\scriptstyle \fix} ;;
     "RU";"TR" **\crv{(8,4)& (0,4)};
     "RB";"BR" **\crv{(8,-4)& (0,-4)};
         (-12,5)*{\scriptstyle \overline{F(x)}};
\endxy
\]
\[
\xy (-10,0)*{};(10,0)*{};
\endxy
\qquad = \qquad \xy (-15,0)*{};(15,0)*{}; (-8,10)*{}="TL";
(0,10)*{}="TR"; (-8,-10)*{}="BL"; (0,-10)*{}="BR"; (8,10)*{}="RU";
(8,16)*{}="RUEND"; (8,-16)*{}="RBEND"; (8,-10)*{}="RB";
     "TL";"BL" **\dir{-};
         ?(.5)*\dir{<};
     "RU";"RUEND" **\dir{-};
         ?(.15)*\dir{>}+(5,2)*{\scriptstyle \overline{F(x)}};
     "RB";"RBEND" **\dir{-};
         ?(0)*\dir{<}+(5,-2)*{\scriptstyle F(\xb)};
     "TL";"TR" **\crv{(-8,18)& (0,18)};
         ?(1)*\dir{>};
     "BL";"BR" **\crv{(-8,-18)& (0,-18)};
         ?(.95)*\dir{<}+(4,8)*{\scriptstyle \fix} ;;
     "RU";"TR" **\crv{(8,4)& (0,4)};
     "RB";"BR" **\crv{(8,-4)& (0,-4)};
             (16,20)*{}="X1";
             (-14,20)*{}="X2";
             (-14,3)*{}="X3";
             (16,3)*{}="X4";
               "X1";"X2" **\dir{.};
               "X2";"X3" **\dir{.};
               "X3";"X4" **\dir{.};
               "X4";"X1" **\dir{.};
               (-12,0)*{\scriptstyle \overline{F(x)}};
\endxy
\]
\[
\xy (-10,0)*{};(10,0)*{};
\endxy
\qquad = \qquad \xy (-15,0)*{};(15,0)*{}; (-8,10)*{}="TL";
(-8,-10)*{}="BL"; (0,-10)*{}="BR"; (8,-16)*{}="RBEND";
(8,-10)*{}="RB";
     "TL";"BL" **\dir{-};
         ?(.5)*\dir{<}+(5,8)*{\scriptstyle \overline{F(x)}};
     "RB";"RBEND" **\dir{-};
         ?(0)*\dir{<}+(5,-2)*{\scriptstyle F(\xb)};
     "BL";"BR" **\crv{(-8,-18)& (0,-18)};
         ?(.95)*\dir{<}+(4,8)*{\scriptstyle \fix} ;;
     "RB";"BR" **\crv{(8,-4)& (0,-4)};
\endxy
\]
Here and in the arguments to come we omit the labels $i_{F(x)}$,
$e_{F(x)}$, $i_{F(x)}^{-1}$, $e_{F(x)}^{-1}$.

Since we have solved for $F_{-1}$ starting from \textbf{H1}, we
have already shown the morphism satisfying this law is unique. We
also know it is an isomorphism, since all morphisms in $C'$ are
invertible. However, we should check that it exists --- that is,
it really does satisfy this coherence law. The proof is a string
diagram calculation:
\[
\xy (-15,0)*{};(15,0)*{}; (-5,12)*{}; (-5,-12)*{};
**\dir{-}?(.5)*\dir{<}+(-4,0)*{\scriptstyle F(x)};
(5,0)*{\scriptstyle F_{-1}}; (5,0)*\xycircle(2.75,2.75){-}="f";
(5,12)**\dir{-} ?(.75)*\dir{>}+(5,-1)*{\scriptstyle
\overline{F(x)}}; "f";(5,-12)**\dir{-}
?(.5)*\dir{<}+(5,-1)*{\scriptstyle F(\xb)};
\endxy
\qquad = \qquad \xy
 (-10,12)*{}; (-10,-12)*{}; **\dir{-}?(.5)*\dir{<}+(-4,0)*{\scriptstyle F(x)};
 (-2,12)*{}="1";
 (6,-7)  *{}="2";
 (14,-7) *{}="3";
 (14,-13)*{}="3'";
 (18,0) *={}="";
   "1";"2" **\crv{(-4,-21)&(5,-15)};
     ?(.15)*\dir{<}+(5,10)*{\scriptstyle \overline{F(x)}};
   "2";"3" **\crv{(6,1)&(14,1)};
     ?(0)*\dir{<}+(4,8)*{\scriptstyle \fix} ;
     ?(.95)*\dir{<}+(5,-2)*{\scriptstyle F(\xb)};
   "3";"3'" **\dir{-}

\endxy
\]
\[
\xy (-15,0)*{};(15,0)*{};
\endxy
\qquad = \qquad \xy
 (-10,12)*{}; (-10,-12)*{}; **\dir{-}?(.5)*\dir{<}+(-4,0)*{\scriptstyle F(x)};
 (-2,12)*{}="1";
 (6,-7)  *{}="2";
 (14,-7) *{}="3";
 (14,-13)*{}="3'";
 (18,0) *={}="";
   "1";"2" **\crv{(-4,-21)&(5,-15)};
     ?(.15)*\dir{<}+(5,10)*{\scriptstyle \overline{F(x)}};
   "2";"3" **\crv{(6,1)&(14,1)};
     ?(0)*\dir{<}+(4,8)*{\scriptstyle \fix} ;
     ?(.95)*\dir{<}+(5,-2)*{\scriptstyle F(\xb)};
   "3";"3'" **\dir{-};
 (-2,6)*{}="X1";
 (-10,6)*{}="X2";
     "X1";"X2" **\crv{~*=<2pt>{.}(-2,-1)&(-10,-1)};
 (-2,-6)*{}="XA";
 (-10,-6)*{}="XB";
     "XA";"XB" **\crv{~*=<2pt>{.}(-2,1)&(-10,1)};
\endxy
\]
\[
\xy (-15,0)*{};(15,0)*{};
\endxy
\qquad = \qquad \xy
 (-10,12)*{}="0";
 (-2,12)*{}="1";
 (6,-7)  *{}="2";
 (14,-7) *{}="3";
 (14,-13)*{}="3'";
 (18,0) *={}="";
 (-2,6)*{}="X1";
 (-10,6)*{}="X2";
 (-2,-6)*{}="XA";
 (-10,-6)*{}="XB";
 (-10,-12)*{}="XB'";
     "XB";"XB'" **\dir{-}; ?(.4)*\dir{>};
     "X1";"1" **\dir{-}; ?(0)*\dir{>};
     "X2";"0" **\dir{-}; ?(0)*\dir{<};
   "XA";"2" **\crv{(-1,-15)&(6,-15)};
     ?(0)*\dir{<};
   "2";"3" **\crv{(6,1)&(14,1)};
     ?(0)*\dir{<}+(4,8)*{\scriptstyle \fix} ;
     ?(.95)*\dir{<}+(5,-2)*{\scriptstyle F(\xb)};
   "3";"3'" **\dir{-};
     "X1";"X2" **\crv{(-2,-1)&(-10,-1)};
     "XA";"XB" **\crv{(-2,1)&(-10,1)};
\endxy
\]
\[
\xy (-15,0)*{};(15,0)*{};
\endxy
\qquad = \qquad \xy
 (-10,12)*{}="0";
 (-2,12)*{}="1";
 (6,-7)  *{}="2";
 (14,-7) *{}="3";
 (14,-13)*{}="3'";
 (18,0) *={}="";
 (-2,6)*{}="X1";
 (-10,6)*{}="X2";
 (-2,-6)*{}="XA";
 (-10,-6)*{}="XB";
 (-10,-12)*{}="XB'";
     "XB";"XB'" **\dir{-}; ?(.4)*\dir{>};
     "X1";"1" **\dir{-}; ?(0)*\dir{>};
     "X2";"0" **\dir{-}; ?(0)*\dir{<};
   "XA";"2" **\crv{(-1,-15)&(6,-15)};
     ?(0)*\dir{<};
   "2";"3" **\crv{(6,1)&(14,1)};
     ?(0)*\dir{<}+(4,8)*{\scriptstyle \fix} ;
     ?(.95)*\dir{<}+(5,-2)*{\scriptstyle F(\xb)};
   "3";"3'" **\dir{-};
     "X1";"X2" **\crv{(-2,-1)&(-10,-1)};
     "XA";"XB" **\crv{(-2,1)&(-10,1)};
   (-13,0)*{}="x1";
   (-13,-13)*{}="x2";
   (7,0)*{}="x4";
   (7,-13)*{}="x3";
     "x1";"x2"; **\dir{.};
     "x2";"x3"; **\dir{.};
     "x3";"x4"; **\dir{.};
     "x4";"x1"; **\dir{.};

\endxy
\]
\[
\xy (-15,0)*{};(15,0)*{};
\endxy
\qquad = \qquad \xy
 (-10,12)*{}="0";
 (-2,12)*{}="1";
 (6,-7)  *{}="2";
 (6,-13)  *{}="2'";
 (14,-7) *{}="3";
 (14,-13)*{}="3'";
 (18,0) *={}="";
 (-2,6)*{}="X1";
 (-10,6)*{}="X2";
     "X1";"1" **\dir{-}; ?(0)*\dir{>};
     "X2";"0" **\dir{-}; ?(0)*\dir{<};
   "2";"3" **\crv{(6,1)&(14,1)};
     ?(0)*\dir{<}+(4,8)*{\scriptstyle \fix} ;
     ?(.95)*\dir{<}+(5,-2)*{\scriptstyle F(\xb)};
   "3";"3'" **\dir{-};
   "2";"2'" **\dir{-};
     "X1";"X2" **\crv{(-2,-1)&(-10,-1)};
\endxy
\]
\[
\xy (-15,0)*{};(15,0)*{};
\endxy
\qquad = \qquad \xy
 (-10,12)*{}="0";
 (-2,12)*{}="1";
 (6,-7)  *{}="2";
 (6,-13)  *{}="2'";
 (14,-7) *{}="3";
 (14,-13)*{}="3'";
 (18,0) *={}="";
 (-2,6)*{}="X1";
 (-10,6)*{}="X2";
     "X1";"1" **\dir{-}; ?(0)*\dir{>};
     "X2";"0" **\dir{-}; ?(0)*\dir{<};
   "2";"3" **\crv{(6,1)&(14,1)};
     ?(0)*\dir{<}+(4,8)*{\scriptstyle \fix} ;
     ?(.95)*\dir{<}+(5,-2)*{\scriptstyle F(\xb)};
   "3";"3'" **\dir{-};
   "2";"2'" **\dir{-};
     "X1";"X2" **\crv{(-2,-1)&(-10,-1)};
     (-6,5)*{}="x1";
     (10,-5)*{}="x2";
     "x1";"x2" **\dir{--};
\endxy
\]
\[
\xy (-15,0)*{};(15,0)*{};
\endxy
\qquad = \qquad \xy
 (0,12)*{}="0";
 (8,12)*{}="1";
 (0,-8)  *{}="2";
 (0,-13)  *{}="2'";
 (8,-8) *{}="3";
 (8,-13)*{}="3'";
 (18,0) *={}="";
 (-10,0) *={}="";
 (8,8)*{}="X1";
 (0,8)*{}="X2";
     "X1";"1" **\dir{-}; ?(.0)*\dir{>};
     "X2";"0" **\dir{-}; ?(0)*\dir{<};
   "2";"3" **\crv{(0,-2)&(8,-2)};
     ?(1)*\dir{<}+(5,-2)*{\scriptstyle F(\xb)};
   "3";"3'" **\dir{-};
   "2";"2'" **\dir{-};?(.25)*\dir{>}+(4,7.5)*{\scriptstyle \fix} ;
     "X1";"X2" **\crv{(8,2)&(0,2)};
\endxy
\]

To conclude, we must show that $F_{-1}$ also satisfies the
coherence law \textbf{H2}.  In string notation, this law says:
\[
\xy (-15,0)*{};(15,0)*{}; (5,13)*{}; (5,-13)*{};
**\dir{-}?(.5)*\dir{<}+(4,0)*{\scriptstyle F(x)};
(-5,0)*{\scriptstyle F_{-1}}; (-5,0)*\xycircle(2.75,2.75){-}="f";
(-5,13)**\dir{-} ?(.75)*\dir{>}+(-5,-1)*{\scriptstyle
\overline{F(x)}}; "f";(-5,-13)**\dir{-}
?(.5)*\dir{<}+(-5,-1)*{\scriptstyle F(\xb)};
\endxy
\qquad = \qquad \xy
 (-5,13)*{};
 (5,13)*{}
   **\crv{(-5,2) & (5,2)} ?(.12)*\dir{>} ?(.92)*\dir{>};
 (5,-15)*{};
 (-5,-15)*{}
   **\crv{(5,-4) & (-5,-4)} ?(.12)*\dir{>} ?(.92)*\dir{>};
(1,2.5)*{\scriptstyle e_{F(x)}};
     (1,-5)*{\scriptstyle \feix };
\endxy
\]
where we set
\[
 \fex = F_2 \circ F(e_x) \circ F_0^{-1} \maps F(\xb) \ten F(x) \to 1' .
\]

Again, the proof is a string diagram calculation. Here we need the
fact that $(F(x),F(\xb),\fix,\fex)$ is an adjunction.  This allows
us to use a zig-zag identity for $\fix$ and $\fex$ below:
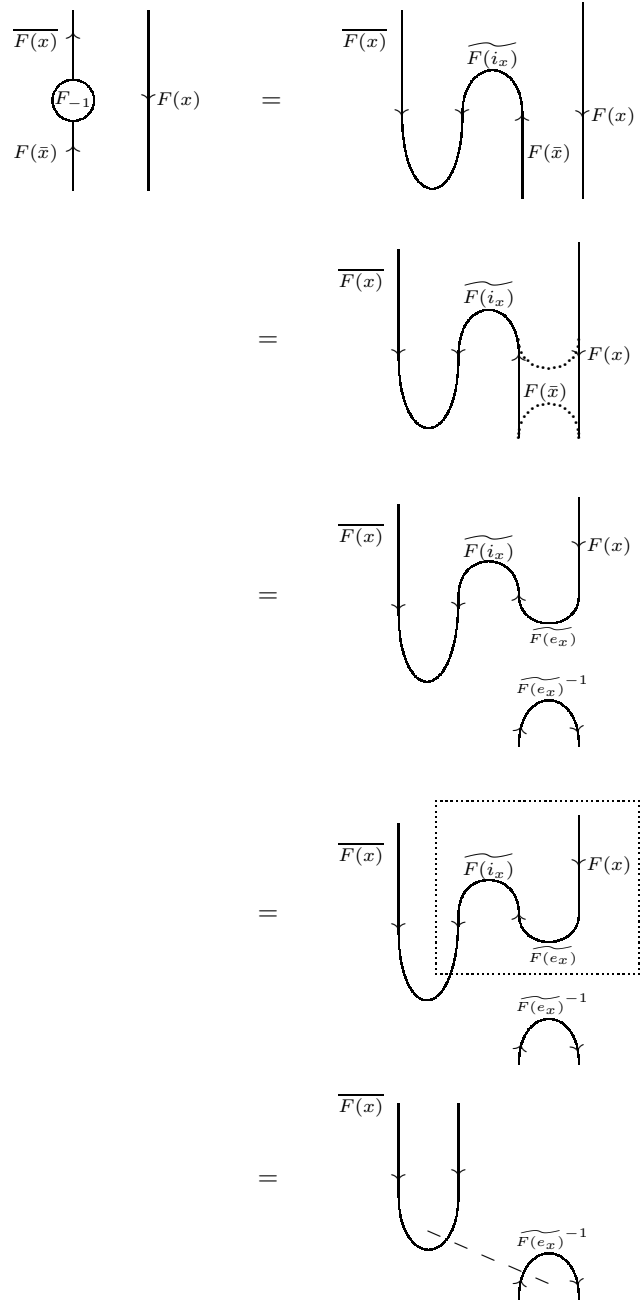
\[
\xy (-12,0)*{};(12,0)*{}; (5,12)*{}; (5,-12)*{};
    **\dir{-}?(.5)*\dir{<}+(4,0)*{\scriptstyle F(x)};
(-5,0)*{\scriptstyle F_{-1}}; (-5,0)*\xycircle(2.75,2.75){-}="f";
(-5,12)**\dir{-} ?(.75)*\dir{>}+(-5,-1)*{\scriptstyle
\overline{F(x)}}; "f";(-5,-12)**\dir{-}
?(.5)*\dir{<}+(-5,-1)*{\scriptstyle F(\xb)};
\endxy
\qquad = \qquad \xy (-15,0)*{};(22,0)*{};
 (14,13)*{}; (14,-13)*{}; **\dir{-}?(.42)*\dir{<}+(4,0)*{\scriptstyle F(x)};
 (-10,12)*{}="1'";
 (-10,-2)*{}="1";
 (-2,-2)  *{}="2";
 (6,-2) *{}="3";
 (6,-13)*{}="3'";
 (18,0) *={}="";
   "1";"2" **\crv{(-10,-15)&(-2,-15)};
   "2";"3" **\crv{(-2,6)&(6,6)};
     ?(0)*\dir{<}+(4,8)*{\scriptstyle \fix} ;
     ?(1)*\dir{<}+(3.5,-5)*{\scriptstyle F(\xb)};
   "3";"3'" **\dir{-};
   "1";"1'" **\dir{-};
      ?(0)*\dir{<}+(-5,10)*{\scriptstyle \overline{F(x)}};
\endxy
\]
\[
\xy (-12,0)*{};(12,0)*{};
\endxy
\qquad = \qquad \xy (-15,0)*{};(22,0)*{};
 (14,13)*{}; (14,-13)*{}; **\dir{-}?(.42)*\dir{<}+(4,0)*{\scriptstyle F(x)};
 (-10,12)*{}="1'";
 (-10,-2)*{}="1";
 (-2,-2)  *{}="2";
 (6,-2) *{}="3";
 (6,-13)*{}="3'";
 (18,0) *={}="";
   "1";"2" **\crv{(-10,-15)&(-2,-15)};
   "2";"3" **\crv{(-2,6)&(6,6)};
     ?(0)*\dir{<}+(4,8)*{\scriptstyle \fix} ;
     ?(1)*\dir{<}+(3.5,-5)*{\scriptstyle F(\xb)};
   "3";"3'" **\dir{-};
   "1";"1'" **\dir{-};
      ?(0)*\dir{<}+(-5,10)*{\scriptstyle \overline{F(x)}};
  (6,0)*{}="Tl";
 (14,0)*{}="Tr";
 (6,-13)*{}="Bl";
 (14,-13)*{}="Br";
     "Tl";"Tr" **\crv{~*=<2pt>{.}(6,-5)&(14,-5)};
     "Bl";"Br" **\crv{~*=<2pt>{.}(6,-7)&(14,-7)};
\endxy
\]

\[
\xy (-12,0)*{};(12,0)*{};
\endxy
\qquad = \qquad \xy (-15,0)*{};(22,0)*{};
 (14,13)*{}; (14,0)*{};
    **\dir{-}?(.5)*\dir{<}+(4,0)*{\scriptstyle F(x)};
 (-10,12)*{}="1'";
 (-10,-2)*{}="1";
 (-2,0)  *{}="2";
 (6,0) *{}="3";
 (6,-15)*{}="3'";
 (18,0) *={}="";
   "1";"2" **\crv{(-10,-15)&(-2,-15)};
   ?(.96)*\dir{<}+(4,8)*{\scriptstyle \fix} ;
   "2";"3" **\crv{(-2,6)&(6,6)};
     ?(1)*\dir{<}+(4,8)*{} ;
   "1";"1'" **\dir{-};
      ?(.0)*\dir{<}+(-5,10)*{\scriptstyle \overline{F(x)}};
 (6,0)*{}="Tl";
 (14,0)*{}="Tr";
 (6,-20)*{}="Bl";
 (14,-20)*{}="Br";
     "Tl";"Tr" **\crv{(6,-5)&(14,-5)};
     "Bl";"Br" **\crv{(6,-12)&(14,-12)};
     ?(.92)*\dir{>}; ?(.12)*\dir{>};
     (10.5,-5.5)*{\scriptscriptstyle \fex};
     (10.5,-12)*{\scriptscriptstyle \feix};
\endxy
\]

\[
\xy (-12,0)*{};(12,0)*{};
\endxy
\qquad = \qquad \xy (-15,0)*{};(22,0)*{};
 (14,13)*{}; (14,0)*{};
    **\dir{-}?(.5)*\dir{<}+(4,0)*{\scriptstyle F(x)};
 (-10,12)*{}="1'";
 (-10,-2)*{}="1";
 (-2,0)  *{}="2";
 (6,0) *{}="3";
 (6,-15)*{}="3'";
 (18,0) *={}="";
   "1";"2" **\crv{(-10,-15)&(-2,-15)};
   ?(.96)*\dir{<}+(4,8)*{\scriptstyle \fix} ;
   "2";"3" **\crv{(-2,6)&(6,6)};
     ?(1)*\dir{<}+(4,8)*{} ;
   "1";"1'" **\dir{-};
      ?(.0)*\dir{<}+(-5,10)*{\scriptstyle \overline{F(x)}};
 (6,0)*{}="Tl";
 (14,0)*{}="Tr";
 (6,-20)*{}="Bl";
 (14,-20)*{}="Br";
     "Tl";"Tr" **\crv{(6,-5)&(14,-5)};
     "Bl";"Br" **\crv{(6,-12)&(14,-12)};
     ?(.92)*\dir{>}; ?(.12)*\dir{>};
     (10.5,-5.5)*{\scriptscriptstyle \fex};
     (10.5,-12)*{\scriptscriptstyle \feix};
         (-5,15)*{}="X1";
         (22,15)*{}="X3";
         (22,-8)*{}="X4";
         (-5,-8)*{}="X2";
             "X1";"X3" **\dir{.};
             "X2";"X1" **\dir{.};
             "X3";"X4" **\dir{.};
             "X4";"X2" **\dir{.};
\endxy
\]
\[
\xy (-12,0)*{};(12,0)*{};
\endxy
\qquad = \qquad \xy (-15,0)*{};(22,0)*{};
 (-10,10)*{}="1'";
 (-2,10)*{}="2'";
 (-10,-2)*{}="1";
 (-2,-2)  *{}="2";
 (18,0) *={}="";
   "1";"2" **\crv{(-10,-12)&(-2,-12)};
   "1";"1'" **\dir{-};
      ?(.15)*\dir{<}+(-5,10)*{\scriptstyle \overline{F(x)}};
   "2";"2'" **\dir{-};
     ?(.2)*\dir{<};
 (6,-16)*{}="Bl";
 (14,-16)*{}="Br";
     "Bl";"Br" **\crv{(6,-8)&(14,-8)};
     ?(.92)*\dir{>}; ?(.12)*\dir{>};
     (10.5,-8)*{\scriptscriptstyle \feix};
 (-6,-7)*{}; (10,-14)*{}; **\dir{--};
\endxy
\]

\[
\xy (-15,0)*{};(15,0)*{};
\endxy
\qquad = \qquad \xy
 (-15,0)*{};(22,0)*{};
 (0,12)*{}="0";
 (8,12)*{}="1";
 (0,-8)  *{}="2";
 (0,-13)  *{}="2'";
 (8,-8) *{}="3";
 (8,-13)*{}="3'";
 (18,0) *={}="";
 (-10,0) *={}="";
 (8,8)*{}="X1";
 (0,8)*{}="X2";
     "X1";"1" **\dir{-}; ?(0)*\dir{>};
     "X2";"0" **\dir{-}; ?(0)*\dir{<};
   "2";"3" **\crv{(0,-2)&(8,-2)};
     ?(1)*\dir{<}+(5,-2)*{\scriptstyle F(\xb)};
   "3";"3'" **\dir{-};
   "2";"2'" **\dir{-};?(.15)*\dir{>}+(4,9.5)*{\scriptscriptstyle \feix} ;
     "X1";"X2" **\crv{(8,2)&(0,2)};
\endxy
\]
\hskip 30em \qed

In short, we do not need to include $F_{-1}$ and its coherence
laws in the definition of a coherent 2-group homomorphism; we get
them `for free'.

\section{Internalization} \label{internalizationsection}

The concept of `group' was born in the category $\Set$, but groups
can live in other categories too.  This vastly enhances
the power of group theory: for example, we have `topological groups',
`Lie groups', `affine group schemes', and so on --- each with their own
special features, but all related.

The theory of 2-groups has a similar flexibility.  Since 2-groups
are categories, we have implicitly defined the concept of 2-group
in the 2-category $\Cat$.   However, as noted by Joyal and Street,
this concept can generalized to other 2-categories as well \cite{JS}.
This makes it possible to define `topological
2-groups', `Lie 2-groups', `affine 2-group schemes' and the like.
In this section we describe how this generalization works.  In the
next section, we give many examples of Lie 2-groups.

`Internalization' is an efficient method of generalizing concepts
from the category of sets to other categories.  To internalize a
concept, we need to express it in a purely diagrammatic form.
Mac Lane illustrates this in his classic text \cite{MacLane3} by
internalizing the concept of a `group'.  We can define this notion
using commutative diagrams by specifying:
\begin{itemize}
 \item a set $G$,
\end{itemize}
together with
\begin{itemize}
 \item a multiplication function $m \maps G \times G \to G$,
 \item an identity element for the multiplication given by the function
 $\id \maps I \to G$ where $I$ is the terminal object in $\Set$,
 \item a function $\inv \maps G \to G$,
\end{itemize}
such that the following diagrams commute:
\begin{itemize}
\item the associative law:
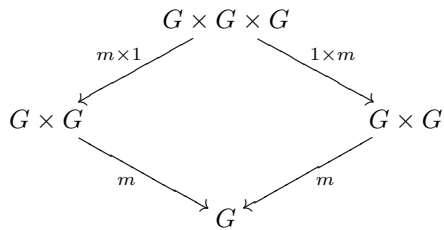
\[ \vcenter{
\xymatrix{ &   G \times G \times G \ar[dr]^{1 \times m}
   \ar[dl]_{m \times 1} \\
 G \times G \ar[dr]_{m}
&&  G \times G \ar[dl]^{m}  \\
&  G }}
\]
\item the right and left unit laws:
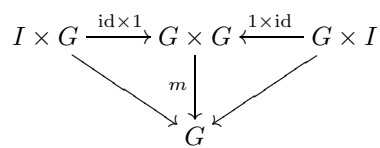
\[ \vcenter{
\xymatrix{
 I \times G \ar[r]^{\id \times 1} \ar[dr]
& G \times G \ar[d]_{m}
& G \times I \ar[l]_{1 \times \id} \ar[dl] \\
& G }}
\]
\item the right and left inverse laws:
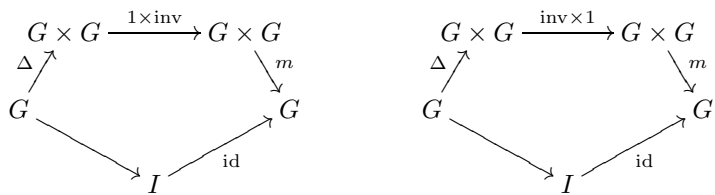
\[
\xy (-12,10)*+{G \times G}="TL"; (12,10)*+{G \times G}="TR";
(-18,0)*+{G}="ML"; (18,0)*+{G}="MR"; (0,-10)*+{I}="B";
     {\ar_{} "ML";"B"};
     {\ar^{\Delta} "ML";"TL"};
     {\ar_{\id} "B";"MR"};
     {\ar^{m} "TR";"MR"};
     {\ar^{1 \times \inv } "TL";"TR"};
\endxy
\qquad \qquad \xy (-12,10)*+{G \times G}="TL"; (12,10)*+{G \times
G}="TR"; (-18,0)*+{G}="ML"; (18,0)*+{G}="MR"; (0,-10)*+{I}="B";
     {\ar_{} "ML";"B"};
     {\ar^{\Delta} "ML";"TL"};
     {\ar_{\id} "B";"MR"};
     {\ar^{m} "TR";"MR"};
     {\ar^{\inv \times 1} "TL";"TR"};
\endxy
\]
\end{itemize}
where $\Delta \maps G \to G \times G$ is the diagonal map.

To internalize the concept of group we simply replace the set $G$
by an {\it object} in some category $K$ and replace the functions
$m, \id,$ and $\inv$ by {\it morphisms} in this category.  Since
the definition makes use of the Cartesian product $\times$, the
terminal object $I$, and the diagonal $\Delta$, the category $K$
should have finite products. Making these substitutions in the
above definition, we arrive at the definition of a {\bf group
object in $K$}.  We shall usually call this simply a {\bf group in
$K$}.

In the special case where $K = \Set$, a group in $K$ reduces to an
ordinary group.  A topological group is a group in $\Top$, a Lie
group is a group in $\Diff$, and a affine group scheme is a group
in $\CommRing^\op$, usually called the category of `affine
schemes'.  Indeed, for any category $K$ with finite products,
there is a category $K\Grp$ consisting of groups in $K$ and
homomorphisms between these, where a {\bf homomorphism} $f \maps G
\to G'$ is a morphism in $K$ that preserves multiplication,
meaning that this diagram commutes:
\[
\xymatrix{ G \times G
 \ar[rr]^{m}
 \ar[dd]_{f \times f}
  && G
 \ar[dd]^{f} \\ \\
  G' \times G'
 \ar[rr]^{m'}
  && G'}
\]
As usual, this implies that $f$ also preserves the identity and
inverses.

Following Joyal and Street \cite{JS}, let us now internalize the
concept of coherent 2-group and define a 2-category of `coherent 2-groups
in $K$' in a similar manner.  For this, one must first define a coherent
2-group using only commutative diagrams.   However, since the usual group
axioms hold only up to natural isomorphism in a coherent 2-group, these
will be 2-categorical rather than 1-categorical diagrams. As a
result, the concept of coherent 2-group will make sense in any
2-category with finite products, $K$.  For simplicity we shall
limit ourselves to the case where $K$ is a strict 2-category.

To define the concept of coherent 2-group using commutative
diagrams, we start with a category $C$ and equip it with a
multiplication functor $m \maps C \times C \to C$ together with an
identity object for multiplication given by the functor $\id \maps
I \to C$, where $I$ is the terminal category.  The functor mapping
each object to its specified weak inverse is a bit more subtle!
One can try to define a functor $* \maps C \to C$ sending each
object $x \in C$ to its specified weak inverse $\bar{x}$, and
acting on morphisms as follows:
\[ * \maps \quad
\xy (0,0)*{f}; (0,0)*\xycircle(2.65,2.65){-}="f"; (0,11)**\dir{-}
?(.5)*\dir{<}+(3,0)*{\scriptstyle x}; "f";(0,-11)**\dir{-}
?(.75)*\dir{>}+(3,0)*{\scriptstyle y};
\endxy
\quad  \mapsto \quad \xy
 (0,0)*{f};
 (0,0)*\xycircle(2.65,2.65){-}="M";
 (-8,-6)*{}="L";
 (8,6)*{}="R";
 (-8,10)*{}="LB";
 (8,-10)*{}="RT";
 (0,6)*{}="MB";
 (0,-6)*{}="MT";
     "L";"MT" **\crv{(-8,-12)&(0,-12)}; ?(.96)*\dir{<};
     "MB";"R" **\crv{(0,12)&(8,12)};  ?(0)*\dir{<};
         "L";"LB" **\dir{-}; ?(.45)*\dir{>};
         "R";"RT" **\dir{-}; ?(.25)*\dir{<};
         "M";"MT" **\dir{-};
         "M";"MB" **\dir{-};
             (-4,-12)*{\scriptstyle e_y};
             (4,12)*{\scriptstyle i_x};
\endxy
\]
However, $*$ is actually a {\it contravariant} functor.  To see
this, we consider composable morphisms $f \maps x \to y$ and $g
\maps y \to z$ and check that $(fg)^* = g^* f^*$.  In string
diagram form, this equation says:
\[\xy
 (0,0)*{fg};
 (0,0)*\xycircle(2.65,2.65){-}="M";
 (-8,-6)*{}="L";
 (8,6)*{}="R";
 (-8,10)*{}="LB";
 (8,-10)*{}="RT";
 (0,6)*{}="MB";
 (0,-6)*{}="MT";
     "L";"MT" **\crv{(-8,-12)&(0,-12)}; ?(.96)*\dir{<};
     "MB";"R" **\crv{(0,12)&(8,12)};  ?(0)*\dir{<};
         "L";"LB" **\dir{-}; ?(.45)*\dir{>};
         "R";"RT" **\dir{-}; ?(.25)*\dir{<};
         "M";"MT" **\dir{-};
         "M";"MB" **\dir{-};
             (-4,-12)*{\scriptstyle e_z};
             (4,12)*{\scriptstyle i_x};
\endxy
\quad = \quad  \xy
 (0,0)*{g};
 (0,0)*\xycircle(2.65,2.65){-}="M";
 (-8,-6)*{}="L";
 (8,6)*{}="R";
 (-8,10)*{}="LB";
 (8,-6)*{}="RT";
 (0,6)*{}="MB";
 (0,-6)*{}="MT";
     "L";"MT" **\crv{(-8,-12)&(0,-12)}; ?(.96)*\dir{<};
     "MB";"R" **\crv{(0,12)&(8,12)};  ?(0)*\dir{<};
         "L";"LB" **\dir{-}; ?(.45)*\dir{>};
         "R";"RT" **\dir{-}?(.5)*\dir{<};
         "M";"MT" **\dir{-};
         "M";"MB" **\dir{-};
             (-4,-12)*{\scriptstyle e_z};
             (4,12)*{\scriptstyle i_y};
 (16,0)*{f};
 (16,0)*\xycircle(2.65,2.65){-}="M";
 (8,-6)*{}="L";
 (24,6)*{}="R";
 (8,6)*{}="LB";
 (24,-10)*{}="RT";
 (16,6)*{}="MB";
 (16,-6)*{}="MT";
     "L";"MT" **\crv{(8,-12)&(16,-12)}; ?(.96)*\dir{<};
     "MB";"R" **\crv{(16,12)&(24,12)};  ?(0)*\dir{<};
         "L";"LB" **\dir{-};
         "R";"RT" **\dir{-}; ?(.25)*\dir{<};
         "M";"MT" **\dir{-};
         "M";"MB" **\dir{-};
             (12,-12)*{\scriptstyle e_y};
             (20,12)*{\scriptstyle i_x};
\endxy
\]
This equation holds if and only if
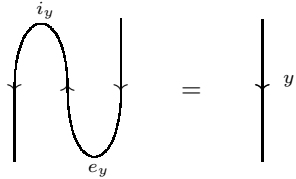
\[
\xy0;/r.28pc/: (-6,-8)*{}="1E"; (-6,0)*{}="1"; (0,0)*{}="2";
(6,0)*{}="3"; (6,8)*{}="3B";
 "2";"1" **\crv{(0,10)& (-6,10)}
     ?(.02)*\dir{>}  ?(1)*\dir{>};
 "3";"2" **\crv{(6,-10)& (0,-10)} ;
 "1";"1E" **\dir{-};
 "3B";"3" **\dir{-}?(1)*\dir{>} ;
(-2.5,9)*{\scriptstyle i_y}; (3.5,-9)*{\scriptstyle e_y};
\endxy
\qquad = \xy0;/r.28pc/: (-6,8)*{}; (0,8)*{}; (0,-8)*{}; **\dir{-}
?(.5)*\dir{<}; (6,-8)*{}?(.5)*\dir{<}+(3,1)*{\scriptstyle y};
\endxy
\]
But this is merely the first zig-zag identity!

Contravariant functors are a bit annoying since they are not
really morphisms in $\Cat$. Luckily, there is also another
contravariant functor ${}^{-1} \maps C \to C$ sending each
morphism to its inverse, expressed diagrammatically as
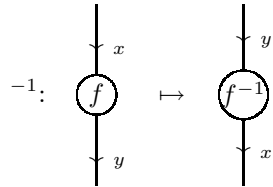
\[ ^{-1} \maps \quad
\xy (0,0)*{f}; (0,0)*\xycircle(2.65,2.65){-}="f"; (0,12)**\dir{-}
?(.5)*\dir{<}+(3,0)*{\scriptstyle x}; "f";(0,-12)**\dir{-}
?(.75)*\dir{>}+(3,0)*{\scriptstyle y};
\endxy
\quad \mapsto \quad \xy (0,0)*{f^{-1}};
(0,0)*\xycircle(3.25,3.25){-}="f"; (0,12)**\dir{-}
?(.59)*\dir{<}+(3,0)*{\scriptstyle y}; "f";(0,-12)**\dir{-}
?(.65)*\dir{>}+(3,0)*{\scriptstyle x};
\endxy
\]
If we compose the contravariant functor $*$ with this, we obtain a
covariant functor $\inv \maps C \to C$ given by
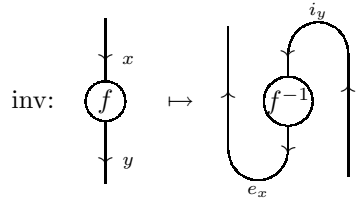
\[ \inv \maps \quad
\xy (0,0)*{f}; (0,0)*\xycircle(2.65,2.65){-}="f"; (0,11)**\dir{-}
?(.5)*\dir{<}+(3,0)*{\scriptstyle x}; "f";(0,-11)**\dir{-}
?(.75)*\dir{>}+(3,0)*{\scriptstyle y};
\endxy
\quad  \mapsto \quad \xy
 (0,0)*{f^{-1}};
 (0,0)*\xycircle(3.25,3.25){-}="M";
 (-8,-6)*{}="L";
 (8,6)*{}="R";
 (-8,10)*{}="LB";
 (8,-10)*{}="RT";
 (0,6)*{}="MB";
 (0,-6)*{}="MT";
     "L";"MT" **\crv{(-8,-12)&(0,-12)}; ?(.96)*\dir{<};
     "MB";"R" **\crv{(0,12)&(8,12)};  ?(0)*\dir{<};
         "L";"LB" **\dir{-}; ?(.45)*\dir{>};
         "R";"RT" **\dir{-}; ?(.25)*\dir{<};
         "M";"MT" **\dir{-};
         "M";"MB" **\dir{-};
             (-4,-12)*{\scriptstyle e_x};
             (4,12)*{\scriptstyle i_y};
\endxy
\]

Thus, we can try to write the definition of a coherent 2-group in
terms of:
\begin{itemize}
 \item the category $C$,
\end{itemize}
together with
\begin{itemize}
 \item the functor $m \maps C\times C \to C$, where we
 write $m(x,y)=x \ten y$ and $m(f,g)=f \ten g$ for objects $x, y,
 \in C$ and morphisms $f, g$ in $C$,
 \item the functor $\id \maps I \to C$ where $I$ is the
 terminal category, and the object in the range of
 this functor is $1 \in C$,
 \item the functor $ \inv \maps C \to C$,
\end{itemize}
together with the following natural isomorphisms:
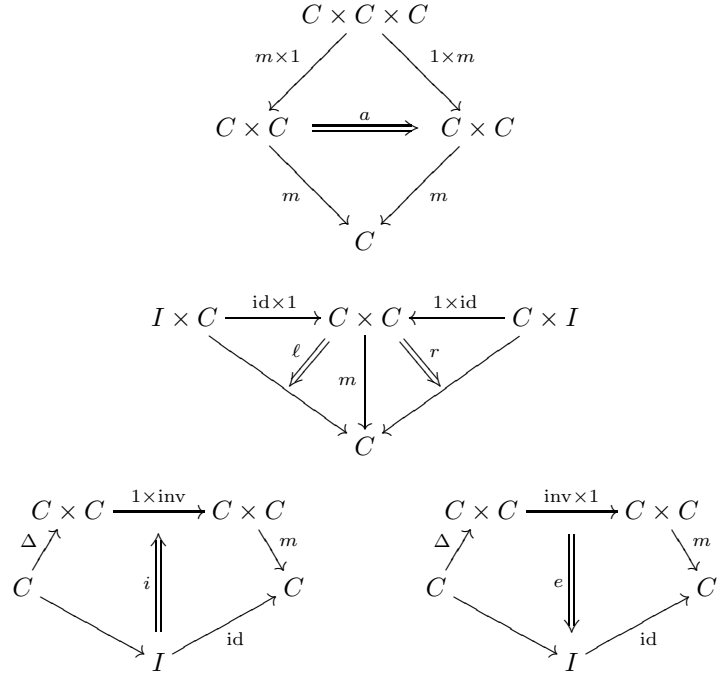
\[
\xy
 (0,15)*+{C \times C \times C}="T";
 (-15,0)*+{C \times C}="L";
 (15,0)*+{C \times C}="R";
 (0,-15)*+{C}="B";
 (-7,0)*{}="ML";
 (7,0)*{}="MR";
     {\ar^{1 \times m} "T";"R"};
     {\ar_{m \times 1} "T";"L"};
     {\ar^{m} "R";"B"};
     {\ar_{m} "L";"B"};
         {\ar@{=>}^{a} "ML";"MR"};
\endxy
\]
\[
\xy
  (-24,5)*+{I \times C}="L";
  (0,5)*+{C \times C}="M";
  (24,5)*+{C \times I}="R";
  (0,-12)*+{C}="B";
  (-5,2)*{}="TL";
  (5,2)*{}="TR";
  (10,-4)*{}="BR";
  (-10,-4)*{}="BL";
     {\ar^{\id \times 1} "L";"M"};
     {\ar_{1 \times \id} "R";"M"};
     {\ar_>>>>>>>{m} "M";"B"};
     {\ar_{} "L";"B"};
     {\ar_{} "R";"B"};
         {\ar@{=>}_{\ell} "TL";"BL"};
         {\ar@{=>}^{r} "TR";"BR"};
\endxy
\]
\[
\xy (-12,10)*+{C \times C}="TL"; (12,10)*+{C \times C}="TR";
(-18,0)*+{C}="ML"; (18,0)*+{C}="MR"; (0,-10)*+{I}="B";
     {\ar_{} "ML";"B"};
     {\ar^{\Delta} "ML";"TL"};
     {\ar_{\id} "B";"MR"};
     {\ar^{m} "TR";"MR"};
     {\ar^{1 \times \inv } "TL";"TR"};
 (0,-6)*{}="X";
 (0,7)*{}="x";
 {\ar@{=>}^{i} "X";"x"};
\endxy
\qquad \qquad \xy (-12,10)*+{C \times C}="TL"; (12,10)*+{C \times
C}="TR"; (-18,0)*+{C}="ML"; (18,0)*+{C}="MR"; (0,-10)*+{I}="B";
     {\ar_{} "ML";"B"};
     {\ar^{\Delta} "ML";"TL"};
     {\ar_{\id} "B";"MR"};
     {\ar^{m} "TR";"MR"};
     {\ar^{\inv \times 1} "TL";"TR"};
 (0,-6)*{}="x";
 (0,7)*{}="X";
 {\ar@{=>}_{e} "X";"x"};
\endxy
\]
and finally the coherence laws satisfied by these isomorphisms.
But to do this, we must write the coherence laws in a way that
does not explicitly mention objects of $C$.  For example, we must
write the pentagon identity
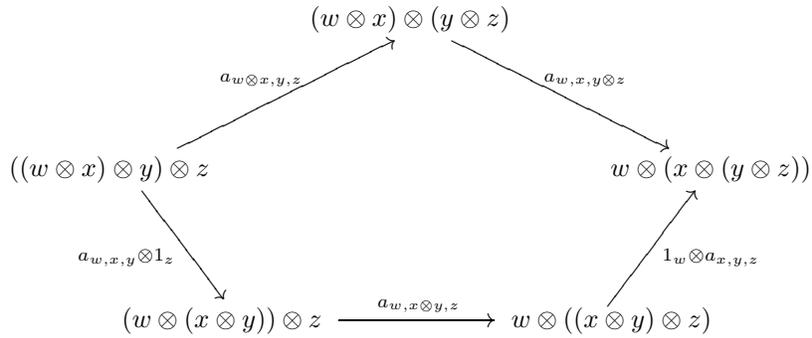
\[
\xy
 (0,20)*+{(w \ten x) \ten (y \ten z)}="1";
 (40,0)*+{w \ten (x \ten (y \ten z))}="2";
 (25,-20)*{ \quad w \ten ((x \ten y) \ten z)}="3";
 (-25,-20)*+{(w \ten (x \ten y)) \ten z}="4";
 (-40,0)*+{((w \ten x) \ten y) \ten z}="5";
     {\ar^{a_{w,x,y \ten z}}     "1";"2"}
     {\ar_{1_w \ten a _{x,y,z}}  "3";"2"}
     {\ar^{a _{w,x \ten y,z}}    "4";"3"}
     {\ar_{a _{w,x,y} \ten 1_z}  "5";"4"}
     {\ar^{a _{w \ten x,y,z}}    "5";"1"}
\endxy
\\
\]
without mentioning the objects $w,x,y,z \in C$. We can do this by
working with (for example) the functor $(1 \times 1 \times m)
\circ (1 \times m) \circ m$ instead of its value on the object
$(x,y,z,w) \in C^4$, namely $x \tensor (y \tensor (z \tensor w))$.
If we do this, we see that the diagram becomes 3-dimensional!  It
is a bit difficult to draw, but it looks something like this:
\\
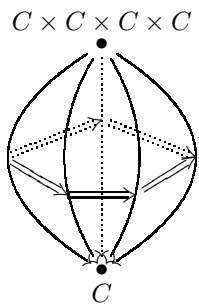
\[
\xy
 (0,15)*+{\bullet}="1";
 (0,-15)*+{\bullet}="2";
   {\ar@{.>} "1";"2"};
   "1";"2" **\crv{(10,0)} ?(.97)*\dir{>};
   "1";"2" **\crv{(-10,0)} ?(.97)*\dir{>};
   "1";"2" **\crv{(25,0)} ?(.97)*\dir{>};
   "1";"2" **\crv{(-25,0)} ?(.97)*\dir{>};
     (4.5,-5)*+{}="D";
     (-4.5,-5)*+{}="C";
     (-13,0)*+{}="B";
     (12.5,0)*+{}="E";
     (0,5)*+{}="A";
       {\ar@{=>} "C";"D"};
       {\ar@{=>} "B";"C"};
       {\ar@{=>} "D";"E"};
       {\ar@{:>} "B";"A"};
       {\ar@{:>} "A";"E"};
   (0,18)*{C \times C \times C \times C};
   (0,-18)*{C};
\endxy
\]
where the downwards-pointing single arrows are functors from $C^4$
to $C$, while the horizontal double arrows are natural
transformations between these functors, forming a commutative
pentagon.  Luckily we can also draw this pentagon in a
2-dimensional way, as follows:
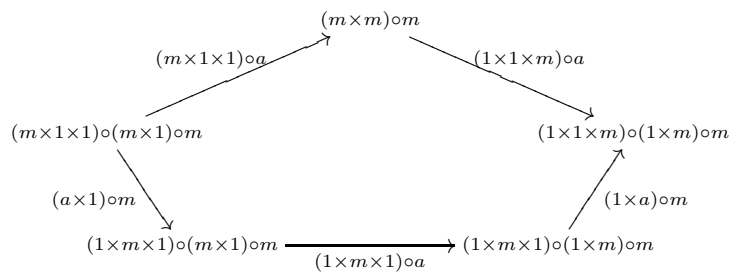
\[
\def\objectstyle{\scriptstyle}
\def\labelstyle{\scriptstyle}
\xy
 (-35,0)*+{(m \times 1 \times 1) \circ (m \times 1) \circ m}="TL";
 (-25,-15)*+{ (1 \times m \times 1)  \circ (m \times 1) \circ m}="BL";
 (35,0)*+{(1 \times 1 \times m) \circ (1 \times m) \circ m}="TR";
 (25,-15)*+{  (1 \times m \times 1) \circ ( 1 \times m ) \circ m}="BR";
 (0,15)*+{(m \times m) \circ m}="T";
     {\ar_{(a \times 1) \circ m} "TL";"BL"};
     {\ar^{(m \times 1 \times 1) \circ a} "TL";"T"};
     {\ar_{(1 \times m \times 1) \circ a} "BL";"BR"};
     {\ar_{(1 \times a) \circ m} "BR";"TR"};
     {\ar^{(1 \times 1 \times m) \circ a} "T";"TR"};
\endxy
\]

\noindent Using this idea we can write the definition of `coherent
2-group' using only the structure of $\Cat$ as a 2-category with
finite products.  We can then internalize this definition, as
follows:

\begin{defn} \label{coherentobject} \et
Given a 2-category $K$ with finite products, a
{\bf coherent 2-group in} $K$ consists of:
\begin{itemize}
 \item an object $C \in K$,
\end{itemize}
together with:
\begin{itemize}
 \item a {\bf multiplication} morphism $m \maps C\times C \to C$,
 \item an {\bf identity-assigning}
  morphism $\id \maps I \to C$ where $I$ is the
 terminal object of $K$,
 \item an {\bf inverse} morphism $ \inv \maps C \to C$,
\end{itemize}
together with the following 2-isomorphisms:
\begin{itemize}
\item the {\bf associator}:
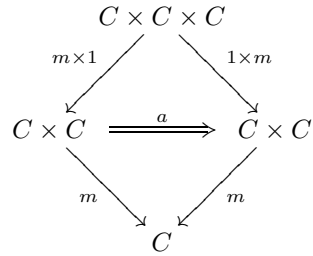
\[
\xy
 (0,15)*+{C \times C \times C}="T";
 (-15,0)*+{C \times C}="L";
 (15,0)*+{C \times C}="R";
 (0,-15)*+{C}="B";
 (-7,0)*{}="ML";
 (7,0)*{}="MR";
     {\ar^{1 \times m} "T";"R"};
     {\ar_{m \times 1} "T";"L"};
     {\ar^{m} "R";"B"};
     {\ar_{m} "L";"B"};
         {\ar@{=>}^{a} "ML";"MR"};
\endxy
\]
\item the {\bf left and right unit laws}:
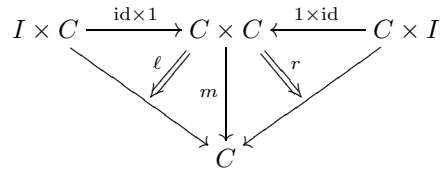
\[
\xy
  (-24,5)*+{I \times C}="L";
  (0,5)*+{C \times C}="M";
  (24,5)*+{C \times I}="R";
  (0,-12)*+{C}="B";
  (-5,2)*{}="TL";
  (5,2)*{}="TR";
  (10,-4)*{}="BR";
  (-10,-4)*{}="BL";
     {\ar^{\id \times 1} "L";"M"};
     {\ar_{1 \times \id} "R";"M"};
     {\ar_>>>>>>>{m} "M";"B"};
     {\ar_{} "L";"B"};
     {\ar_{} "R";"B"};
         {\ar@{=>}_{\ell} "TL";"BL"};
         {\ar@{=>}^{r} "TR";"BR"};
\endxy
\]
\item the {\bf unit} and {\bf counit}:
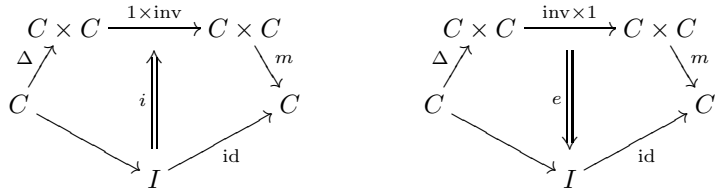
\[
\xy (-12,10)*+{C \times C}="TL"; (12,10)*+{C \times C}="TR";
(-18,0)*+{C}="ML"; (18,0)*+{C}="MR"; (0,-10)*+{I}="B";
     {\ar_{} "ML";"B"};
     {\ar^{\Delta} "ML";"TL"};
     {\ar_{\id} "B";"MR"};
     {\ar^{m} "TR";"MR"};
     {\ar^{1 \times \inv } "TL";"TR"};
 (0,-6)*{}="X";
 (0,7)*{}="x";
 {\ar@{=>}^{i} "X";"x"};
\endxy
\qquad \qquad \xy (-12,10)*+{C \times C}="TL"; (12,10)*+{C \times
C}="TR"; (-18,0)*+{C}="ML"; (18,0)*+{C}="MR"; (0,-10)*+{I}="B";
     {\ar_{} "ML";"B"};
     {\ar^{\Delta} "ML";"TL"};
     {\ar_{\id} "B";"MR"};
     {\ar^{m} "TR";"MR"};
     {\ar^{\inv \times 1} "TL";"TR"};
 (0,-6)*{}="x";
 (0,7)*{}="X";
 {\ar@{=>}_{e} "X";"x"};
\endxy
\]
\end{itemize}
such that the following diagrams commute:
\begin{itemize}
 \item the {\bf pentagon identity} for the associator:
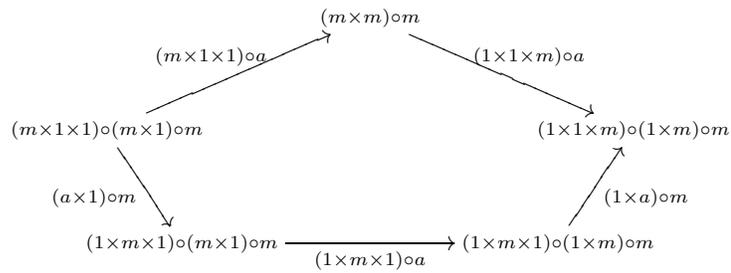
\[
\def\objectstyle{\scriptstyle}
\def\labelstyle{\scriptstyle}
\xy
 (-35,0)*+{(m \times 1 \times 1) \circ (m \times 1) \circ m}="TL";
 (-25,-15)*+{ (1 \times m \times 1)  \circ (m \times 1) \circ m}="BL";
 (35,0)*+{(1 \times 1 \times m) \circ (1 \times m) \circ m}="TR";
 (25,-15)*+{  (1 \times m \times 1) \circ ( 1 \times m ) \circ m}="BR";
 (0,15)*+{(m \times m) \circ m}="T";
     {\ar_{(a \times 1) \circ m} "TL";"BL"};
     {\ar^{(m \times 1 \times 1) \circ a} "TL";"T"};
     {\ar_{(1 \times m \times 1) \circ a} "BL";"BR"};
     {\ar_{(1 \times a) \circ m} "BR";"TR"};
     {\ar^{(1 \times 1 \times m) \circ a} "T";"TR"};
\endxy
\]
\item the {\bf triangle identity} for the left and right unit
laws:
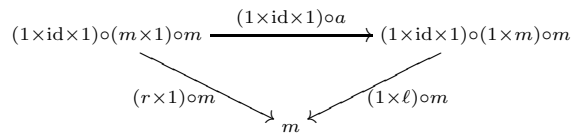
\[
\def\objectstyle{\scriptstyle}
\def\labelstyle{\scriptstyle}
\xymatrix{
  (1 \times \id \times 1) \circ (m \times 1) \circ m
      \ar[rr]^{(1 \times \id \times 1) \circ a}
      \ar[dr]_{(r \times 1) \circ m}
  &&  (1 \times \id \times 1) \circ(1 \times m) \circ m
      \ar[dl]^{(1 \times \ell) \circ m}     \\
  & m   }
\]
 \item the {\bf first zig-zag identity}:
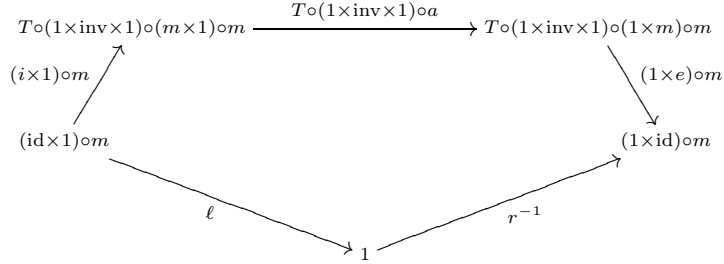
\[
\def\objectstyle{\scriptstyle}
\def\labelstyle{\scriptstyle}
\xy
 (55,0)*{};
 (-55,0)*{};
 (-40,0)*+{(\id \times 1) \circ m}="ML";
 (0,-15)*+{1}="B";
 (40,0)*+{(1 \times \id) \circ m }="MR";
 (31,15)*+{T \circ (1 \times \inv \times 1)
             \circ (1 \times m) \circ m }="TR";
 (-31,15)*+{T \circ (1 \times \inv \times 1)
              \circ (m \times 1) \circ m}="TL";
     {\ar_{\ell } "ML";"B"};
     {\ar^{(i \times 1) \circ m} "ML";"TL"};
     {\ar_{r^{-1}} "B";"MR"};
     {\ar^{(1 \times e)\circ m } "TR";"MR"};
     {\ar^{T \circ (1 \times \inv \times 1) \circ a} "TL";"TR"};
\endxy
\]
 \item the {\bf second zig-zag identity}:
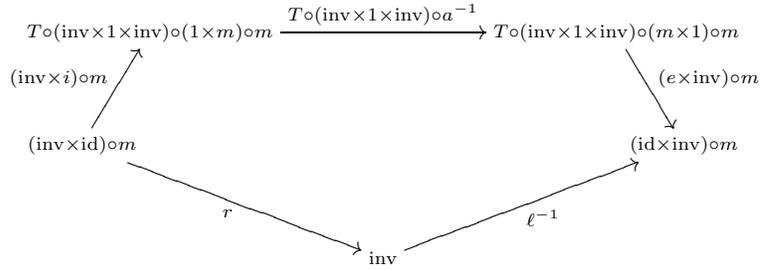
\[
\def\objectstyle{\scriptstyle}
\def\labelstyle{\scriptstyle}
\xy
 (55,0)*{};
 (-55,0)*{};
 (-40,0)*+{(\inv \times \id) \circ m}="ML";
 (0,-15)*+{\inv}="B";
 (40,0)*+{(\id \times \inv) \circ m }="MR";
 (31,15)*+{T \circ (\inv \times 1 \times \inv)
             \circ (m \times 1) \circ m}="TR";
 (-31,15)*+{T \circ (\inv \times 1 \times \inv)
              \circ (1 \times m) \circ m}="TL";
     {\ar_{r } "ML";"B"};
     {\ar^{(\inv \times i) \circ m} "ML";"TL"};
     {\ar_{\ell^{-1}} "B";"MR"};
     {\ar^{(e \times \inv) \circ m } "TR";"MR"};
     {\ar^{T \circ (\inv \times 1 \times \inv) \circ a^{-1}} "TL";"TR"};
\endxy
\]
\end{itemize}
where $T \maps C \to C^3$ is built using the diagonal functor.
\end{defn}

\begin{prop} \et
A coherent 2-group in $\Cat$ is the same as a coherent 2-group.
\end{prop}

\textbf{Proof. } Clearly any coherent 2-group gives a coherent
2-group in $\Cat$.  Conversely, suppose $C$ is a coherent 2-group
in $\Cat$.  It is easy to check that $C$ is a weak monoidal
category and that for each object $x \in C$ there is an adjoint
equivalence $(x,\xb,i_x,e_x)$ where $\xb = \inv(x)$.  This permits
the use of string diagrams to verify the one remaining point,
which is that all morphisms in $C$ are invertible.

To do this, for any morphism $f \maps x \to y$ we define a
morphism $f^{-1} \maps y \to x$ by
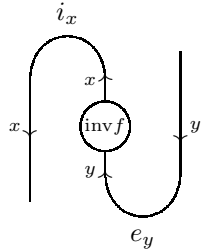
\[
\xy
 (0,0)*{\scriptstyle \inv f};
 (0,0)*\xycircle(3.3,3.3){-}="M";
 (-10,6)*{}="L";
 (10,-6)*{}="R";
 (-10,-10)*{}="LB";
 (10,10)*{}="RT";
 (0,-6)*{}="MB";
 (0,6)*{}="MT";
     "L";"MT" **\crv{(-10,14)&(0,14)}; ?(.96)*\dir{<};
     "MB";"R" **\crv{(0,-14)&(10,-14)};  ?(0)*\dir{<};
         "L";"LB" **\dir{-}; ?(.45)*\dir{>};
         "R";"RT" **\dir{-}; ?(.25)*\dir{<};
         "M";"MT" **\dir{-};
         "M";"MB" **\dir{-};
             (-12,0)*{\scriptstyle x};
             (-2,6)*{\scriptstyle x};
             (-2,-6)*{\scriptstyle y};
             (12,0)*{\scriptstyle y};
             (-5,15)*{ i_x};
             (5,-15)*{ e_y};
\endxy
\]

\noindent To check that $f^{-1}f$ is the identity, we use the fact
that $i$ is a natural isomorphism to note that this square
commutes:
\[
\xy (-13,8)*+{x \ten \xb}="TL"; (13,8)*+{y \ten \bar{y}}="TR";
(-13,-8)*+{1}="BL"; (13,-8)*+{1}="BR";
   {\ar^{f \ten \inv(f)} "TL";"TR"};
   {\ar_{i_x} "TL";"BL"};
   {\ar^{i^{-1}_y} "TR";"BR"};
   {\ar^{1_1} "BL";"BR"};
\endxy
\]

\noindent In string notation this says that:
\[
\xy (-12,0)*{f}; (-12,0)*\xycircle(2.65,2.65){-}="L";
(0,0)*{\scriptstyle \inv f}; (0,0)*\xycircle(3.3,3.3){-}="M";
(-12,6)*{}="LT"; (-12,-6)*{}="LB"; (0,6)*{}="MT"; (0,-6)*{}="MB";
   "L";"LT" **\dir{-};
     ?(.75)*\dir{<};
   "L";"LB" **\dir{-};
     ?(1)*\dir{>};
   "M";"MT" **\dir{-}+(-6,10)*{ i_x};
   "M";"MB" **\dir{-}+(-6,-10)*{ i_y^{-1}};
     ?(.85)*\dir{<};
   "LT";"MT" **\crv{(-12,16)&(0,16)};
             ?(.97)*\dir{<};
   "LB";"MB" **\crv{(-12,-16)&(0,-16)};
         (-15,6)*{\scriptstyle x};
         (-15,-6)*{\scriptstyle y};
         (3,6)*{\scriptstyle x};
         (3,-6)*{\scriptstyle y};
\endxy
\qquad = \qquad \qquad \qquad
\]

\noindent and we can use this equation to verify that $f^{-1}f =
1_y$:

\[ 
\xy (0,-5)*{f}; (-5,0)*{}; (5,0)*{}; (0,5)*{f^{-1}};
(0,5)*\xycircle(3.3,3.3){-}="f"; (0,15)**\dir{-}
?(.5)*\dir{<}+(3,0)*{\scriptstyle y}; "f";
(0,-5)*\xycircle(2.65,2.65){-}="g"; **\dir{-}
?(.4)*\dir{<}+(3,0)*{\scriptstyle x}; "f";"g";(0,-15)**\dir{-}
?(.75)*\dir{>}+(3,0)*{\scriptstyle y};
\endxy
\qquad = \qquad \xy
  (17,0)*{};
(-17,0)*{};
 (0,0)*{\scriptstyle \inv f};
 (0,0)*\xycircle(3.3,3.3){-}="M";
 (-12,6)*{}="L";
 (12,-6)*{}="R";
 (-12,-8)*{f};
 (-12,-8)*\xycircle(2.65,2.65){-}="LB";
 (-12,-16)*{}="LBB";
 (12,12)*{}="RT";
 (0,-6)*{}="MB";
 (0,6)*{}="MT";
     "L";"MT" **\crv{(-12,14)&(0,14)}; ?(.96)*\dir{<};
     "MB";"R" **\crv{(0,-14)&(12,-14)};  ?(0)*\dir{<};
         "L";"LB" **\dir{-}; ?(.45)*\dir{>};
         "R";"RT" **\dir{-}; ?(.25)*\dir{<};
         "M";"MT" **\dir{-};
         "M";"MB" **\dir{-};
             (-14,0)*{\scriptstyle x};
             (-2,6)*{\scriptstyle x};
             (-2,-6)*{\scriptstyle y};
             (14,0)*{\scriptstyle y};
             (-14,-13)*{\scriptstyle y};
             (-6,15)*{ i_x};
             (6,-15)*{ e_y};
             "LB";"LBB" **\dir{-}; ?(.75)*\dir{>};
\endxy
\]
\[
\xy (-6,0)*{};
 (6,0)*{};
\endxy
\qquad = \qquad \xy (17,0)*{}; (-17,0)*{}; (-12,0)*{f};
(-12,0)*\xycircle(2.65,2.65){-}="L"; (0,0)*{\scriptstyle \inv f};
(0,0)*\xycircle(3.3,3.3){-}="M"; (12,0)*{1_y};
(12,0)*\xycircle(2.65,2.65){-}="R"; (-12,6)*{}="LT";
(-12,-15)*{}="LB"; (0,6)*{}="MT"; (0,-6)*{}="MB"; (12,15)*{}="RT";
(12,-6)*{}="RB";
   "L";"LT" **\dir{-};
     ?(.75)*\dir{<};
   "L";"LB" **\dir{-};
   "M";"MT" **\dir{-};
   "M";"MB" **\dir{-};
   "R";"RT" **\dir{-};
     ?(.5)*\dir{<};
   "R";"RB" **\dir{-};
   "LT";"MT" **\crv{(-12,22)&(0,22)}; ?(.97)*\dir{<};
   "MB";"RB" **\crv{(0,-22)&(12,-22)}; ?(.97)*\dir{<};
         (-15,6)*{\scriptstyle x};
         (-2,6)*{\scriptstyle x};
         (15,6)*{\scriptstyle y};
         (15,-6)*{\scriptstyle y};
         (-15,-8)*{\scriptstyle y};
         (2.5,-8)*{\scriptstyle y};
   (1,-13)*{}="XTM";
   (-12,-13)*{}="XTL";
   (0,-4)*{}="XBM";
   (-12,-4)*{}="XBL";
   "XTM";"XTL" **\crv{~*=<3pt>{.} (0,-8)& (-12,-8)};
   "XBM";"XBL" **\crv{~*=<3pt>{.} (0,-9)&(-12,-9)};
\endxy
\]
\[
\xy (-6,0)*{};
 (6,0)*{};
\endxy
\qquad = \qquad \xy (17,0)*{}; (-17,0)*{}; (-12,2)*{f};
(-12,2)*\xycircle(2.65,2.65){-}="L"; (0,2)*{\scriptstyle \inv f};
(0,2)*\xycircle(3.3,3.3){-}="M"; (12,0)="R"; (-12,6)*{}="LT";
(-12,-20)*{}="LB"; (0,6)*{}="MT"; (0,-6)*{}="MB"; (12,17)*{}="RT";
(12,-15)*{}="RB"; (0,-4)*{}="XTM"; (-12,-4)*{}="XTL";
(1,-15)*{}="XBM"; (-12,-15)*{}="XBL";
   "L";"LT" **\dir{-};
     ?(1)*\dir{<};
   "L";"XTL" **\dir{-}; ?(1)*\dir{>};
   "XBL";"LB" **\dir{-}; ?(.2)*\dir{>};
   "M";"MT" **\dir{-};
   "M";"XTM" **\dir{-}; ?(.75)*\dir{<};
   "R";"RT" **\dir{-};
     ?(0)*\dir{<};
   "R";"RB" **\dir{-};
   "LT";"MT" **\crv{(-12,22)&(0,22)}; ?(.93)*\dir{<};
   "XBM";"RB" **\crv{(1,-24)&(12,-24)};
   "XTM";"XTL" **\crv{ (-1,-9)& (-11,-9)};
   "XBM";"XBL" **\crv{ (1,-8)&(-12,-8)}; ?(.04)*\dir{>};
         (-15,8)*{\scriptstyle x};
         (-15,-4)*{\scriptstyle y};
         (3,8)*{\scriptstyle x};
         (3,-4)*{\scriptstyle y};
         (15,0)*{\scriptstyle y};
         (-15,-15)*{\scriptstyle y};
\endxy
\]
\[
\xy (-6,0)*{};
 (6,0)*{};
\endxy
\qquad = \qquad \xy (17,0)*{}; (-17,0)*{}; (-12,2)*{f};
(-12,2)*\xycircle(2.65,2.65){-}="L"; (0,2)*{\scriptstyle \inv f};
(0,2)*\xycircle(3.3,3.3){-}="M"; (12,0)="R"; (-12,6)*{}="LT";
(-12,-20)*{}="LB"; (0,6)*{}="MT"; (0,-6)*{}="MB"; (12,17)*{}="RT";
(12,-15)*{}="RB"; (0,-4)*{}="XTM"; (-12,-4)*{}="XTL";
(1,-15)*{}="XBM"; (-12,-15)*{}="XBL";
   "L";"LT" **\dir{-};
     ?(1)*\dir{<};
   "L";"XTL" **\dir{-}; ?(1)*\dir{>};
   "XBL";"LB" **\dir{-}; ?(.2)*\dir{>};
   "M";"MT" **\dir{-};
   "M";"XTM" **\dir{-}; ?(.75)*\dir{<};
   "R";"RT" **\dir{-};
     ?(0)*\dir{<};
   "R";"RB" **\dir{-};
   "LT";"MT" **\crv{(-12,22)&(0,22)}; ?(.93)*\dir{<};
   "XBM";"RB" **\crv{(1,-24)&(12,-24)};
   "XTM";"XTL" **\crv{ (-1,-9)& (-11,-9)};
   "XBM";"XBL" **\crv{ (1,-8)&(-12,-8)}; ?(.04)*\dir{>};
         (-15,8)*{\scriptstyle x};
         (-15,-4)*{\scriptstyle y};
         (-15,-15)*{\scriptstyle y};
         (3,8)*{\scriptstyle x};
         (3,-4)*{\scriptstyle y};
         (15,0)*{\scriptstyle y};
             (-17,-9)*{}="X1";
             (-17,-22)*{}="X4";
             (17,-9)*{}="X2";
             (17,-22)*{}="X3";
             "X1";"X2" **\dir{.};
             "X2";"X3" **\dir{.};
             "X3";"X4" **\dir{.};
             "X4";"X1" **\dir{.};
\endxy
\]
\[
\xy (-6,0)*{};
 (6,0)*{};
\endxy
\qquad = \qquad \xy (17,0)*{}; (-17,0)*{}; (-12,0)*{f};
(-12,0)*\xycircle(2.65,2.65){-}="L"; (0,0)*{\scriptstyle \inv f};
(0,0)*\xycircle(3.3,3.3){-}="M"; (12,0)="R"; (-12,6)*{}="LT";
(-12,-6)*{}="LB"; (0,6)*{}="MT"; (0,-6)*{}="MB"; (12,13)*{}="RT";
(12,-13)*{}="RB";
   "L";"LT" **\dir{-};
     ?(.75)*\dir{<};
   "L";"LB" **\dir{-};
     ?(1)*\dir{>};
   "M";"MT" **\dir{-}+(-6,10)*{ i_x};
   "M";"MB" **\dir{-}+(-6,-10)*{ i_y^{-1}};
     ?(.85)*\dir{<};
   "R";"RT" **\dir{-};
     ?(0)*\dir{<};
   "R";"RB" **\dir{-};
   "LT";"MT" **\crv{(-12,16)&(0,16)};
             ?(.97)*\dir{<};
   "LB";"MB" **\crv{(-12,-16)&(0,-16)};
         (-15,6)*{\scriptstyle x};
         (-15,-6)*{\scriptstyle y};
         (-2,6)*{\scriptstyle x};
         (-2,-6)*{\scriptstyle y};
         (15,0)*{\scriptstyle y};
\endxy
\]
\[
\xy (-6,0)*{};
 (6,0)*{};
\endxy
\qquad = \qquad \xy (17,0)*{}; (-17,0)*{}; (0,0)="R";
(0,11)*{}="RT"; (0,-11)*{}="RB";
   "R";"RT" **\dir{-};
     ?(0)*\dir{<};
   "R";"RB" **\dir{-};
         (3,0)*{\scriptstyle y};
\endxy
\]
\noindent The proof that $ff^{-1} = 1_x$ is similar, based on the
fact that $e$ is a natural isomorphism. \qed

Given a 2-category $K$ with finite products, we can also
define homomorphisms between coherent 2-groups in $K$, and
2-homomorphisms between these, by internalizing the
definitions of `weak monoidal functor' and `monoidal
natural transformation':

\begin{defn} \et
Given coherent 2-groups $C,C'$ in $K$, a {\bf homomorphism} $F
\maps C \to C'$ consists of:
\begin{itemize}
\item a morphism $F \maps C \to C'$
\end{itemize}
equipped with:
\begin{itemize}
\item a 2-isomorphism
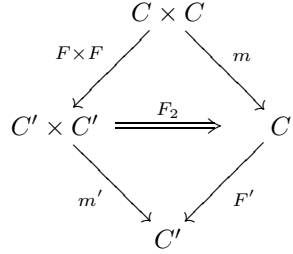
\[
\xy
 (0,15)*+{C \times C}="T";
 (-15,0)*+{C' \times C'}="L";
 (15,0)*+{C}="R";
 (0,-15)*+{C'}="B";
 (-7,0)*{}="ML";
 (7,0)*{}="MR";
     {\ar^{m} "T";"R"};
     {\ar_{F \times F} "T";"L"};
     {\ar^{F'} "R";"B"};
     {\ar_{m'} "L";"B"};
         {\ar@{=>}^{F_2} "ML";"MR"};
\endxy
\]
\item a 2-isomorphism
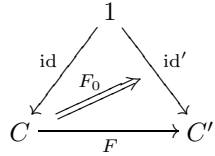
\[
 \xy
 (-12,0)*+{C}="L";
 (12,0)*+{C'}="R";
 (0,16)*+{1}="T";
    {\ar_{\id} "T";"L"};
    {\ar^{\id'} "T";"R"};
    {\ar_{F} "L";"R"};
    {\ar@{=>}^{F_0} (-7,2);(4,7)}
 \endxy
\]
\end{itemize}
such that diagrams commute expressing these laws:
\begin{itemize}
\item compatibility of $F_2$ with the associator:
\[ \makebox[0em]{
 \xy
   (-55,10)*+{(F \times F \times F)(m' \times 1)m'}="tl";
    (0,10)*+{(m \times 1)(F \times F)m'}="tm";
    (45,10)*+{(m \times 1)mF}="tr";
    (-55,-10)*+{(F \times F \times F)(1 \times m')m'}="bl";
    (0,-10)*+{(1 \times m)(F \times F)m'}="bm";
    (45,-10)*+{(1 \times m)mF}="br";
        {\ar^>>>>>>>>>>>{(F_2 \times F) \circ m'} "tl";"tm"};
        {\ar^>>>>>>>>>>>{(m \times 1) \circ F_2} "tm";"tr"};
        {\ar^{a \circ F} "tr";"br"};
        {\ar_{(F \times F \times F)\circ a} "tl";"bl"};
        {\ar_>>>>>>>>>>>{(F \times F_2) \circ m'} "bl";"bm"};
        {\ar_>>>>>>>>>>>{(1 \times m)\circ F_2} "bm";"br"};
 \endxy }
\]
\item compatibility of $F_0$ with the left unit law:
\[
 \xy
 (-25,8)*+{(\id' \times F)m'}="tl";
 (20,8)*+{F}="tr";
 (-25,-8)*+{(\id \times 1)(F \times F) m'}="bl";
 (20,-8)*+{(\id \times 1)m F}="br";
    {\ar^{F \circ \ell'} "tl";"tr"};
    {\ar_{\ell \circ F} "br";"tr"};
    {\ar_{(F_0 \times F) \circ m'} "tl";"bl"};
    {\ar_{(\id \times 1) \circ F_2} "bl";"br"};
 \endxy
\]
\item compatibility of $F_0$ with the right unit law:
\[
 \xy
 (-25,8)*+{(F \times \id')m'}="tl";
 (20,8)*+{F}="tr";
 (-25,-8)*+{(1 \times \id)(F \times F) m'}="bl";
 (20,-8)*+{(1 \times \id)m F}="br";
    {\ar^{F \circ r'} "tl";"tr"};
    {\ar_{r \circ F}"br";"tr"};
    {\ar_{(F \times F_0) \circ m'} "tl";"bl"};
    {\ar_{(1 \times \id) \circ F_2} "bl";"br"};
 \endxy
\]
\end{itemize}
\end{defn}

\begin{defn} \et
Given homomorphisms $F,G \maps C \to C'$ between coherent 2-groups
$C,C'$ in $K$, a {\bf 2-homomorphism} $\theta \maps F \To G$ is a
2-morphism such that the following diagrams commute:
\begin{itemize}
\item compatibility with $F_2$ and $G_2$:
\[
\xymatrix{
 (F \times F)m'
  \ar[rr]^{(\theta \times \theta)\circ m'}
  \ar[d]_{F_2}
&& (G \times G)m'
   \ar[d]^{G_2}     \\
  mF
   \ar[rr]^{m \circ \theta}
&&  mG }
\]
\item compatibility with $F_0$ and $G_0$:
\[
 \xy
 (-12,0)*+{\id F}="L";
 (12,0)*+{\id G}="R";
 (0,14)*+{\id'}="T";
    {\ar_{F_{0}} "T";"L"};
    {\ar^{G_{0}} "T";"R"};
    {\ar_{\id \circ \theta} "L";"R"};
 \endxy
\]
\end{itemize}
\end{defn}

\noindent It is straightforward to define a 2-category
\textbf{\textit{K}$\!$C2G} of
coherent 2-groups in $K$, homomorphisms between these, and
2-homomomorphisms between those.  We leave this to the reader,
who can also check that when $K = \Cat$, this 2-category
$K\cg$ reduces to $\cg$ as already defined.

To define concepts such as `topological 2-group', `Lie 2-group'
and `affine 2-group scheme' we need to consider coherent 2-group
objects in a special sort of 2-category which is defined by a
further process of internalization.  This is the 2-category of
`categories in $K$', where $K$ itself is a category.
A category in $K$ is usually called an `internal category'.  This
concept goes back to Ehresmann \cite{Ehresmann}, but a more
accessible treatment can be found in Borceux's handbook
\cite{Borceux}.  For completeness, we recall the definition here:

\begin{defn} \et Let $K$ be a category.
An {\bf internal category} or {\bf category in $K$}, say $X$,
consists of:
\begin{itemize}
\item an {\bf object of objects} $X_{0} \in K,$ \item an
{\bf object of morphisms} $X_{1} \in K,$
\end{itemize}
together with
\begin{itemize}
\item {\bf source} and {\bf target} morphisms $s,t \maps X_{1}
\rightarrow X_{0},$ \item a {\bf identity-assigning} morphism $i
\maps X_{0} \rightarrow X_{1},$ \item a {\bf composition} morphism
$\circ \maps X_{1} \times _{X_{0}} X_{1} \rightarrow X_{1}$
\end{itemize}
such that the following diagrams commute, expressing the usual
category laws:
\begin{itemize}
\item laws specifying the source and target of identity morphisms:
\[
\xymatrix{
 X_{0}
   \ar[r]^{i}
   \ar[dr]_{1}
   & X_{1}
   \ar[d]^{s} \\
  & X_{0} }
\hspace{.2in} \xymatrix{
   X_{0}
   \ar[r]^{i}
   \ar[dr]_{1}
   & X_{1}
   \ar[d]^{t} \\
  & X_{0}}
\]
\item laws specifying the source and target of composite
morphisms:
\[
\xymatrix{ X_{1} \times _{X_{0}} X_{1}
  \ar[rr]^{\circ}
  \ar[dd]_{p_{1}}
  && X_{1}
  \ar[dd]^{s} \\ \\
X_{1}
  \ar[rr]^{s}
  && X_{0} }
  \hspace{.2in}
\xymatrix{ X_{1} \times_{X_{0}} X_{1}
  \ar[rr]^{\circ}
  \ar[dd]_{p_{2}}
   && X_{1}
  \ar[dd]^{t} \\ \\
   X_{1}
  \ar[rr]^{t}
   && X_{0} }
\]
\item the associative law for composition of morphisms:
\[
\xymatrix{ X_{1} \times _{X_{0}} X_{1} \times _{X_{0}} X_{1}
  \ar[rr]^{\circ \times_{X_{0}} 1}
  \ar[dd]_{1 \times_{X_{0}} \circ}
   && X_{1} \times_{X_{0}} X_{1}
  \ar[dd]^{\circ} \\ \\
   X_{1} \times _{X_{0}} X_{1}
  \ar[rr]^{\circ}
   && X_{1} }
\]
\item the left and right unit laws for composition of morphisms:
\[
\xymatrix{ X_{0} \times _{X_{0}} X_{1}
  \ar[r]^{i \times 1}
  \ar[ddr]_{p_2}
   & X_{1} \times _{X_{0}} X_{1}
  \ar[dd]^{\circ}
   & X_{1} \times_{X_{0}} X_{0}
  \ar[l]_{1 \times i}
  \ar[ddl]^{p_1} \\ \\
   & X_{1} }
\]
\end{itemize}
\end{defn}

The pullbacks used in this definition should be obvious from the
usual definition of category; for example, composition should be
defined on pairs of morphisms such that the target of one is the
source of the next, and the object of such pairs is the pullback
$X_0 \times_{X_0} X_1$.  Notice that inherent to the definition is
the assumption that the pullbacks involved actually exist.  This
automatically holds if $K$ is a category with finite limits, but
there are some important examples like $K = \Diff$ where this is
not the case.

\begin{defn} \et Let $K$ be a category.
Given categories $X$ and $X'$ in $K$, an {\bf internal functor}
or {\bf functor in $K$}, say
$F \maps X \to X'$, consists of:
\begin{itemize}
\item a morphism $F_{0} \maps X_{0} \to X_{0}'$, \item a morphism
$F_{1} \maps X_{1} \rightarrow X_{1}'$
\end{itemize}
such that the following diagrams commute, corresponding to the
usual laws satisfied by a functor:
\begin{itemize}
\item preservation of source and target:
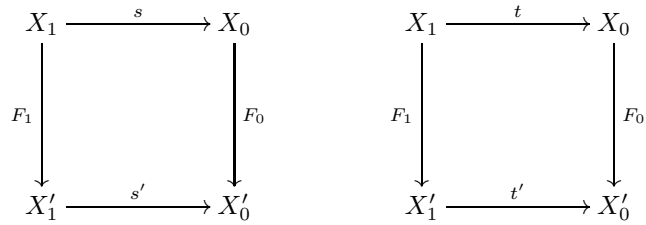
\[
\xymatrix{ X_{1} \ar[rr]^{s} \ar[dd]_{F_{1}}
 && X_{0}
\ar[dd]^{F_{0}} \\ \\
 X_{1}'
\ar[rr]^{s'}
 && X_{0}' }
\qquad \qquad \xymatrix{ X_{1} \ar[rr]^{t} \ar[dd]_{F_{1}}
 && X_{0}
\ar[dd]^{F_{0}} \\ \\
 X_{1}'
\ar[rr]^{t'}
 && X_{0}' }
\]
\item preservation of identity morphisms:
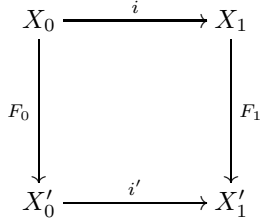
\[
\xymatrix{
 X_{0}
\ar[rr]^{i} \ar[dd]_{F_{0}}
 && X_{1}
\ar[dd]^{F_{1}} \\ \\
 X_{0}'
\ar[rr]^{i'}
 && X_{1}' }
\]
\item preservation of composite morphisms:
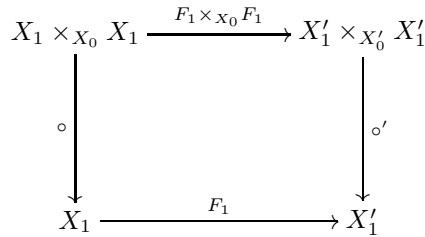
\[
\xymatrix{ X_{1} \times _{X_{0}} X_{1}
 \ar[rr]^{F_{1} \times_{X_0} F_{1}}
 \ar[dd]_{\circ}
  && X_{1}' \times_{X_{0}'} X_{1}'
 \ar[dd]^{\circ'} \\ \\
  X_{1}
 \ar[rr]^{F_{1}}
  && X_{1}' }
\]
\end{itemize}
\end{defn}

\begin{defn} \et Let $K$ be a
category.  Given categories $X, X'$ in $K$ and functors $F,G \maps
X \to X'$, an {\bf internal natural transformation}
or {\bf natural transformation in $K$}, say $\theta
\maps F \To G$, is a morphism $\theta \maps X_0 \to X'_1$ for
which the following diagrams commute, expressing the usual laws
satisfied by a natural transformation:
\begin{itemize}
\item laws specifying the source and target of the natural
transformation:
\[
 \xymatrix{X_0 \ar[dr]^F \ar[d]_{\theta} \\ X'_1 \ar[r]_s & X'_0 }
 \qquad \qquad
 \xymatrix{X_0 \ar[dr]^G \ar[d]_{\theta} \\ X'_1 \ar[r]_t & X'_0 }
\]
\item the commutative square law:
\[  \xymatrix{
X_1
 \ar[rr]^{\Delta (s\theta \times G)}
 \ar[dd]_{\Delta (F \times t\theta)}
  && X'_1 \times_{X_0} X'_1
 \ar[dd]^{\circ'} \\ \\
  X'_1 \times_{X_0} X'_1
 \ar[rr]^{\circ'}
  && X'_1
}
\]
\end{itemize}
\end{defn}

Given any category $K$, there is a strict 2-category $K\Cat$ whose
objects are categories in $K$, whose morphisms are functors in
$K$, and whose 2-morphisms are natural transformations in $K$. Of
course, a full statement of this result requires defining how to
compose functors in $K$, how to vertically and horizontally
compose natural transformations in $K$, and so on.  We shall not
do this here; the details can be found in Borceux's handbook
\cite{Borceux} or HDA6 \cite{HDA6}.

One can show that if $K$ is a category with finite products,
$K\Cat$ also has finite products.  This allows us to define
coherent 2-groups in $K\Cat$.  For example:

\begin{defn} \et
A {\bf topological category} is a category in $\Top$, the category
of topological spaces and continuous maps. A {\bf topological
2-group} is a coherent 2-group in $\Top\Cat$.
\end{defn}

\begin{defn} \label{lie.2-group} \et
A {\bf smooth category} is a category in $\Diff$, the category of
smooth manifolds and smooth maps. A {\bf Lie 2-group} is a
coherent 2-group in $\Diff\Cat$.
\end{defn}

\begin{defn} \et
An {\bf affine category scheme} is a category in $\CommRing^\op$,
the opposite of the category of commutative rings and ring
homomorphisms.  An {\bf affine 2-group scheme} is a coherent
2-group in $\CommRing^\op\Cat$.
\end{defn}

\noindent In the next section we shall give some examples of these
things.  For this, it sometimes handy to use an internalized
version of the theory of crossed modules.

As mentioned in the Introduction, a strict 2-group is essentially
the same thing as a {\bf crossed module}: a quadruple
$(G,H,t,\alpha)$ where $G$ and $H$ are groups, $t \maps H \to G$
is a homomorphism, and $\alpha \maps G \times H \to H$ is an
action of $G$ as automorphisms of $H$ such that $t$ is
$G$-equivariant:
\[   t( \alpha(g,h)) = g \, t(h)\, g^{-1} \]
and $t$ satisfies the so-called {\bf Peiffer identity}:
\[    \alpha(t(h),h') = hh'h^{-1} . \]
To obtain a crossed module from a strict 2-group $C$ we let $G =
C_0$, let $H = \ker s \subseteq C_1$, let $t \maps H \to G$ be the
restriction of the target map $t \maps C_1 \to C_0$ to $H$, and
set
\[       \alpha(g,h) = i(g)\, h \, i(g)^{-1}  \]
for all $g \in G$ and $h \in H$.  (In this formula multiplication
and inverses refer to the group structure of $H$, not composition
of morphisms in Conversely, we can build a strict 2-group from a
crossed module $(G,H,t,\alpha)$ as follows.  First we let $C_0 =
G$ and let $C_1$ be the semidirect product $H \rtimes G$ in which
multiplication is given by
\[    (h,g) (h',g') = (h \alpha(g,h'), gg') .\]
Then, we define source and target maps $s,t \maps C_1 \to C_0$ by:
\[   s(h,g) = g , \qquad   t(h,g) = t(h) g  ,\]
define the identity-assigning map $i \maps C_0 \to C_1$ by:
\[  i(g) = (1,g), \]
and define the composite of morphisms
\[  (h,g) \maps g  \to g'  , \qquad \quad
(h',g') \maps g' \to g'' \] to be:
\[   (hh',g) \maps g \to g'' .\]
For a proof that these constructions really work, see the
expository paper by Forrester-Barker \cite{FB}.

Here we would like to internalize these constructions in order to
build `strict 2-groups in $K\Cat$' from `crossed modules in $K$'
whenever $K$ is any category satisfying suitable conditions. Since
the details are similar to the usual case where $K = \Set$, we
shall be brief.

\begin{defn} \et A {\bf strict 2-group} in a 2-category with finite
products is a coherent 2-group in this 2-category such that
$a,i,e,l,r$ are all identity 2-morphisms --- or equivalently, a
group in the underlying category of this 2-category.
\end{defn}

\begin{defn} \et Given a category $K$ with finite products
and a group $G$ in $K$, an {\bf action} of $G$ on an object $X \in
K$ is a morphism $\alpha \maps G \times X \to X$ such that the
following diagrams commute:
\[
 \xymatrix{
 G \times G \times X
     \ar[rr]^{m \times 1_X}
     \ar[d]_{1_G \times \alpha}
 && G \times X    \ar[d]^{\alpha} \\
 G \times X     \ar[rr]^{\alpha} && X}
\]
\[
 \xymatrix{
 I \times X \ar[rr]^{\id \times 1_X} \ar[drr]_{\iso}
 && G \times X \ar[d]^{\alpha}\\
 && X }
\]
If $X$ is a group in $K$, we say $\alpha$ is an action of $G$ {\bf
as automorphisms of $X$} if this diagram also commutes:
\[
\xy
 (-50,8)*+{G \times X \times X}="L";
 (40,8)*+{X}="R";
 (-50,-8)*+{G \times G \times X \times X}="BL";
 (40,-8)*+{X \times X}="BR";
 (0,8)*+{G \times X}="T";
  (0,-8)*+{G \times X \times G \times X}="M";
     {\ar^{ 1_G \times m} "L";"T"};
     {\ar_{(\Delta_G \times 1_{X \times X})} "L";"BL"};
     {\ar^{ \alpha} "T";"R"};
      {\ar^{ (1_G \times S_{G,X} \times 1_X)} "BL";"M"};
     {\ar^{ \alpha \times \alpha} "M";"BR"};
     {\ar_{ m} "BR";"R"};
\endxy
\]
where $S_{G,X}$ stands for the `switch map' from $G \times X$ to
$X \times G$.

\end{defn}

\begin{defn} \label{internal.crossed.module}
\et  Given a category $K$ with finite products, a {\bf crossed
module in $K$} is a quadruple $(G,H,t,\alpha)$ with $G$ and $H$
being groups in $K$, $t \maps H \to G$ a homomorphism, and $\alpha
\maps G \times H \to H$ an action of $G$ as automorphisms of $H$,
such that diagrams commute expressing the $G$-equivariance of $t$:
\[
\xy (-34,8)*+{G \times H}="TL"; (0,8)*+{H}="TM"; (25,8)*+{G}="TR";
(-34,-8)*+{G \times H \times G }="BL"; (0,-8)*+{G \times G \times
G}="BM"; (25,-8)*+{G \times G}="BR"; {\ar_{\scs (\Delta_G \times
1_H)(1_G \times S_{G,H})} "TL";"BL"}; {\ar^{\scs \alpha}
"TL";"TM"}; {\ar^{\scs t} "TM";"TR"}; {\ar^<<<<<<<{\scs 1_G \times
t \times 1_G} "BL";"BM"}; {\ar^<<<<<{\scs m \times 1} "BM";"BR"};
{\ar_{\scs m} "BR";"TR"};
\endxy
\]
and the Peiffer identity:
\[
\xy (-34,8)*+{H \times H}="TL"; (25,8)*+{G \times H}="TR";
(-34,-8)*+{H \times H \times H }="BL"; (25,-8)*+{H}="BR";
(0,-8)*+{H \times H}="BM"; {\ar_{\scs (\Delta_H \times 1_H)(1_H
\times S_{H,H})} "TL";"BL"}; {\ar^{\scs t \times 1_H} "TL";"TR"};
{\ar^>>>>>>>>>>{\scs m \times \inv} "BL";"BM"}; {\ar^{\scs m}
"BM";"BR"}; {\ar^{\scs \alpha} "TR";"BR"};
\endxy
\]
\end{defn}

Next, consider a strict 2-group $C$ in the 2-category $K\Cat$,
where $K$ is a category with finite products.  This is the same as a
group in the underlying category of $K\Cat$.  By `commutativity of
internalization', this is the same as a category in $K\Grp$.  So,
$C$ consists of:
\begin{itemize}
\item a group $C_0$ in $K$, \item a group $C_1$ in $K$, \item {\bf
source} and {\bf target} homomorphisms $s,t \maps C_1 \to C_0$,
\item an {\bf identity-assigning} homomorphism $i \maps C_0 \to
C_1$, \item a {\bf composition} homomorphism $\circ \maps C_{1}
\times_{C_{0}} C_{1} \rightarrow C_{1}$
\end{itemize}
such that the usual laws for a category hold:
\begin{itemize}
\item laws specifying the source and target of identity morphisms,
\item laws specifying the source and target of composite
morphisms,
\item the associative law for composition,
\item the left and right unit laws for composition of morphisms.
\end{itemize}
We shall use this viewpoint in the following:

\begin{prop} \label{crossed.module} \et Let $K$ be a category with
finite products such that $K\Grp$ has finite limits.  Given a strict
2-group $C$ in $K\Cat$, there is a crossed module $(G,H,t,\alpha)$
in $K$ such that
\[           G = C_0  , \qquad          H = \ker s , \]
such that
\[           t \maps H \to G  , \]
is the restriction of $t \maps C_1 \to C_0$ to the subobject $H$,
and such that
\[           \alpha \maps G \times H \to H   \]
makes this diagram commute:
\[
\xymatrix{ G \times H
 \ar[rrr]^{\alpha} \ar[d]_{(\Delta \times 1_H)(1_G \times S_{H,G})} &&& H \\
G \times H \times G \ar[d]_{i \times 1_H \times i} &&& H \times H \ar[u]_{m} \\
H \times H \times H \ar[rrr]^{1_G \times 1_H \times \inv_H}
 &&& H \times H \times H \ar[u]_{m \times 1_H}
}
\]

Conversely, given a crossed module $(G,H,t,\alpha)$ in $K$, there
is a strict 2-group $C$ in $K$ for which
\[          C_0 = G  \]
and for which
\[          C_1 = H \times G \]
is made into a group in $K$ by taking the semidirect product using
the action $\alpha$ of $G$ as automorphisms on $H$. In this strict
2-group we define source and target maps $s,t \maps C_1 \to C_0$
so that these diagrams commute:
\[
\xymatrix{
 H \times G \ar[r]^{s} \ar[d]_{1_{H \times G}}& G \\
 H \times G \ar[ur]_{\pi_2}
}
\]

\[
\xymatrix{
H \times G \ar[r]^{t} \ar[d]_{t \times 1_G}& G \\
G \times G \ar[ur]_{m} }
\]
define the identity-assigning map $\id \maps C_0 \to C_1$ so that
this diagram commutes:
\[
\xymatrix{
G \ar[r]^{\id} \ar[d]_{\iso}& H \times G \\
I \times G \ar[ur]_{i_H \times 1_G} }
\]
and define composition $\circ \maps C_1 \times_{C_0} C_1 \to C_1$
such that this commutes:
\[
\xymatrix{ (H \times G) \times_{\scriptscriptstyle G} (H \times G)
 \ar[rr]^{\circ} \ar[d]_{\pi_{123}}
&& H \times G \\
H \times G \times H \ar[rr]_{1_H \times S_{G,H}} && H \times H
\times G \ar[u]_{m \times 1_G} }
\]
where $\pi_{123}$ projects onto the product of the first, second
and third factors.
\end{prop}

\textbf{Proof. }  The proof is modeled directly after the case $K
= \Set$; in particular, the rather longwinded formula for $\alpha$
reduces to
\[       \alpha(g,h) = i(g)\, h \, i(g)^{-1}  \]
in this case.  Note that to define $\ker s$ we need $K\Grp$ to
have finite limits, while to define $C_1$ and make it into a group
in $K$, we need $K$ to have finite products. \hbox{\hskip 30em}
\qed

When the category $K$ satisfies the hypotheses of this
proposition, one can go further and show that strict 2-groups in
$K\Cat$ are indeed `the same' as crossed modules in $K$. To do
this, one should first construct a 2-category of strict 2-groups
in $K\Cat$ and a 2-category of crossed modules in $K$, and then
prove these 2-categories are equivalent.  We leave this as an
exercise for the diligent reader.

\section{Examples} \label{examplesection}

\subsection{Automorphism 2-groups}

Just as groups arise most naturally from the consideration of {\it
symmetries}, so do 2-groups. The most basic example of a group is
a permutation group, or in other words, the automorphism group of
a set. Similarly, the most basic example of a 2-group consists of
the automorphism group of a category.  More generally, we can talk
about the automorphism group of an object in any category.
Likewise, we can talk about the `automorphism 2-group' of an
object in any 2-category.

We can make this idea precise in somewhat different ways depending
on whether we want a strict, weak, or coherent 2-group.  So, let
us consider various sorts of `automorphism 2-group' for an object
$x$ in a 2-category $K$.

The simplest sort of automorphism 2-group is a strict one:

\begin{example} \et {\rm
For any strict 2-category $K$ and object $x \in K$ there is a
strict 2-group $\Aut_s(x)$, the {\bf strict automorphism 2-group}
of $x$.  The objects of this 2-group are isomorphisms $f \maps x
\to x$, while the morphisms are 2-isomorphisms between these.
Multiplication in this 2-group comes from composition of morphisms
and horizontal composition of 2-morphisms.   The identity object
$1 \in \Aut_s(x)$ is the identity morphism $1_x \maps x \to x$.
}\end{example}

To see what this construction really amounts to, take $K = \Cat$
and let $M \in K$ be a category with one object. A category with
one object is secretly just a monoid, with the morphisms of the
category corresponding to the elements of the monoid.   An
isomorphism $f \maps M \to M$ is just an automorphism of this
monoid.  Given isomorphisms $f,f' \maps M \to M$, a 2-isomorphism
from $f$ to $f'$ is just an invertible element of the monoid, say
$\alpha$, with the property that $f$ conjugated by $\alpha$ gives
$f'$:
\[       f'(m) = \alpha^{-1}f(m)\alpha   \]
for all elements $m \in M$.  This is just the usual commuting
square law in the definition of a natural isomorphism, applied to
a category with one object.  So, $\Aut_s(M)$ is a strict 2-group
that has automorphisms of $x$ as its objects and `conjugations' as
its morphisms.

Of course the automorphisms of a monoid are its symmetries in the
classic sense, and these form a traditional group.  The new
feature of the automorphism 2-group is that it keeps track of the
{\it symmetries between symmetries}: the conjugations carrying one
automorphism to another.  More generally, in an $n$-group, we
would keep track of symmetries between symmetries between
symmetries between... and so on to the $n$th degree.

The example we are discussing here is especially well-known when
the monoid is actually a group, in which case its automorphism
2-group plays an important role in nonabelian cohomology and the
theory of nonabelian gerbes \cite{Breen,Breen2,Giraud}. In fact,
given a group $G$, people often prefer to work, not with
$\Aut_s(G)$, but with a certain weak 2-group that is equivalent to
$\Aut_s(G)$ as an object of $\wg$.  The objects of this group are
called `$G$-bitorsors'.  They are worth understanding, in part
because they illustrate how quite different-looking weak 2-groups
can actually be equivalent.

Given a group $G$, a \textbf{\textit{G}-bitorsor} $X$ is a set
with commuting left and right actions of $G$, both of which are
free and transitive.  We write these actions as $g \cdot x$ and $x
\cdot g$, respectively.   A morphism between $G$-bitorsors $f
\maps X \to Y$ is a map which is equivariant with respect to both
these actions.  The tensor product of $G$-bitorsors $X$ and $Y$ is
defined to be the space
\[   X \tensor Y = X \times Y/((x\cdot g,y) \sim (x,g\cdot y)) , \]
which inherits a left $G$-action from $X$ and a right $G$-action
from $Y$.  It is easy to check that $X \tensor Y$ is again a
bitorsor.  Accompanying this tensor product of bitorsors there is
an obvious notion of the tensor product of morphisms between
bitorsors, making $G$-bitorsors into a weak monoidal category
which we call \textbf{\textit{G}-Bitors}.

The identity object of $G$-$\Bitors$ is just $G$, with its
standard left and right action on itself.  This is the most
obvious example of a $G$-bitorsor, but we can get others as
follows. Suppose that $f \maps G \to G$ is any automorphism.  Then
we can define a $G$-bitorsor $G_f$ whose underlying set is $G$,
equipped with the standard left action of $G$:
\[   g \cdot h = gh , \qquad \qquad g \in G, h \in G_f \]
but with the right action twisted by the automorphism $f$:
\[   h \cdot g = hf(g), \qquad \qquad g \in G, h \in G_f .\]

The following facts are easy to check.  First, every $G$-bitorsor
is isomorphic to one of the form $G_f$.  Second, every morphism
from $G_f$ to $G_{f'}$ is of the form
\[     h \mapsto h \alpha \]
for some $\alpha \in G$ with
\[       f'(g) = \alpha^{-1}f(g)\alpha   \]
for all $g \in G$.  Third, the tensor product of $G_f$ and
$G_{f'}$ is isomorphic to $G_{ff'}$.

With the help of these facts, one can show that $G$-$\Bitors$ is
equivalent as a weak monoidal category to $\Aut_s(G)$.  The point
is that the objects of $G$-$\Bitors$ all correspond, up to
isomorphism, to the objects of $\Aut_s(G)$: namely, automorphisms
of $G$. Similarly, the morphisms of $G$-$\Bitors$ all correspond
to the morphisms of $\Aut_s(G)$: namely, `conjugations'.  The
tensor products agree as well, up to isomorphism.

Since $\Aut_s(G)$ is a strict 2-group, it is certainly a weak one
as well.  Since $G$-$\Bitors$ is equivalent to $\Aut_s(G)$ as a
weak monoidal category, it too is a weak 2-group, and it is
equivalent to $\Aut_s(G)$ as an object of the 2-category $\wg$.

In this particular example, the `strict automorphism 2-group'
construction seems quite useful. But for some applications, this
construction is overly strict.  First, we may be interested in
automorphism 2-group of an object in a weak 2-category
(bicategory), rather than a strict one.  Second, given objects
$x,y$ in a weak 2-category $K$, it is often unwise to focus
attention on the isomorphisms $f \maps x \to y$.  A more robust
concept is that of a {\bf weakly invertible morphism}: a morphism
$f \maps x \to y$ for which there exists a morphism $\bar{f} \maps
y \to x$ and 2-isomorphisms $\iota \maps 1_x \To f\bar{f}$,
$\epsilon \maps \bar{f}f \To 1_y$. Using weakly invertible
morphisms as a substitute for isomorphisms gives a weak version of
the automorphism 2-group:

\begin{example} \et {\rm
For any weak 2-category $K$ and object $x \in K$ there is a weak
2-group $\Aut_w(x)$, the {\bf weak automorphism 2-group} of $x$.
The objects of this 2-group are weakly invertible morphisms $f
\maps x \to x$, while the morphisms are 2-isomorphisms between
these.  Multiplication in this 2-group comes from composition of
morphisms and horizontal composition of 2-morphisms.  The identity
object $1 \in \Aut_w(C)$ is the identity functor. }\end{example}

A weakly invertible morphism $f \maps x \to y$ is sometimes called
an `equivalence'.  Here we prefer to define an {\bf equivalence}
from $x$ to $y$ to be a morphism $f \maps x \to y$ with a {\it
specified} weak inverse $\bar{f} \maps y \to x$ and {\it
specified} 2-isomorphisms $\iota_f \maps 1_x \To f\bar{f}$,
$\epsilon_f \maps \bar{f}f \To 1_y$.  An equivalence from $x$ to
$y$ is thus a quadruple $(f,\bar{f},\iota_f,\epsilon_f)$. We can
make a coherent 2-group whose objects are equivalences from $x$ to
itself:

\begin{example} \label{autoequivalence} \et {\rm
For any weak 2-category $K$ and object $x \in K$ there is a
coherent 2-group $\Aut_{eq}(x)$, the {\bf autoequivalence 2-group}
of $x$.  The objects of $\Aut_{eq}(x)$ are equivalences from $x$
to $x$. A morphism in $\Aut_w(x)$ from
$(f,\bar{f},\iota_f,\epsilon_f)$ to
$(g,\bar{g},\iota_g,\epsilon_g)$ consists of 2-isomorphisms
\[  \alpha \maps f \To g , \qquad \bar{\alpha} \maps
\bar{f} \To \bar{g}   \] such that the following diagrams commute:
\[
 \xy
 (-10,0)*+{1}="L";
 (0,10)*+{f\bar{f}}="T";
 (0,-10)*+{g\bar{g}}="B";
    {\ar^{\iota_f} "L";"T"};
    {\ar_{\iota_g} "L";"B"};
    {\ar^{\alpha\cdot\bar{\alpha}} "T";"B"};
 \endxy
 \qquad \qquad
  \xy
 (10,0)*+{1}="L";
 (0,10)*+{\bar{f}f}="T";
 (0,-10)*+{\bar{g}g}="B";
    {\ar^{\epsilon_f} "T";"L"};
    {\ar_{\epsilon_g} "B";"L"};
    {\ar_{\bar{\alpha}\cdot\alpha} "T";"B"};
 \endxy
\]
Multiplication in this 2-group comes from the standard way of
composing equivalences, together with horizontal composition of
2-morphisms.  The identity object $1 \in \Aut_{eq}(x)$ is the
equivalence $(1_x,1_x,1_{1_x},1_{1_x})$.

One can check that $\Aut_{eq}(x)$ is a weak 2-group, because every
object $F = (f,\bar{f},\iota_f,\epsilon_f)$ of $\Aut_{eq}(x)$ has
the weak inverse $\bar{F} =
(\bar{f},f,\epsilon_f^{-1},\iota_f^{-1})$. But in fact, the proof
of this involves constructing isomorphisms
\[    i_F \maps 1_x \To F\bar{F}, \qquad
   e_F \maps \bar{F}F \To 1_x   \]
from the data at hand, and these isomorphisms can easily be chosen
to satisfy the zig-zag identities, so $\Aut_{eq}(x)$ actually
becomes a coherent 2-group. }\end{example}

An equivalence $(f,\bar{f},\iota_f,\epsilon_f)$ is an {\bf adjoint
equivalence} if it satisfies the zig-zag identities. We can also
construct a coherent 2-group whose objects are adjoint
equivalences from $x$ to itself:

\begin{example} \et {\rm For any weak 2-category $K$ and object $x \in K$
there is a coherent 2-group $\Aut_{ad}(x)$, the {\bf adjoint
autoequivalence group} of $x$. The objects of this 2-group are
adjoint equivalences from $x$ to $x$, while the morphisms are
defined as in $\Aut_{eq}(x)$.  Multiplication in this 2-group
comes from the usual way of composing equivalences (using the fact
that composite of adjoint equivalences is again an adjoint
equivalence) together with horizontal composition of 2-morphisms.
The identity object $1 \in \Aut_{ad}(x)$ is the adjoint
equivalence $(1_x,1_x,1_{1_x},1_{1_x})$. $\Aut_{ad}(x)$ becomes a
coherent 2-group using the fact that every object $F$ of
$\Aut_{ad}(x)$ becomes part of an adjunction $(F,\bar{F},i_F,e_F)$
as in Example \ref{autoequivalence}. }\end{example}

\subsection{The fundamental 2-group}

Another source of 2-groups is topology: for any topological space
$X$ and any point $x \in X$ there is a coherent 2-group
$\Pi_2(X,x)$ called the `fundamental 2-group' of $X$ based at $x$.
The fundamental 2-group is actually a watered-down version of what
Hardie, Kamps and Kieboom \cite{HKK} call the `homotopy
bigroupoid' of $X$, denoted by $\Pi_2(X)$.  This is a weak
2-category whose objects are the points of $X$, whose morphisms
are paths in $X$, and whose 2-morphisms are homotopy classes of
paths-of-paths. More precisely, a morphism $f \maps x \to y$ is a
map $f \maps [0,1] \to X$ with $f(0) = x$ and $f(1) = y$, while a
2-morphism $\alpha\maps f \To g$ is an equivalence class of maps
$\alpha\maps [0,1]^2 \to X$ with
\[
\begin{array}{ccl}
\alpha(s,0) &=& f(s) \\
\alpha(s,1) &=& g(s) \\
\alpha(0,t) &=& x \\
\alpha(1,t) &=& y
\end{array}
\]
for all $s,t \in [0,1]$, where the equivalence relation is that
$\alpha\sim \alpha'$ if there is a map $H \maps [0,1]^3 \to X$
with
\[
\begin{array}{ccl}
H(s,t,0) &=& \alpha(s,t)  \\
H(s,t,1) &=& \alpha'(s,t)  \\
H(s,0,u) &=& f(s)  \\
H(s,1,u) &=& g(s) \\
H(0,t,u) &=& x \\
H(1,t,u) &=& y
\end{array}
\]
for all $s,t,u \in [0,1]$.  This becomes a weak 2-category in a
natural way, with composition of paths giving composition of
morphisms, and two ways of composing paths-of-paths giving
vertical and horizontal composition of 2-morphisms:
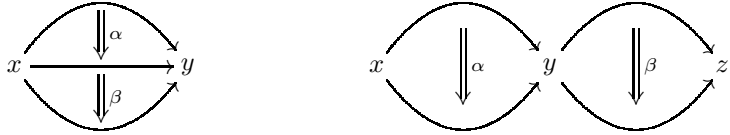
\[
\xymatrix{
  x
   \ar@/^2pc/[rr]^{}_{}="0"
   \ar[rr]^{}="1"
   \ar@/_2pc/[rr]_{}_{}="2"
   \ar@{=>}"0";"1" ^{\alpha}
   \ar@{=>}"1";"2" ^{\beta}
&&  y } \qquad \qquad \qquad \xymatrix{
 x
   \ar@/^2pc/[rr]^{}_{}="0"
   \ar@/_2pc/[rr]_{}_{}="2"
   \ar@{}"0";"2"|(.2){\,}="7"
  \ar@{}"0";"2"|(.8){\,}="8"
   \ar@{=>}"7";"8" ^{\alpha}          
&&  y
   \ar@/^2pc/[rr]^{}_{}="0"
   \ar@/_2pc/[rr]_{}="2"
   \ar@{}"0";"2"|(.2){\,}="7"
  \ar@{}"0";"2"|(.8){\,}="8"
   \ar@{=>}"7";"8" ^{\beta}
&&  z }
\]

Hardie, Kamps and Kieboom show that every 2-morphism in $\Pi_2(X)$
is invertible, and they construct an adjoint equivalence
$(f,\bar{f},\iota_f,\epsilon_f)$ for every morphism $f$ in
$\Pi_2(X)$.  This is why they call $\Pi_2(X)$ a `bigroupoid'. One
might also call this a `coherent 2-groupoid', since such a thing
with one object is precisely a coherent 2-group.  Regardless of
the terminology, this implies that for any point $x \in X$ there
is a coherent 2-group whose objects are morphisms $f \maps x \to
x$ in $\Pi_2(X)$, and whose morphisms are those 2-morphisms
$\alpha \maps f \To g$ in $\Pi_2(X)$ for which $f,g \maps x \to
x$.   We call this coherent 2-group the {\bf fundamental 2-group}
$\Pi_2(X,x)$.

In fact, a fundamental 2-group is a special case of an
`autoequivalence 2-group', as defined in Example
\ref{autoequivalence}.   A point $x \in X$ is an object of the
weak 2-category $\Pi_2(X)$, and the autoequivalence 2-group of
this object is precisely the fundamental 2-group $\Pi_2(X,x)$.
Even better, we can turn this idea around: there is a way to see
any autoequivalence 2-group as the fundamental 2-group of some
space, at least up to equivalence!  Unfortunately, proving this
fact would take us too far out of our way here.  However, the
relation between 2-groups and topology is so important that we
should at least sketch the basic idea.

Suppose $K$ is a weak 2-category, and let $K_0$ be its underlying
coherent 2-groupoid --- that is, the weak 2-category with the same
objects as $K$, but with the adjoint equivalences in $K$ as its
morphisms and the invertible 2-morphisms of $K$ as its
2-morphisms.  Let $|K_0|$ be the geometric realization of the
nerve of $K_0$ as defined by Duskin \cite{Duskin}.  Then any
object $x \in K$ gives a point $x \in |K_0|$, and the
autoequivalence 2-group $\Aut_{eq}(x)$ is equivalent to
$\Pi_2(|K_0|,x)$.

In fact, something much stronger than this should be true.
According to current thinking on $n$-categories and homotopy
theory \cite{BD}, 2-groups should really be `the same' as
connected pointed homotopy 2-types.  For example, we should be
able to construct a 2-category ${\rm Conn2Typ_*}$ having connected
pointed CW complexes with $\pi_n = 0$ for $n > 2$ as objects,
continuous basepoint-preserving maps as morphisms, and homotopy
classes of basepoint-preserving homotopies between such maps as
2-morphisms.  The fundamental 2-group construction should give a
2-functor:
\[
\begin{array}{rcl}
\Pi_2 \maps {\rm Conn2Typ_*} &\to& \cg  \\
              (X,x) &\mapsto & \Pi_2(X,x)
\end{array}
\]
while the geometric realization of the nerve should give a
2-functor going the other way:
\[
\begin{array}{rcl}
\Pi_2^{-1} \maps \cg &\to& {\rm Conn2Typ_*}  \\
            C &\mapsto & (|C|, 1)
\end{array}
\]
and these should extend to a biequivalence of 2-categories.
To the best of our knowledge, nobody has yet written up a proof of
this result.  However, a closely related higher--dimensional result 
has been shown by C.\ Berger \cite{Berger}: the model category
of homotopy 3-types is Quillen equivalent to a suitably defined
model category of weak 3-groupoids.  

\subsection{Classifying 2-groups using group cohomology}
\label{classificationsection}

In this section we sketch how a coherent 2-group is determined, up
to equivalence, by four pieces of data:
\begin{itemize}
\item a group $G$,
\item an abelian group $H$,
\item an action $\alpha$ of $G$ as automorphisms of $H$,
\item an element $[a]$ of the cohomology group $H^3(G,H)$,
\end{itemize}
where the last item comes from the associator.  Various versions
of this result have been known to experts at least
since Sinh's thesis \cite{Sinh}, but since this thesis was
unpublished they seem to have spread largely in the form of
`folk theorems'.  A very elegant treatment can be found in the 1986
draft of Joyal and Street's paper on braided tensor categories
\cite{JS}, but not in the version that was finally published in 1993.
So, it seems worthwhile to provide a precise statement and proof
here.

One way to prove this result would be to take a detour through
topology.  Using the ideas sketched at end of the previous section,
equivalence classes of coherent 2-groups should be in one-to-one
correspondence with homotopy equivalence classes of connected
pointed CW complexes having homotopy groups that vanish above
$\pi_2$.  The latter, in turn, can be classified in terms of their
`Postnikov data': the group $G = \pi_1$, the abelian group $H =
\pi_2$, the action of $\pi_1$ on $\pi_2$, and the Postnikov
$k$-invariant, which is an element of $H^3(\pi_1,\pi_2)$.  The
advantage of this approach is that it generalizes to $n$-groups
for higher $n$, and clarifies their relation to topology.  The
disadvantage is that it is indirect and relies on results that
themselves take some work to prove.  Besides the relation between
coherent 2-groups and homotopy 2-types, one needs the theory of
Postnikov towers in the case where $\pi_1$ acts nontrivially on the
higher homotopy groups \cite{Robinson}.

To avoid all this, we take a more self-contained approach.  First
we show that every coherent 2-group is equivalent to a `special' one:

\begin{defn} \et A coherent 2-group is {\bf skeletal} if its
underlying category is skeletal: that is, if any pair of isomorphic
objects in this category are equal.
\end{defn}

\begin{defn} \et A {\bf special 2-group} is a skeletal coherent 2-group
such that the left unit law $\ell$, the right unit law $r$,
the unit $i$ and the counit $e$ are identity natural transformations.
\end{defn}

\noindent
We then show that any special 2-group determines a quadruple
$(G,H,\alpha,a)$.  The objects of a special 2-group form a
group $G$.  The automorphisms of the unit object form an abelian
group $H$.  There is an action $\alpha$ of $G$ on $H$, defined just as
in the construction of a crossed module from a strict 2-group.  The
associator gives rise to a map $a \maps G^3 \to H$.  Furthermore, the
pentagon identity and other properties of monoidal categories imply
that $a$ is a `normalized 3-cocycle' on $G$ with values in the
$G$-module $H$.  When we work through this in detail, it will also
become clear that conversely, any quadruple $(G,H,\alpha,a)$ of this
sort determines a special 2-group.

Following Joyal and Street, we exploit these results by
constructing a 2-category of special 2-groups that is
biequivalent to $\cg$, for which not only the objects but also
the morphisms and 2-morphisms can be described using group cohomology.
As a corollary, it will follow that coherent 2-groups are classified
up to equivalence by quadruples $(G,H,\alpha,[a])$, where
$[a] \in H^3(G,H)$ is the cohomology class of the 3-cocycle $a$.

We begin by proving the following fact:

\begin{prop} \label{equivalence} \et Every coherent 2-group is
equivalent in $\cg$ to a special 2-group.
\end{prop}

\textbf{Proof. } First suppose that $C$ is a coherent 2-group.
Note that $C$ is equivalent, as an object of $\cg$, to a skeletal
coherent 2-group.  To see this, recall that every category is
equivalent to a skeletal one: we can take this to be any full
subcategory whose objects include precisely one representative
from each isomorphism class.  Using such an equivalence of
categories, we can transfer the coherent 2-group structure from
$C$ to a skeletal category $C_0$.  It is then routine to check
that $C$ and $C_0$ are equivalent as objects of $\cg$.

Next suppose $C$ is a skeletal coherent 2-group.  We shall
construct a special 2-group $\tilde{C}$ that is equivalent to
$C$.  As a category, $\tilde{C}$ will be precisely the same as $C$,
so it will still be skeletal. However, we shall adjust the tensor
product, left and right unit laws, unit and counit, and associator
to ensure that $\tilde{\ell}, \tilde{r}, \tilde{\i}$ and
$\tilde{e}$ are identity natural transformations.  We do this
using the following lemma:

\begin{lem} \et If $C$ is a coherent 2-group, and for each $x,y \in C$
we choose an isomorphism $\gamma_{x,y} \maps x \tilde{\tensor} y
\to x \tensor y$ for some object $x \tilde{\tensor} y \in C$, then
there exists a unique way to make the underlying category of $C$
into a coherent 2-group $\tilde C$ such that:
\begin{enumerate}
\item the tensor product of any pair of objects $x,y$ in
$\tilde{C}$ is $x \tilde{\tensor} y$, \item there is a
homomorphism of coherent 2-groups $F \maps C \to \tilde{C}$ whose
underlying functor is the identity, for which $F_0$ and $F_{-1}$
are the identity, and for which
\[           (F_2)_{x,y} = \gamma_{x,y}   \]
for every $x,y \in C$.
\end{enumerate}
Moreover, $F \maps C \to \tilde{C}$ is an equivalence in $\cg$.
\end{lem}

\textbf{Proof. }  First we show uniqueness.  The tensor product of
objects in $\tilde{C}$ is determined by item 1.  For $F$ as in
item 2 to be a weak monoidal functor we need $F_2$ to be natural,
so the tensor product $f \tilde{\tensor} g$ of morphisms $f \maps
x \to x'$, $g \maps y \to y'$ in $\tilde{C}$ is determined by the
requirement that
\[
\xymatrix{
 x \tilde{\ten} y
   \ar[r]^{\gamma_{x,y}}
   \ar[d]_{f \tilde{\ten} g}
&   x \ten y
  \ar[d]^{f \ten g}
  \\
 x' \tilde{\ten} y'
   \ar[r]^{\gamma_{x',y'}}
&  x' \ten y' }
\]
commute.  The unit object of $\tilde{C}$ must be the same as that
of $C$, since $F_0$ is the identity.   The unit $\tilde{\i}$ and
counit $\tilde{e}$ of $\tilde C$ are determined by the coherence
laws \textbf{H1} and \textbf{H2} in Section
\ref{preservationsection}.  The associator $\tilde{a}$ of $\tilde
C$ is determined by this coherence law in the definition of `weak
monoidal functor':
\[
\xymatrix@!C{
 (F(x) \tilde{\ten} F(y)) \tilde{\ten} F(z)
   \ar[r]^>>>>>>>{F_{2} \tilde{\ten} 1}
   \ar[d]_{\tilde{a}_{F(x), F(y), F(z)}}
& F(x \ten y) \tilde{\ten} F(z)
   \ar[r]^{F_{2}}
& F((x \ten y) \ten z)
   \ar[d]^{F(a_{x,y,z})}   \\
 F(x) \tilde{\ten} (F(y) \ten F(z))
   \ar[r]^>>>>>>>{1 \tilde{\ten} F_{2}}
& F(x) \tilde{\ten} F(y \ten z)
   \ar[r]^{F_{2}}
& F(x \ten (y \ten z))
 }
\]
Similarly, the left and right unit laws $\tilde{l},\tilde{r}$ of
$\tilde C$ are determined by the other two coherence laws in this
definition:
\[
\xymatrix{
 1 \tilde{\ten} F(x)
   \ar[r]^{\tilde{\ell}_{F(x)}}
   \ar[d]_{F_{0} \tilde{\ten} 1}
&  F(x)  \\
 F(1) \tilde{\ten} F(x)
   \ar[r]^{F_{2}}
&  F(1 \ten x)
   \ar[u]_{F(\ell_{x})}
}
\]
\vskip 1em
\[
\xymatrix{
 F(x) \tilde{\ten} 1
   \ar[r]^{\tilde{r}_{F(x)}}
   \ar[d]_{1 \tilde{\ten} F_{0}}
&  F(x)    \\
 F(x) \tilde{\ten} F(1)
   \ar[r]^{F_{2}}
&  F(x \ten 1)
   \ar[u]_{F(r_x)}
}
\]
It is then an exercise to check that with these choices,
$\tilde{C}$ really does become a coherent 2-group, that $F \maps C
\to \tilde{C}$ is a homomorphism, and that $F$ is an equivalence
of coherent 2-groups. \qed

We now apply this lemma, taking
$$         x \tilde{\tensor} y =
\left\{
\begin{array}{cl}
y        & {\rm if \; } x = 1 \\
x        & {\rm if \; } y = 1 \\
1        & {\rm if \; } x = \yb \\
1        & {\rm if \; } y = \xb \\
x \ten y & {\rm otherwise}
\end{array}   \right.
\]
and taking
\[         \gamma_{x,y} =
\left\{ \begin{array}{cl}
\ell_y^{-1} & {\rm if \; } x = 1 \\
r_x^{-1}    & {\rm if \; } y = 1 \\
e_y^{-1}    & {\rm if \; } x = \yb \\
i_x         & {\rm if \; } y = \xb \\
1_{x \ten y}& {\rm otherwise}
\end{array}  \right.
\]
Calculations then show that with these choices, $\tilde{\ell},
\tilde{r}, \tilde{\i}$ and $\tilde{e}$ are identity natural
transformations.   For example, to show that $\tilde{\i}$ is the
identity we use coherence law \textbf{H1}, which says this diagram
commutes:
\[
 \\
 \vcenter{
 \xymatrix@!C{
    x \tilde{\ten} \overline{F(x)}
    \ar[r]^<<<<<{1 \tilde{\ten} F_{-1}}
  & F(x) \tilde{\ten} F(\xb)
       \ar[r]^{F_{2}}
  & F(x \ten \xb)   \\
     1
       \ar[u]^{\tilde{\i}_{F(x)}}
       \ar[rr]^{F_{0}}
  && F(1)
       \ar[u]_{F(i_x)}
  }}
\]
By the definition of $F, F_0$, and $F_2 = \gamma$ together with
the fact that $x \tilde{\ten} \xb = 1$, this diagram reduces to
\[
 \\
 \vcenter{
 \xymatrix@!C{
    1
    \ar[r]^<<<<<{1_1}
  & 1
       \ar[r]^{i_x}
  & x \ten \xb   \\
     1
       \ar[u]^{\tilde{\i}_x}
       \ar[rr]^{1_1}
  && 1
       \ar[u]_{i_x}
  }}
\]
which implies that $\tilde{\i}_x = 1_1$.  Similarly, to show that
$\tilde{e}_x$ is the identity we use \textbf{H2}, and to show
$\tilde{\ell}_x$ and $\tilde{r}_x$ are identities we use the
coherence laws for the left and right unit laws in the definition
of `weak monoidal functor'.  \qed

We now describe in a bit more detail how to get a quadruple
$(G,H,\alpha,a)$ from a special 2-group $C$.
In general, the objects of a 2-group need not form a group under
multiplication, since we only have isomorphisms
\[    (x \tensor y) \tensor z \iso x \tensor (y \tensor z) , \]
\[          1 \tensor x \iso x, \qquad  x \tensor 1 \iso x  ,\]
\[          \xb \tensor x \iso 1, \qquad  x \tensor \xb \iso 1  .\]
However, in a special 2-group, isomorphic objects are equal,
so the objects form a group.  This is our group $G$.

The Eckmann-Hilton argument shows that in any weak monoidal
category, the endomorphisms of the unit object form a commutative
monoid under tensor product or, what is the same, composition:
\[
\begin{array}{ccl}
h \tensor h' &=& (h 1_1) \tensor (1_1 h') \\
          &=& (h \tensor 1_1)(1_1 \tensor h') \\
          &=& hh' \\
          &=& (1_1 \tensor h)(h' \tensor 1_1) \\
          &=& (1_1 h') \tensor (h 1_1) \\
          &=& h' \tensor h
\end{array}
\]
for all $h,h' \maps 1 \to 1$.  Applied to 2-groups this implies
that the automorphisms of the object $1$ form an abelian group.
This is our abelian group $H$.

There is an action $\alpha$ of $G$ as automorphisms of $H$ given
by
\[   \alpha(g,h) = (1_g \tensor h) \tensor 1_{\overline{g}}  . \]
This is the same formula for $\alpha$ as in the crossed module
construction of Section \ref{internalizationsection}; we are just
writing it a bit differently now because a coherent 2-group is not
a category object in $\Grp$.  Here we need to be a bit more
careful to check that $\alpha$ is an action as automorphisms,
since the associator is nontrivial.

Finally, since our 2-group is skeletal, we do not need to
parenthesize tensor products of objects, and the associator gives
an automorphism
\[        a_{g_1,g_2,g_3} \maps g_1 \tensor g_2 \tensor g_3
                      \to  g_1 \tensor g_2 \tensor g_3.  \]
For any object $x \in G$ we identify $\Aut(x)$ with $\Aut(1) = H$
by tensoring with $\xb$ on the right: if $f \maps x \to x$ then $f
\tensor \xb \maps 1 \to 1$, since $x \tensor \xb = 1$.  By this
trick the associator can be thought of as a map from $G^3$ to $H$,
and by abuse of language we denote this map by:
\[
\begin{array}{rcl}
a \maps G^3 &\to& H \\
(g_1,g_2,g_3) &\mapsto & a(g_1,g_2,g_3) := a_{g_1,g_2,g_3} \tensor
\overline{g_1 \tensor g_2 \tensor g_3}.
\end{array}
\]
The pentagon identity implies that this map satisfies
\[   g_0 a(g_1,g_2,g_3) - a(g_0g_1, g_2, g_3) + a(g_0, g_1g_2, g_3) -
      a(g_0, g_1, g_2g_3) + a(g_0, g_1, g_2) = 0 \]
for all $g_0,g_1,g_2,g_3 \in G$, where the first term is defined
using the action of $G$ on $H$, and we take advantage of the
abelianness of $H$ to write its group operation as addition.  In
the language of group cohomology \cite{MacLane1}, this says precisely
that $a$ is a `3-cocycle' on $G$ with coefficients in the
$G$-module $H$.  Mac Lane's coherence theorem for monoidal categories
also implies that $a$ is a `normalized' 3-cocycle, meaning
that $a(g_1,g_2,g_3) = 1$ whenever $g_1,g_2$ or $g_3$ equals $1$.

This completes the construction of a quadruple $(G,H,\alpha,a)$
from any special 2-group.  Conversely, any such
quadruple determines a unique 2-group of this sort.  Since proving
this is largely a matter of running the previous construction
backwards, we leave this as an exercise for the reader.

Having shown that every coherent 2-group is equivalent
to one that can be described using group cohomology,
we now proceed to do the same thing for homomorphisms between
coherent 2-groups.

\begin{defn} \et  A {\bf special} homomorphism $F \maps C \to C'$
is a homomorphism between special 2-groups such that $F_0$ is an
identity morphism.
\end{defn}

\begin{prop} \label{isomorphism} \et  Any homomorphism between
special 2-groups is isomorphic in $\cg$ to a special homomorphism.
\end{prop}

\textbf{Proof. }  It suffices to show that for any weak monoidal
functor $F \maps C \to C'$ between weak monoidal categories, there is
a weak monoidal natural isomorphism $\theta \maps F \To F'$ where
$F'_0 \maps 1 \to F'(1)$ is an identity isomorphism.  We leave this as
an exercise for the reader.  \qed

To give a cohomological description of special homomorphisms,
let $F \maps C \to C'$ be a special homomorphism and
let $(G,H,\alpha,a)$ and $(C,H,\alpha',a')$ be the quadruples
corresponding to $C$ and $C'$, respectively.  The functor
$F$ maps objects to objects and preserves tensor products up
to isomorphism, so it gives a group homomorphism
\[         \phi \maps G \to G'. \]
For similar reasons, $F$ also gives a group homomorphism
\[         \psi \maps H \to H'  ,\]
and in fact this is a morphism of modules in the following sense:
\[      \psi(\alpha(g) h) = \alpha'(g) \psi(h) \]
for all $g \in G$ and $h \in H$.  As a weak monoidal functor, $F$
also comes equipped with an natural isomorphism
from $F(g_1) \tensor F(g_2)$ to $F(g_1 \tensor g_2)$ for all
$g_1,g_2 \in G$.  Since $C'$ is skeletal, this is an automorphism:
\[ (F_2)_{g_1,g_2} \maps F(g_1) \tensor F(g_2)
\to F(g_1) \tensor F(g_2). \]
Copying what we did for the associator, we define a map
\[
\begin{array}{rcl}
 k \maps G^2 &\to& H'   \\
          (g_1,g_2) &\mapsto & k(g_1,g_2) :=
(F_2)_{g_1,g_2} \tensor \overline{F(g_1) \tensor F(g_2)} .
\end{array}
\]
Using the fact that $F_0$ is the identity, the coherence
laws for the left and right unit laws in the definition of
a weak monoidal functor imply that $k(g_1,g_2) = 1$ whenever
$g_1$ or $g_2$ equals 1.  In the language of group
cohomology, $k$ is thus a `normalized 2-cocycle' on $G$ with
values in $H'$.  Furthermore,
the coherence law for the associator in the definition of
a weak monoidal functor implies that
\[    \psi(a(g_0,g_1,g_2)) - a'(\phi(g_0), \phi(g_1), \phi(g_2))  = \]
\[ \phi(g_0) k(g_1,g_2) - k(g_0g_1, g_2) + k(g_0, g_1 g_2) - k(g_1,g_2)
\]
for all $g_0,g_1,g_2 \in G$.  This says precisely that $\psi a$ and
$a' \phi^3$ differ by the coboundary of $k$:
\[        \psi a - a' \phi^3 = dk . \]

In short, a special homomorphism $F \maps C \to C$ gives a
triple $(\phi,\psi,k)$ where $\phi \maps G \to G'$ is a
group homomorphism, $\psi \maps H \to H'$ is a module
homomorphism, and $k$ is a normalized 2-cocycle on $G$ with
values in $H'$ such that $dk = \psi a - a' \phi^3$.  Conversely,
it is not hard to see that any such triple gives a
special homomorphism from $C$ to $C'$.

Finally, we give a cohomological description of 2-homomorphisms
between special homomorphisms.
Let $F,F' \maps C \to C'$ be special homomorphisms
with corresponding triples $(\phi,\psi,k)$ and $(\phi',\psi',k')$,
respectively.   A 2-homomorphism $\theta \maps F \To F'$
is just a monoidal natural transformation, so it gives a map
\[
\begin{array}{rcl}
 p \maps G &\to& H'   \\
         g &\mapsto & k(g) := \theta_g \tensor \overline{F(g)} .
\end{array}
\]
The condition that $\theta$ be natural turns out to have
no implications for $p$: it holds no matter what $p$ is.
However, the condition that $\theta$ be monoidal is equivalent
to the equations $p(1) = 1$ and
\[    k(g_1,g_2) - k'(g_1,g_2) =
\phi'(g_1) p(g_2) - p(g_1 g_2) + p(g_1) \]
for all $g_1,g_2 \in G$.  In the language of group cohomology,
these equations say that $p$ is a 1-cochain on $G$ with values in
$H'$ such that $dp = k - k'$.  So, 2-homomorphisms between special
homomorphisms are in one-to-one correspondence with 1-cochains of
this sort.

Summarizing all this, we obtain:

\begin{thm} \label{classification} \et  The 2-category
$\cg$ is biequivalent to the sub-2-category {\bf S2G} for which
the objects are special 2-groups,
the morphisms are special homomorphisms between these, and
the 2-morphisms are arbitrary 2-homomorphisms between those.
Moreover:
\begin{itemize}
\item There is a one-to-one correspondence between
special 2-groups $C$ and quadruples $(G,H,\alpha,a)$
consisting of:
\begin{itemize}
\item a group $G$,
\item an abelian group $H$,
\item an action $\alpha$ of $G$ as automorphisms of $H$,
\item a normalized 3-cocycle $a \maps G^3 \to H$.
\end{itemize}
\item Given special 2-groups $C,C'$ with corresponding
quadruples $(G,H,\alpha,a)$ and $(G',H',\alpha',a')$,
there is a one-to-one correspondence between special
homomorphisms $F \maps C \to C'$ and triples $(\phi,\psi,k)$
consisting of:
\begin{itemize}
\item a homomorphism of groups $\phi \maps G \to G'$,
\item a homomorphism of modules $\psi \maps H \to H'$,
\item a normalized 2-cochain $k \maps G^2 \to H'$
such that $dk = \psi a - a' \phi^3$.
\end{itemize}
\item Given special homomorphisms $F,F' \maps C \to C'$
with corresponding triples $(\phi,\psi,k)$ and $(\phi',\psi',k')$,
there is a one-to-one correspondence between 2-homomorphisms
$\theta \maps F \To F'$ and normalized 1-cochains $p \maps G \to H'$
with $dp = k - k'$.
\end{itemize}
\end{thm}

\textbf{Proof. }  The fact that $\cg$ is biequivalent
to the sub-2-category $\sg$ follows from the fact that
every object of $\cg$ is equivalent to an object in
$\sg$ (Proposition \ref{equivalence}) and every morphism of
$\cg$ is isomorphic to a morphism in $\sg$ (Proposition
\ref{isomorphism}).  The cohomological descriptions of objects,
morphisms and 2-morphisms in $\sg$ were deduced above.  \qed

We could easily use this theorem to give a complete description of
the 2-category $\sg$ in terms of group cohomology, but we
prefer to extract a simple corollary:

\begin{cor} \label{classification.cor} \et
There is a 1-1 correspondence between
equivalence classes of coherent 2-groups, where equivalence is as
objects of the 2-category $\cg$, and isomorphism classes of
quadruples $(G,H,\alpha,[a])$ consisting of:
\begin{itemize}
\item a group $G$,
\item an abelian group $H$,
\item an action $\alpha$ of $G$ as automorphisms of $H$,
\item an element $[a]$ of the cohomology group $H^3(G,H)$,
\end{itemize}
where an isomorphism from $(G,H,\alpha,[a])$ to
$(G',H',\alpha',[a'])$ consists of an isomorphism from $G$ to $G'$
and an isomorphism from $H$ to $H'$, carrying $\alpha$ to
$\alpha'$ and $[a]$ to $[a']$.
\end{cor}

\textbf{Proof. }  This follows directly from Theorem
\ref{classification}, together with the fact that group
cohomology can be computed using normalized cochains.
\qed

Though the main use of Proposition \ref{equivalence} is to help
prove Theorem \ref{classification}, it has some interest in its
own right, because it clarifies the extent to which any coherent
2-group can be made simultaneously both skeletal and strict.  Any
coherent 2-group is equivalent to a skeletal one in which $\ell,r,i$
and $e$ are identity natural transformations --- but not the associator,
unless the invariant $[a] \in H^3(G,H)$ vanishes.
On the other hand, if we drop our insistence on making a 2-group
skeletal, we can make it completely strict:

\begin{prop} \label{strictification} \et  Every coherent 2-group is
equivalent in $\cg$ to a strict one --- that is,
one for which $\ell,r,i,e$ and $a$ are identity natural
transformations.
\end{prop}

\textbf{Proof. }  Let $C$ be a coherent 2-group.  By a theorem
of Mac Lane \cite{MacLane2}, there is a strict monoidal category $C'$
that is equivalent to $C'$ as a monoidal category.  We can use this
equivalence to transfer the coherent 2-group structure from $C$ to
$C'$, making $C'$ into a coherent 2-group for which $\ell,r,$ and $a$ are
identity natural transformations, but not yet $i$ and $e$.

As a strict monoidal category, $C'$ is an object of $\Cat{\rm Mon}$, the
category of `monoids in $\Cat$'.  There is a pair of adjoint functors
consisting of the forgetful functor $U \maps \Cat\Grp \to \Cat{\rm Mon}$
and its left adjoint $F \maps \Cat{\rm Mon} \to \Cat\Grp$.   Thus
$C'' = F(C')$ is a group in $\Cat$, or in other words a strict 2-group.
It suffices to show that $C'$ is equivalent to $C''$ as an object
of $\cg$.

The unit of the adjunction between $U$ and $F$ gives a strict monoidal
functor $i_{C'} \maps C' \to U(F(C'))$, which by Theorem \ref{preservation}
determines a 2-group homomorphism from $C'$ to $C'' = F(C')$.
One can check that this is extends to an equivalence in $\cg$; we leave
this to the reader.

An alternative approach uses Proposition \ref{equivalence} to note
that $C$ is equivalent to a special 2-group $C'$.  From
the quadruple $(G,H,\alpha,a)$ corresponding to this special
2-group one can construct a crossed module (see
Mac Lane \cite{MacLane1} or, for a more readable treatment,
Ken Brown's text on group cohomology \cite{KBrown}).  This crossed
module in turn gives a strict 2-group $C''$, and one can check that
$C''$ is equivalent to $C$ in $\cg$.  The details for this approach
can be found in the 1986 draft of Joyal and Street's paper on braided
tensor categories \cite{JS}.  \qed

This result explains why Mac Lane and Whitehead \cite{MW} were able to
use strict 2-groups (or actually crossed modules) to describe arbitrary
connected pointed homotopy 2-types, instead of needing the more general
coherent 2-groups.

\subsection{Strict Lie 2-groups}

It appears that just as Lie groups describe continuous symmetries
in geometry, Lie 2-groups describe continuous symmetries in
categorified geometry.  In Definition \ref{lie.2-group} we said that
Lie 2-groups are coherent 2-groups in $\Diff\Cat$, the 2-category of
smooth categories.  In this section we shall give some examples,
but only `strict' ones, for which the associator, left and right
unit laws, unit and counit are all identity 2-morphisms.  We discuss
the challenge of finding interesting nonstrict Lie 2-groups
in the next section.

Strict Lie 2-groups make no use of the 2-morphisms in $\Diff\Cat$, so
they are really just groups in the underlying category of $\Diff\Cat$.
By `commutativity of internalization', these are the same as categories
in $\Diff\Grp$, the category of Lie groups.
To see this, note that if $C$ is a strict Lie 2-group, it is first
of all an object in $\Diff\Cat$.  This means it is a category with
a manifold of objects $C_0$ and a manifold of morphisms $C_1$,
with its source, target, identity-assigning and composition maps
all smooth. But since $C$ is a {\it group} in $\Diff\Cat$, $C_0$
and $C_1$ become Lie groups, and all these maps become Lie group
homomorphisms. Thus, $C$ is a category in $\Diff\Grp$.  The
converse can be shown by simply reversing this argument.

Treating strict Lie 2-groups as categories in $\Diff\Grp$ leads
naturally to yet another approach, where we treat them as `Lie
crossed modules'.  Here we use the concept of `crossed module in
$K$', as described in Definition \ref{internal.crossed.module}:

\begin{defn} \et
A {\bf Lie crossed module} is a crossed module in $\Diff$.
\end{defn}

\noindent Concretely, a Lie crossed module is a quadruple
$(G,H,t,\alpha)$ consisting of Lie groups $G$ and $H$, a
homomorphism $t \maps H \to G$, and an action $\alpha$ of $G$ on
$H$ such that $t$ is $G$-equivariant
\[   t(\alpha(g,h) = g \, t(h)\, g^{-1} \]
and $t$ satisfies the Peiffer identity
\[    \alpha(t(h),h') = hh'h^{-1}  \]
for all $g \in G$ and $h,h' \in H$.   Proposition
\ref{crossed.module} shows how we can get a Lie crossed module
from a strict Lie 2-group and vice versa.  Using this, one can
construct a 2-category of strict Lie 2-groups and a 2-category of
Lie crossed modules and show that they are equivalent.  This
equivalence lets us efficiently construct many examples of strict
Lie 2-groups:

\begin{example} \label{crossed.1} \et {\rm
Given any Lie group $G$, abelian Lie group $H$, and homomorphism
$\alpha \maps G \to \Aut(H)$, there is a Lie crossed module with
$t \maps G \to H$ the trivial homomorphism and $G$ acting on $H$
via $\rho$. Because $t$ is trivial, the corresponding strict Lie
2-group $C$ is `skeletal', meaning that any two isomorphic objects
are equal. It is easy to see that conversely, all skeletal strict
Lie 2-groups are of this form. }\end{example}

\begin{example} \label{crossed.2} \et {\rm
Given any Lie group $G$, we can form a Lie crossed module as in
Example \ref{crossed.1} by taking $H = \g$, thought of as an
abelian Lie group, and letting $\alpha$ be the adjoint
representation of $G$ on $\g$.  If $C$ is the corresponding strict
Lie 2-group we have
\[    C_1 = \g \rtimes G \iso TG \]
where $TG$ is the tangent bundle of $G$, which becomes a Lie group
with product
\[    dm \maps TG \times TG \to TG , \]
obtained by differentiating the product
\[    m \maps G \times G \to G .\]
We call $C$ the {\bf tangent 2-group} of $G$ and denote it as $\T
G$.

Another route to the tangent 2-group is as follows.  Given any
smooth manifold $M$ there is a smooth category $\T M$, the {\bf
tangent groupoid} of $M$, whose manifold of objects is $M$ and
whose manifold of morphisms is $TM$.  The source and target maps
$s,t \maps TM \to M$ are both the projection to the base space,
the identity-assigning map $i \maps M \to TM$ is the zero section,
and composition of morphisms is addition of tangent vectors.  In
this category the arrows are actually little arrows --- that is,
tangent vectors!

This construction extends to a functor
\[    \T \maps \Diff \to \Diff\Cat  \]
in an obvious way.  This functor preserves products, so it sends
group objects to group objects.  Thus, if $G$ is a Lie group, its
tangent groupoid $\T G$ is a strict Lie 2-group. }\end{example}

\begin{example} \label{crossed.3} \et {\rm
Similarly, given any Lie group $G$, we can form a Lie crossed
module as in Example \ref{crossed.1} by letting $\alpha$ be the
coadjoint representation on $H = \g^*$. If $C$ is the
corresponding Lie 2-group, we have
\[    C_1 = \g^* \rtimes G \iso T^*G \]
where $T^*G$ is the cotangent bundle of $G$.  We call $C$ the {\bf
cotangent 2-group} of $G$ and denote it as $\T^* G$.
}\end{example}

\begin{example} \label{crossed.4} \et {\rm
More generally, given any representation $\alpha$ of a Lie group
$G$ on a finite-dimensional vector space $V$, we can form a Lie
crossed module and thus a strict Lie 2-group with this data,
taking $H = V$. For example, if $G$ is the Lorentz group
$\SO(n,1)$, we can form a Lie crossed module by letting $\alpha$
be the defining representation of $\SO(n,1)$ on $H = \R^{n+1}$. If
$C$ is the corresponding strict Lie 2-group, we have
\[    C_1 = \R^{n+1} \rtimes \SO(n,1) \iso \ISO(n,1) \]
where $\ISO(n,1)$ is the Poincar\'e group.  We call $C$ the {\bf
Poincar\'e 2-group}.   After this example was introduced by one of
the authors \cite{Baez}, it became the basis of an interesting new
theory of quantum gravity \cite{CS,CY}. }\end{example}

\begin{example} \label{crossed.5} \et {\rm
Given any Lie group $H$, there is a Lie crossed module with $G =
\Aut(H)$, $t \maps H \to G$ the homomorphism assigning to each
element of $H$ the corresponding inner automorphism, and the
obvious action of $G$ as automorphisms of $H$.  We call the
corresponding strict Lie 2-group the {\bf strict automorphism
2-group} of $H$, $\Aut_s(H)$, because its underlying 2-group is
just $\Aut_s(H)$ as defined previously.  }\end{example}

\begin{example} \label{crossed.6} \et {\rm
If we take $H = \SU(2)$ and form $\Aut_s(H)$, we get a strict Lie
2-group with $G = \SO(3)$.   Similarly, if we take $H$ to be the
multiplicative group of nonzero quaternions, $\Aut_s(H)$ is again
a strict Lie 2-group with $G = \SO(3)$. This latter example is
implicit in Thompson's work on `quaternionic gerbes'
\cite{Thompson}.  }\end{example}

\begin{example} \label{crossed.7} \et {\rm
Suppose that $t \maps G \to H$ is a surjective homomorphism of Lie
groups.  Then there exists a Lie crossed module $(G,H,t,\alpha)$
if and only if $t$ is a central extension (that is, the kernel of
$t$ is contained in the center of $G$). Moreover, when this Lie
crossed module exists it is unique. }\end{example}

\begin{example} \label{crossed.8} \et {\rm
Suppose that $V$ is a finite-dimensional real vector space
equipped with an antisymmetric bilinear form $\omega \maps V
\times V \to \R$.  Make $H = V \oplus \R$ into a Lie group with
the product
\[    (v,\alpha)(w,\beta) = (v+w, \alpha + \beta + \omega(v,w)) .\]
This Lie group is called the `Heisenberg group'. Let $G$ be $V$
thought of as a Lie group, and let $t \maps H \to G$ be the
surjective homomorphism given by
\[    t(v,\alpha) = v .\]
Then $t$ is a central extension, so by Example \ref{crossed.7} we
obtain a 2-group which we call the {\bf Heisenberg 2-group} of
$(V,\omega)$. }\end{example}

\subsection{2-Groups from Chern--Simons theory}
\label{chernsimonssection}

We conclude by presenting some interesting examples of 2-groups
built using Chern--Simons theory.  Since the existence of these
2-groups was first predicted using an analogy between the
classifications of 2-groups and Lie 2-algebras, we begin by sketching
this analogy.  We then describe some nonstrict Lie
2-algebras discussed in the companion paper HDA6, and use
this analogy together with some results from Chern--Simons theory
to build corresponding 2-groups.  Naively, one would expect these to
be Lie 2-groups.  However, we prove a `no-go theorem' ruling out
the simplest ways in which this could be true.

The paper HDA6 studies `semistrict Lie 2-algebras'.  These are
categorified Lie algebras in which the Jacobi identity
has been weakened, but not the antisymmetry of the bracket.
A bit more precisely, a semistrict Lie 2-algebra is a category in
$\Vect$, say $L$, equipped with an antisymmetric bilinear functor
called the `bracket':
\[       [\cdot, \cdot] \maps L \times L \to L  , \]
together with a natural isomorphism called the `Jacobiator':
\[        J_{x,y,z} \maps [[x,y],z] \to [x,[y,z]] + [[x,z],y]  \]
satisfying certain coherence laws of its own.

HDA6 gives a classification of semistrict Lie 2-algebras that
perfectly mirrors the classification of 2-groups summarized
in Corollary \ref{classification.cor} above, but with
Lie algebras everywhere replacing groups.  Namely, there is
a 1-1 correspondence between equivalence classes of semistrict Lie
2-algebras $L$ and isomorphism classes of quadruples $(\g,\h,\rho,[j])$
consisting of:
\begin{itemize}
\item a Lie algebra $\g$,
\item an abelian Lie algebra $\h$,
\item a representation $\rho$ of $\g$ as derivations of $\h$,
\item an element $[j]$ of the Lie algebra cohomology group $H^3(\g,\h)$.
\end{itemize}
Here $\g$ is the Lie algebra of objects in a skeletal version of
$L$, $\h$ is the Lie algebra of endomorphisms of the zero object
of $L$, the representation $\rho$ comes from the bracket in $L$, and
the 3-cocycle $j$ comes from the Jacobiator.  Of course, an abelian Lie
algebra is nothing but a vector space, so it adds nothing to say
that in the representation $\rho$ elements of $\g$ act `as derivations'
of $\h$.  We say this merely to make the analogy to Corollary
\ref{classification.cor} as vivid as possible.

Recall that in the classification of 2-groups, the cohomology class
$[a] \in H^3(G,H)$ comes from the associator in a skeletal version
of the 2-group in question.   In fact, this class is the only obstruction
to finding an equivalent 2-group that is both skeletal and strict.  The
situation for Lie 2-algebras is analogous: the cohomology class
$[j] \in H^3(\g,\h)$ comes from the Jacobiator, and gives the
obstruction to finding an equivalent Lie 2-algebra that is both skeletal
and strict.

Using this, in HDA6 we construct some Lie 2-algebras that are not
equivalent to skeletal strict ones.  Suppose $G$ is a connected and
simply-connected compact simple Lie group, and let $\g$ be its Lie
algebra.  Let $\rho$ be the trivial representation of $\g$ on $\u(1)$,
the 1-dimensional abelian Lie algebra over the reals.  Then
\[           H^3(\g,\u(1)) \iso \R  .\]
By the classification of Lie 2-algebras, for any value of $\hbar \in \R$
we obtain a skeletal Lie 2-algebra $\g_\hbar$ having $\g$ as its Lie
algebra of objects and $\u(1)$ as the endomorphisms of its zero object.
When $\hbar = 0$ this Lie 2-algebra is just $\g$ with identity morphisms
adjoined to make it into a strict Lie 2-algebra.  However, when $\hbar \ne 0$,
this Lie 2-algebra is not equivalent to a skeletal strict one.

An interesting question is whether these Lie 2-algebras have corresponding
Lie 2-groups.  There is not a general construction of Lie 2-groups
from Lie 2-algebras, but we can try to build them `by hand'.  We
begin by seeking a skeletal 2-group $G_\hbar$ with $G$ as its group of
objects and $\U(1)$ as the automorphism group of its identity object,
which is strict only at $\hbar = 0$.  To define the associator in
$G_\hbar$, we would like to somehow `exponentiate' the element
of $H^3(\g,\u(1))$ coming from the Jacobiator in $\g_\hbar$
to obtain an element of $H^3(G,\U(1))$.  However,
from experience with affine Lie algebras and central extensions of
loop groups, we expect this to be possible only for elements of
$H^3(\g,\u(1))$ satisfying some sort of integrality condition.

Indeed this is the case: sitting inside the Lie algebra cohomology
$H^3(\g,\u(1)) \iso \R$ there is a lattice $\Lambda$, which we can
identify with $\Z$, that comes equipped with an inclusion
\[      \iota \maps \Lambda \hookrightarrow H^3(G,\U(1))  .\]
This is actually a key result from the papers of Chern--Simons
\cite{ChernS} and Cheeger--Simons \cite{CheegerS} on secondary
characteristic classes.  We describe how this inclusion is constructed
below, but for now we record this:

\begin{thm} \label{Ghbar} \et
Let $G$ be a connected and simply-connected compact simple Lie group.
Then for any $\hbar \in \Z$ there exists a special 2-group $G_\hbar$
having $G$ as its group of objects, $\U(1)$ as the group
of endomorphisms of its unit object, the trivial action of
$G$ on $\U(1)$, and $[a] \in H^3(G,\U(1))$ given by
$\iota(\hbar)$.  The 2-groups $G_\hbar$ are inequivalent for
different values of $\hbar$, and strict only for $\hbar = 0$.
\end{thm}

To give more of a feeling for this result, let us sketch how
the lattice $\Lambda$ and the map $\iota$ can be constructed.
Perhaps the most illuminating approach
uses this commutative diagram:
\[
 \xymatrix{
 & & H^{2n-1}(G,\U(1)) \\
 & H^{2n}_{\rm top}(BG,\R) \ar[d]_{\tau_\R}
 & H^{2n}_{\rm top}(BG,\Z) \ar[u]_{\kappa} \ar[l]_{\iota_{BG}}
\ar[d]^{\tau_\Z}       \\
 H^{2n-1}(\g,\u(1))
 & H^{2n-1}_{\rm top}(G,\R) \ar[l]_{\quad \sim}
 & H^{2n-1}_{\rm top}(G,\Z) \ar[l]_{\iota_G}
}
\]
In this diagram, the subscript `top' refers to the cohomology of
the compact simple Lie group $G$ or its classifying space $BG$ as a
topological space.  The integral cohomology
$H^{2n-1}_{\rm top}(G,\Z)$ maps to a lattice in
the vector space $H^{2n-1}_{\rm top}(G,\R)$, and thus defines
a lattice $\Lambda$ inside the isomorphic vector space
$H^{2n-1}(\g,\u(1))$.  In the case relevant here, namely $n = 2$,
the maps labelled $\tau$ are isomorphisms and the maps labelled
$\iota$ and $\kappa$ are injections.  Thus, in this case the
diagram serves to define an injection
\[      \iota \maps \Lambda \hookrightarrow H^3(G,\U(1))  .\]

Let us say a few words about the maps in this diagram.
The isomorphism from $H^n_{\rm top}(G,\R)$ to $H^n(\g,\u(1))$
is defined using deRham theory: there is a
cochain map given by averaging differential forms on $G$
to obtain left-invariant forms, which can be identified with
cochains in Lie algebra cohomology \cite{CE}.
Since the classifying space has $\Omega(BG) \simeq G$, there are
`transgression' maps $\tau_\Z$ and $\tau_\R$ from the $2n$th
integral or real cohomology $BG$ to the $(2n-1)$st cohomology of
$G$.  These are isomorphisms for $n = 2$, since in general
the transgression map $\tau \maps H^{q+1}_{\rm top}(X,R) \to
H^q_{\rm top}(\Omega X,R)$ is an isomorphism whenever $X$ is
$k$-connected, $q \le 2k-1$ and the coefficient ring $R$ is a
principal ideal domain \cite{Whitehead0}.  Finally,
the change--of--coefficient maps $\iota_{BG}$ and $\iota_G$
map the integral cohomology of either of these spaces to a
full lattice in its real cohomology.  The universal coefficient theorem
implies $\iota_G$ is an injection for $n = 2$ because the 3rd integral
cohomology of a compact simple Lie group is torsion-free, in fact $\Z$.
Similarly, $\iota_{BG}$ is an injection because $H^4_{\rm top}(BG,\Z)
\iso \Z$.

The innovation of Chern, Cheeger and Simons was
the homomorphism $\kappa$, which maps elements of
$H^{2n}_{\rm top}(BG,\Z)$ to certain elements of
$H^{2n-1}(G,\U(1))$ called `secondary characteristic classes'.
This is where some differential geometry enters the story.
For ease of exposition, we describe this map only in the case
we need, namely $n = 2$.  In this particular case we only
need to say what $\kappa$ does to the standard generator of
$H^4_{\rm top}(BG,\Z)$, which is called the `second Chern class' $c_2$.

Since $BG$ is the classifying space for principal $G$-bundles,
any principal $G$-bundle $P$ over a smooth manifold $M$
gives a homotopy class of maps $M \to BG$, which we can use to
pull back $c_2$ to an element of $H^4_{\rm top}(M,\Z)$.  Chern
showed that the corresponding element of $H^4_{\rm top}(M,\R)$
can be described using deRham theory by choosing an arbitrary
connection $A$ on $P$.  We can think of this connection as
$\g$-valued 1-form on $P$, and its curvature
\[     F = dA + A \wedge A  \]
as an $\g$-valued 2-form.  This allows us to define a 4-form on
$P$:
\[      c_2(A) = {1\over 8\pi^2} \tr(F \wedge F) . \]
where `$\tr$' is defined using a suitably normalized invariant
bilinear form on $\g$.  The 4-form $c_2(A)$ is the pullback of a
unique closed 4-form on $M$, which represents the image of
$c_2(P)$ in $H^4_{\rm top}(M,\R)$.

While the 4-form down on $M$ is merely closed, Chern and Simons noted
that $c_2(A)$ itself is actually exact, being the differential of
this 3-form:
\[   {\rm CS}_2(A) = {1 \over 8\pi^2} \tr(A \wedge dA
                     + {2\over 3} A \wedge A \wedge A)  .\]
If the connection $A$ is flat, meaning $F = 0$, then
${\rm CS}_2(A)$ is closed.  It thus represents an element of
$H^3_{\rm top}(P,\R)$.  This element is not canonically
the pullback of an element of $H^3_{\rm top}(M,\R)$, but
it is up to an integral cohomology class.

It follows that ${\rm CS}_2(A)$ canonically gives rise to an
element of $H^3_{\rm top}(M,\R/\Z) \iso H^3_{\rm top}(M,\U(1))$
for any principal $G$-bundle with flat connection over $M$.
Note however that a principal $G$-bundle with
flat connection is the same as a principal $G_\delta$-bundle, where
$G_\delta$ is the group $G$ equipped with the discrete topology.
Since our assignment of cohomology classes to manifolds equipped
with principal $G_\delta$-bundle is functorial, it must be a
characteristic class: in other words, it must come from pulling back
some element of $H^3_{\rm top}(BG_\delta,\U(1))$ along the classifying
map $M \to BG_\delta$.  But $H^3_{\rm top}(BG_\delta,\U(1))$ is just
another way of talking about the group cohomology $H^3(G,\U(1))$.
Thus we obtain an element ${\rm CS}_2 \in H^3(G,\U(1))$.

Since the second Chern class generates $H^4_{\rm top}(BG,\Z)$, we
can define
\[    \kappa \maps H^4_{\rm top}(BG, \Z) \to H^3(G,\U(1))   \]
by
\[      \kappa(c_2) = {\rm CS}_2  .\]
One can show that $\kappa$ is an injection by explicit calculations
\cite{Jeffrey}.

It would be natural to hope the 2-groups $G_\hbar$ are Lie 2-groups
and therefore topological 2-groups.  However, we shall conclude with a
`no-go theorem' saying that $G_\hbar$ can be made into a
topological 2-group with a reasonable topology only in the trivial case
$\hbar = 0$.  For this, we start by internalizing the cohomological
classification of special 2-groups given in Theorem \ref{classification}.
Suppose $K$ is any category with finite products such that
$K\Grp$ has finite limits.  We discussed the concept of `coherent
2-group in $K\Cat$' in Section \ref{internalizationsection}.
We now say what it means for such a 2-group to be `special':

\begin{defn} \et A {\bf special 2-group} $C$ in $K\Cat$ is a
coherent 2-group in $K\Cat$ for which:
\begin{enumerate}
\item its underlying category in $K$ is {\bf skeletal}, meaning
that the source and target morphisms $s,t \maps C_1 \to C_0$ are
equal,
\item the equalizer of the morphisms $s \maps C_1 \to C_0$
and
$ \xymatrix@1{C_1 \ar[r] & I \ar[r]^{i} & C_0 }$
exists,
\item the left unit law $\ell$, the right unit law $r$, the unit
$i$ and the counit $e$ are identity natural transformations.
\end{enumerate}
\end{defn}

Given a special 2-group $C$ in $K\Cat$, we can obtain a
quadruple $(G,H,\alpha,a)$ by internalizing the construction
described in Section \ref{classificationsection}.  We merely
sketch how this works.  The multiplication in $C$ makes
$C_0$ into a group in $K$, even if the associator is nontrivial,
since $C$ is skeletal.  Let $G$ be this group in $K$.
Composition of morphisms makes the equalizer in item 2 into
an abelian group in $K$, thanks to the Eckmann--Hilton argument.
Let $H$ be this abelian group in $K$.  Conjugation in $C$ gives an
action $\alpha$ of $G$ as automorphisms of $H$, and the
associator of $C$ gives a morphism $a \maps G^3 \to H$.  This
morphism $a$ is a normalized 3-cocycle in the cochain complex
for {\bf internal group cohomology}:
\[
 \xymatrix{\hom(G^0,H)
\ar[r]^d &  \hom(G^1,H) \ar[r]^d & \hom(G^2,H) \ar[r]^d &
 \cdots
\hbox{\qquad \qquad}}
\]
where the differential is defined as usual for group cohomology.
It thus defines an element $[a] \in H^3(G,H)$ of internal group
cohomology.  Conversely, given a quadruple $(G,H,\alpha,a)$ of
this form, we can obtain a special 2-group in $K\Cat$.

We have been unable to show that every coherent 2-group
in $K\Cat$ is equivalent to a special one, or even a skeletal one.
After all, to show this for $K = \Set$, we used the axiom of choice
to pick a representative for each isomorphism class of objects in a
given 2-group $C$.  This axiom is special to $\Set$, and
fails in many other categories.  So, the above cohomological
description of special 2-groups in $K\Cat$ may not yield a
complete classification of coherent 2-groups in $K\Cat$.
Nonetheless we can use it to obtain some information about
the problem of making the 2-groups $G_\hbar$ into topological or
Lie 2-groups.

To do this, we also need the concept of `special homomorphisms'
between special 2-groups in $K\Cat$:

\begin{defn} \et  A {\bf special} homomorphism $F \maps C \to C'$
is a homomorphism between special 2-groups such that $F_0$ is an
identity morphism.
\end{defn}

\noindent
Recall that $K\Cat\cg$ is the 2-category of coherent 2-groups
in $K\Cat$.   By a straightforward internalization of Theorem
\ref{classification} we obtain:

\begin{prop} \label{internalclassification} \et
Suppose that $K$ is a category with finite products.  The
2-category ${K\Cat}\cg$ has a sub-2-category
\textbf{\textit{K}{\bf CatS2G}}
for which the objects are special 2-groups,
the morphisms are special homomorphisms between these, and
the 2-morphisms are arbitrary 2-homomorphisms between those.
There is a 1-1 correspondence between
equivalence classes of objects in $K\Cat{\rm S2G}$ and
isomorphism classes of quadruples $(G,H,\alpha,[a])$ consisting of:
\begin{itemize}
\item a group $G$ in $K$,
\item an abelian group $H$,
\item an action $\alpha$ of $G$ as automorphisms of $H$,
\item an element $[a] \in H^3(G,H)$.
\end{itemize}
\end{prop}

Now we consider $K = \Top$.  In this case the internal group
cohomology is usually called `continuous cohomology', and we
shall denote it by $H^n_{\rm cont}(G,H)$ to avoid confusion.

\begin{thm} \label{strict} \et
Let $G$ be a connected compact Lie group and $H$ a connected
abelian Lie group.  Suppose $C$ is
a special topological 2-group having $G$ as its group of objects
and $H$ as the group of endomorphisms of its unit object.  Then
the associator $a$ of $C$ has $[a] = 0$.  Thus $C$ is equivalent
in $\Top\Cat{\rm S2G}$ to a special topological 2-group that is
strict.
\end{thm}

\textbf{Proof. }  The work of Hu \cite{Hu},
van Est \cite{vanEst} and Mostow \cite{Mostow}
on continuous cohomology implies that
$H^3_{\rm cont}(G,H)$ is trivial.  We thus have $[a] = 0$,
and the rest follows from Proposition
\ref{internalclassification}.

For the sake of completeness we sketch the proof that
$H^3_{\rm cont}(G,H) \iso \{0\}$.  First we consider the case
where $H$ is a real vector space equipped with an arbitrary
representation of $G$.  For any
continuous cocycle $f \maps G^n \to H$ with $n \ge 1$ there
is a continuous cochain $F \maps G^{n-1} \to H$ given by
\[  F(g_1,\dots,g_{n-1}) = \int_G f(g_1, \dots g_n) dg_n ,\]
where the integral is done using the normalized Haar measure
on $G$.  A simple calculation shows that $dF = \pm f$.
This implies that $H^n_{\rm cont}(G,H) \iso \{0\}$ for all
$n \ge 1$.

In general, any action of $G$ on a connected abelian Lie
group $H$ lifts uniquely to an action on the universal cover
$\tilde H$, which is a real vector space.  Any normalized continuous
cochain $f \maps G^n \to H$ lifts uniquely to a normalized continuous
cochain $\tilde f \maps G^n \to \tilde H$ for $n \ge 2$, since the
$n$-fold smash product of $G$ with itself is simply-connected in this
case.  Since $d\tilde f = \widetilde{df}$, this implies that
$H^n_{\rm cont}(G,H) \iso H^n_{\rm cont}(G,\tilde H) \iso \{0\}$
for $n \ge 3$.
\qed

Now suppose $C$ is a topological 2-group whose underlying
2-group is isomorphic to a 2-group of the form $G_\hbar$ for
some $\hbar \in \Z$.  Then the objects of $C$ form a topological
group which is isomorphic as a group to $G$, but possibly with some
nonstandard topology, e.g.\ the discrete topology.  Similarly,
the endomorphisms of its identity object form a topological group
which is isomorphic as a group to $\U(1)$, but possibly with some
nonstandard topology.

\begin{cor} \et Let $G$ be a connected and simply-connected compact
simple Lie group.
Suppose $C$ is a topological 2-group whose underlying 2-group
is isomorphic to $G_\hbar$ for some $\hbar \in \Z$.  If the
topological group of objects of $C$ is isomorphic to $G$ with
its usual topology, and the topological group of endomorphisms
of its identity object is isomorphic to $\U(1)$ with its usual
topology, then $\hbar = 0$.
\end{cor}

\textbf{Proof. }  Given the assumptions, $C$ is a special
topological 2-group which fulfills the hypotheses of Theorem
\ref{strict}.  It is thus equivalent in $\Top\Cat{\rm S2G}$
to a strict special topological 2-group, so its underlying
2-group $G_\hbar$ is equivalent in ${\rm C2G}$ to a strict skeletal
2-group.  By Theorem \ref{Ghbar} this happens only for $\hbar = 0$.
\hbox{\hskip 30 em} \qed

In rough terms, this means that for $\hbar \ne 0$,
the 2-group $G_\hbar$ cannot be made into a topological 2-group
with a sensible topology.  However, we have not ruled out the
possibility that it is {\it equivalent} to the underlying 2-group
of some interesting topological 2-group, or even of some Lie 2-group.
Another possibility is that the concept of Lie 2-group needs to be
broadened to handle this case --- perhaps along lines suggested
by Brylinksi's paper on multiplicative gerbes \cite{Brylinski,CJMSW}.

\subsection*{Acknowledgements}

We thank Lawrence Breen, James Dolan, Danny Stevenson
and Ross Street for discussions
on this subject, Miguel Carri\'on-\'Alvarez for suggesting Lemma
\ref{imove}, Alissa Crans for careful editing and discussions about
the relation of Lie 2-groups to Lie 2-algebras, and Andr\'ee Ehresmann
for correspondence on the history of internal categories.

\end{document}